\documentclass{amsart}

\usepackage{amscd,amssymb}
\usepackage[mathcal]{euscript}
\usepackage[matrix,arrow,curve,cmtip]{xy}

\newtheorem{theorem}{Theorem}[subsection]
\newtheorem{lemma}[theorem]{Lemma}
\newtheorem{proposition}[theorem]{Proposition}

\newtheorem{corollary}[theorem]{Corollary}
\theoremstyle{definition}
\newtheorem{definition}[theorem]{Definition}
\newtheorem{remark}[theorem]{Remark}
\newtheorem{example}[theorem]{Example}
\newtheorem{construction}[theorem]{Construction}

\numberwithin{equation}{subsection}

\DeclareMathOperator{\ho}{Ho}
\DeclareMathOperator{\Hom}{Hom}
\DeclareMathOperator{\Map}{Map}
\DeclareMathOperator{\shom}{Map}
\DeclareMathOperator{\ihom}{Hom}
\DeclareMathOperator{\Ev}{Ev}
\DeclareMathOperator{\sym}{Sym}

\newcommand{\wt}[1]{\widetilde{#1}}
\newcommand{\RK}{R^{\infty} K}

\newcommand{\boxprod}{\mathbin\square}
\newcommand{\cN}{\mathbb{N}}
\newcommand{\cat}[1]{\mathcal{#1}}
\newcommand{\cc}{\cat{C}}
\newcommand{\cd}{\cat{D}}
\newcommand{\cm}{\cat{M}}
\newcommand{\cp}{\cat{P}}
\newcommand{\cof}{\text{-cof}}
\newcommand{\fib}{\text{-fib}} 
\newcommand{\cs}{\cat{S}}
\newcommand{\from}{\leftarrow}
\newcommand{\inj}{\text{-inj}}
\newcommand{\llp}{left lifting property with respect to }
\newcommand{\natiso}{\cong}
\newcommand{\ov}{\overline}
\newcommand{\proj}{\text{-proj}}
\newcommand{\rlp}{right lifting property with respect to }
\newcommand{\sep}{0.5ex}
\newcommand{\set}{\text{Set}}
\newcommand{\sm}{\wedge}

\newcommand{\BF}{Sp^{\cN}}
\newcommand{\colim}{\operatorname{colim}}
\newcommand{\seq}{\sset^{\cN}}
\newcommand{\spec}{Sp^{\Sigma }}
\newcommand{\sset}{\mathcal{S}_{*}}
\newcommand{\symseq}{\sset^{\Sigma }}
\renewcommand{\top}{\mathrm{Top}_{*}}
\newcommand{\spaces}{\mathcal{T}_{*}}
\newcommand{\cg}{\mathcal{K}_{*}}
\newcommand{\topspec}{\spec _{\spaces }}
\newcommand{\topsymseq}{\spaces^{\Sigma }}
\newcommand{\topBF}{\BF _{\spaces}}
\newcommand{\georeal}[1]{| #1 |}
\newcommand{\Sing}{\operatorname{Sing}}

\newcommand{\topmonoids}{\mathcal{M}(\topspec )}
\newcommand{\ie}{\textit{i.e., }}
\newcommand{\mathcolon}{\colon\,}
\newcommand{\badmap}{\lambda}
\begin{document}

\title{Symmetric spectra}

\date{\today} 
\author{Mark Hovey} 
\thanks{The first two authors were partially supported by NSF
Postdoctoral Fellowships}  
\address{Department of Mathematics \\ 
Wesleyan University \\ 
Middletown, CT}
\email{hovey@member.ams.org} 
\author{Brooke Shipley}
\address{Department of Mathematics \\
University of Chicago \\ 
Chicago, IL}
\email{bshipley@math.uchicago.edu} 
\author{Jeff Smith}
\thanks{The third author was partially supported by an NSF Grant}
\address{Department of Mathematics \\ 
Purdue University \\ 
West Lafayette, IN} 
\email{jhs@math.purdue.edu}

\maketitle

\tableofcontents

\section*{Introduction} 

Stable homotopy theory studies spectra as the linear approximation to
spaces.  Here, ``stable" refers to the consideration of spaces after
inverting the suspension functor.  This approach is a general one: one
can often create a simpler category by inverting an operation such as
suspension.  In this paper we study a particularly simple model for
inverting such operations which preserves product structures.  The
combinatorial nature of this model means that it is easily transported,
and hence may be useful in extending the methods of stable homotopy
theory to other settings.

The idea of a spectrum is a relatively simple one:
Freudenthal's suspension theorem implies that the sequence of homotopy
classes of maps
\[
[X,Y] \xrightarrow{} [\Sigma X,\Sigma Y] \xrightarrow{} \dots
\xrightarrow{} [\Sigma ^{n}X,\Sigma ^{n}Y]  \xrightarrow{} \dots 
\]
is eventually constant for finite-dimensional pointed CW-complexes $X$
and $Y$, where $\Sigma X=S^{1}\sm X$ is the reduced suspension of $X$.
This suggests forming a stable category where the suspension functor is
an isomorphism.  The standard way to do this is to define a
\emph{spectrum} to be a sequence of pointed topological spaces (or
simplicial sets) $X_{n}$ together with structure maps $S^{1}\sm
X_{n}\xrightarrow{}X_{n+1}$.  This was first done by Lima~\cite{lima}
and later generalized by Whitehead~\cite{whitehead-generalized}.  The
suspension functor is not an isomorphism in the category of spectra, but
becomes an isomorphism when we invert the stable homotopy equivalences.
The resulting homotopy category of spectra is often called the stable
homotopy category and has been extensively studied, beginning with the
work of Boardman~\cite{boardman-vogt} and Adams~\cite{adams-blue} and
continuing to this day.  Notice that this definition of a spectrum can be
applied to any situation where one has an operation on a category that
one would like to invert; however, this simplest construction does not 
preserve the smash product structure coming from spaces.  

One of the stable homotopy category's basic features is that it is
symmetric monoidal.  There is a smash product, built from the smash
product of pointed topological spaces and analogous to the tensor
product of modules, that is associative, commutative, and unital, up to
coherent natural isomorphism.  However, the category of spectra defined
above is not symmetric monoidal.  This has been a sticking point for
almost forty years now.  Indeed, it was long thought that there could be
no symmetric monoidal category of spectra; see~\cite{lewis-convenient},
where it is shown that a symmetric monoidal category of spectra can not
have all the properties one might like.

Any good symmetric monoidal category of spectra allows one to perform
algebraic constructions on spectra that are impossible without such a
category.  This is extremely important, for example, in the algebraic
$K$-theory of spectra. In particular, given a good symmetric monoidal
category of spectra, it is possible to construct a homotopy category of
monoids (ring spectra) and of modules over a given monoid.

In this paper, we describe a symmetric monoidal category of spectra,
called the category of symmetric spectra.  The ordinary category of
spectra as described above is the category of modules over the sphere
spectrum.  The sphere spectrum is a monoid in the category of sequences
of spaces, but it is not a commutative monoid, because the twist map on
$S^{1}\sm S^{1}$ is not the identity.  This explains why the ordinary
category of spectra is not symmetric monoidal, just as in algebra where
the usual internal tensor product of modules is defined only over a
commutative ring.  To make the sphere spectrum a commutative monoid, we
need to keep track of the twist map, and, more generally, of
permutations of coordinates.  We therefore define a symmetric spectrum
to be a sequence of pointed simplicial sets $X_{n}$ together with a
pointed action of the permutation group $\Sigma _{n}$ on $X_{n}$ and
equivariant structure maps $S^{1}\sm X_{n}\xrightarrow{}X_{n+1}$.  We
must also require that the iterated structure maps $S^{p}\sm
X_{n}\xrightarrow{}X_{n+p}$ be $\Sigma _{p}\times \Sigma
_{n}$-equivariant.  This idea is due to the third author; the first and
second authors joined the project later.

At approximately the same time as the third author discovered symmetric
spectra, the team of Elmendorf, Kriz, Mandell, and
May~\cite{elmendorf-kriz-mandell-may} also constructed a symmetric
monoidal category of spectra, called $S$-modules.  Some generalizations
of symmetric spectra appear in~\cite{mandell-may-shipley-schwede-1}.
These many new symmetric monoidal categories of spectra, including
$S$-modules and symmetric spectra, are shown to be equivalent in an
appropriate sense in~\cite{mandell-may-shipley-schwede-1}
and~\cite{mandell-may-shipley-schwede-2}.  Another symmetric monoidal
category of spectra sitting between the approaches
of~\cite{elmendorf-kriz-mandell-may} and of this paper is developed
in~\cite{dwyer-shipley-hyper}.  We also point out that symmetric spectra
are part of a more general theory of localization of model
categories~\cite{hirschhorn}; we have not adopted this approach, but
both~\cite{hirschhorn} and~\cite{kan-model} have influenced us
considerably.

Symmetric spectra have already proved useful.
In~\cite{geisser-hesselholt}, symmetric spectra are used to extend the
definition of topological cyclic homology from rings to schemes.
Similarly, in~\cite{shipley-thh}, B\"{o}kstedt's approach to topological
Hochschild homology~\cite{bokstedt-thh} is extended to symmetric ring
spectra, without connectivity conditions.  And
in~\cite{schwede-shipley-stable}, it is shown that any linear model
category is Quillen equivalent to a model category of modules over a
symmetric ring spectrum.

As mentioned above, since the construction of symmetric spectra is
combinatorial in nature it may be applied in many different situations.
Given any well-behaved symmetric monoidal model category, such as chain
complexes, simplicial sets, or topological spaces, and an endofunctor on
it that respects the monoidal structure, one can define symmetric
spectra.  This more general approach is explored
in~\cite{hovey-stable-model}.  In particular, symmetric spectra may be
the logical way to construct a model structure for Voevodsky's stable
homotopy of schemes~\cite{voevodsky-milnor}.

In this paper, we can only begin the study of symmetric spectra.  The
most significant loose end is the construction of a model category of
commutative symmetric ring spectra; such a model category has been
constructed by the third author in work in progress.  It would also be
useful to have a stable fibrant replacement functor, as the usual
construction $QX$ does not work in general.  A good approximation to
such a functor is constructed in~\cite{shipley-thh}.

At present the theory of $S$-modules
of~\cite{elmendorf-kriz-mandell-may} is considerably more developed than
the theory of symmetric spectra.  Their construction appears to be
significantly different from symmetric spectra; however, there is work
in progress~\cite{mandell-may-shipley-schwede-2} showing that the two
approaches define equivalent stable homotopy categories and equivalent
homotopy categories of monoids and modules, as would be expected.  Each
approach has its own advantages.  The category of symmetric spectra is
technically much simpler than the $S$-modules
of~\cite{elmendorf-kriz-mandell-may}; this paper is almost entirely
self-contained, depending only on some standard results about simplicial
sets and topological spaces.  As discussed above, symmetric spectra can
be built in many different circumstances, whereas $S$-modules appear to
be tied to the category of topological spaces.  There are also technical
differences reflecting the result of~\cite{lewis-convenient} that there
are limitations on any symmetric monoidal category of spectra. For
example, the sphere spectrum $S$ is cofibrant in the category of
symmetric spectra, but is not in the category of $S$-modules.  On the
other hand, every $S$-module is fibrant, a considerable technical
advantage.  Also, the $S$-modules of~\cite{elmendorf-kriz-mandell-may}
are very well suited to the varying universes that arise in equivariant
stable homotopy theory, whereas we do not yet know how to realize
universes in symmetric spectra.  For a first step in this direction see
\cite{schwede-shipley-stable}.

{\em Organization.} 
The paper is organized as follows.  We choose to work in the category of
simplicial sets until Section~\ref{sec-topological-spectra}, where we
discuss topological symmetric spectra.  This is a significant technical
simplification; while it is possible to develop symmetric spectra in
topological spaces from scratch, technical issues arise that are not
present when working with simplicial sets.  In the first section, we
define symmetric spectra, give some examples, and establish some basic
properties.  In Section~\ref{sec-smash-product} we describe the closed
symmetric monoidal structure on the category of symmetric spectra, and
explain why such a structure can not exist in the ordinary category of
spectra.  In Section~\ref{sec-homotopy} we study the stable homotopy
theory of symmetric spectra.  This section is where the main subtlety of
the theory of symmetric spectra arises: we cannot define stable
equivalence by using stable homotopy isomorphisms.  Instead, we define a
map to be a stable equivalence if it is a cohomology isomorphism for all
cohomology theories.  The main result of this section is that symmetric
spectra, together with stable equivalences and suitably defined classes
of stable fibrations and stable cofibrations, form a model category.  As
expected, the fibrant objects are the $\Omega $-spectra; \ie symmetric
spectra $X$ such that each $X_{n}$ is a Kan complex and the adjoint
$X_{n}\xrightarrow{}X_{n+1}^{S^{1}}$ of the structure map is a weak
equivalence.  In Section~\ref{sec-comparison}, we prove that the stable
homotopy theories of symmetric spectra and ordinary spectra are
equivalent.  More precisely, we construct a Quillen equivalence of model
categories between symmetric spectra and the model category of ordinary
spectra described in~\cite{bousfield-friedlander}.

In Section~\ref{sec-misc} we discuss some of the properties of symmetric
spectra.  In particular, in Section~\ref{subsec-level-model}, we tie up
a loose end from Section~\ref{sec-homotopy} by establishing two
different model categories of symmetric spectra where the weak
equivalences are the level equivalences.  We characterize the stable
cofibrations of symmetric spectra in Section~\ref{subsec-stable-cofib}.
In Section~\ref{subsec-pushout-smash}, we show that the smash product of
symmetric spectra interacts with the model structure in the expected
way.  This section is crucial for the applications of symmetric spectra,
and, in particular, is necessary to be sure that the smash product of
symmetric spectra does define a symmetric monoidal structure on the
stable homotopy category.  We establish that symmetric spectra are a
proper model category in Section~\ref{subsec-proper}, and use this to
verify the monoid axiom in Section~\ref{subsec-monoids}.  The monoid
axiom is required to construct model categories of monoids and of
modules over a given monoid, see~\cite{schwede-shipley-monoids}.  In
Section~\ref{sec-semistable}, we define semistable spectra, which are
helpful for understanding the difference between stable equivalences and
stable homotopy equivalences.  Finally, we conclude the paper by
considering topological symmetric spectra in
Section~\ref{sec-topological-spectra}.

{\em Acknowledgments.}  
The authors would like to thank Dan Christensen, Bill Dwyer, Phil
Hirschhorn, Dan Kan, Haynes Miller, John Palmieri, Charles Rezk, and
Stefan Schwede for many helpful conversations about symmetric spectra.  
We would also like to thank Gaunce Lewis and Peter May for pointing out
that topological spaces are more complicated than we had originally
thought.  

{\em Notation.} 
We now establish some notation we will use throughout the paper.  Many
of the categories in this paper have an enriched Hom as well as a set
valued Hom. To distinguish them: in a category $\cc$, the set of maps
from $X$ to $Y$ is denoted $\cc(X,Y)$; in a simplicial category $\cc$,
the simplicial set of maps from $X$ to $Y$ is denoted $\shom_{\cc}(X,Y)$
or $\shom(X,Y)$; in a category $\cc$ with an internal Hom, the object in
$\cc$ of maps from $X$ to $Y$ is denoted $\ihom_{\cc}(X,Y)$ or
$\ihom(X,Y)$.  In case $\cc $ is the category of modules over a
commutative monoid $S$, we also use $\ihom _{S}(X,Y)$ for the internal
Hom.

\section{Symmetric spectra}\label{sec-spectra}

In this section we construct the category of symmetric spectra over
simplicial sets. We discuss the related category of topological
symmetric spectra in Section~\ref{sec-topological-spectra}. We begin
this section by recalling the basic facts about simplicial sets in
Section~\ref{subsec-ssets}, then we define symmetric spectra in Section
~\ref{subsec-sym-spectra}. We describe the simplicial structure on the
category of symmetric spectra in
Section~\ref{subsec-simplicial-structure}. The homotopy category of
symmetric $\Omega$-spectra is described in
Section~\ref{subsec-stable-category}.

\subsection{Simplicial sets}\label{subsec-ssets}

With the exception of Section~\ref{sec-topological-spectra} dealing with
topological spectra, this paper is written using simplicial sets. We
recall the basics. Consult \cite{may-simplicial-objects}
or \cite{curtis} for more details.

The category $\Delta$ has the ordered sets $[n]=\{0,1,\dots,n\}$ for
$n\ge0$ as its objects and the order preserving functions $[n]\to[m]$
as its maps. The category of \emph{simplicial sets},
denoted $\cs$, is the category of functors from $\Delta^{\textup{op}}$
to the category of sets. The set of $n$-simplices of the simplicial
set $X$, denoted $X_n$, is the value of the functor $X$ at
$[n]$. The standard $n$-simplex $\Delta [n]$ is the contravariant
functor $\Delta(-,[n])$.  Varying $n$ gives a covariant functor
$\Delta[-]\mathcolon\Delta\to\cs$.  By the Yoneda lemma,
$\cs(\Delta[n],X)=X_n$ and the contravariant functor
$\cs(\Delta[-],X)$ is naturally isomorphic to $X$. 

Let $G$ be a discrete group. The category of \emph{$G$-simplicial sets}
is the category $\cs^G$ of functors from $G$ to $\cs$, where $G$ is
regarded as a category with one object.  A $G$-simplicial set is
therefore a simplicial set $X$ with a left simplicial $G$-action, \ie a
homomorphism $G\to\cs(X,X)$.

A \emph{basepoint} of a simplicial set $X$ is a distinguished
$0$-simplex $*\in X_0$. The category of pointed simplicial sets and
basepoint preserving maps is denoted $\sset$. The simplicial set
$\Delta[0]=\Delta(-,[0])$ has a single simplex in each degree and is the
terminal object in $\cs $.  A basepoint of $X$ is the same as a map
$\Delta[0]\to X$. The disjoint union $X_+=X\amalg\Delta[0]$ adds a
disjoint basepoint to the simplicial set $X$. For example, the
$0$-sphere is $S^0=\Delta[0]_+$. A basepoint of a $G$-simplicial set $X$
is a $G$-invariant $0$-simplex of $X$. The category of pointed
$G$-simplicial sets is denoted $\sset^G$.

The smash product $X\sm Y$ of the pointed simplicial sets $X$ and $Y$
is the quotient $(X\times Y)/(X\vee Y)$ that collapses the simplicial
subset $X\vee Y=X\times *\cup *\times Y$ to a point. For pointed
$G$-simplicial sets $X$ and $Y$, let $X\sm_G Y$ be the quotient of
$X\sm Y$ by the diagonal action of $G$.  For pointed simplicial sets
$X$, $Y$, and $Z$, there are natural isomorphisms $(X\sm Y)\sm Z\natiso
X\sm(Y\sm Z)$, $X\sm Y\natiso Y\sm X$ and $X\sm S^0 \natiso X$. In the
language of monoidal categories, the smash product is a symmetric
monoidal product on the category of pointed simplicial sets. We recall
the definition of symmetric monoidal product, but for more details see
\cite[VII]{maclane-categories} or \cite[6.1]{borceux}.

\begin{definition}\label{def-monoidal-cat}
A \emph{symmetric monoidal product} on a category $\cc$ is: a bifunctor
$\otimes\mathcolon\cc \times\cc \to\cc $; a unit $U\in\cc $; and
coherent natural isomorphisms $(X\otimes Y)\otimes Z\natiso X\otimes
(Y\otimes Z)$ (the associativity isomorphism), $X\otimes Y\natiso
Y\otimes X$ (the twist isomorphism), and $U\otimes X\natiso X$ (the unit
isomorphism). The product is \emph{closed} if the functor $X\otimes (-)$
has a right adjoint $\ihom(X,-)$ for every $X\in\cc$. A
\emph{\textrm{(}closed\textrm{)} symmetric monoidal category} is a
category $\cc$ with a (closed) symmetric monoidal product.
\end{definition}

Coherence of the natural isomorphisms means that all reasonable diagrams
built from the natural isomorphisms also commute
\cite{maclane-categories}.  When the product is closed, the pairing
$\Hom(X,Y)\mathcolon \cc^{\textup{op}}\times\cc \to \cc$ is an internal
Hom. For example, the smash product on the category $\sset$ of pointed
simplicial sets is closed.  For $X,Y\in\sset$, the pointed simplicial
set of maps from $X$ to $Y$ is
$\shom_{\sset}(X,Y)=\sset(X\sm\Delta[-]_+,Y)$. For pointed
$G$-simplicial sets $X$ and $Y$, the simplicial subset of
$G$-equivariant pointed maps is
$\shom_G(X,Y)=\sset^G(X\sm\Delta[-]_+,Y)$.

\subsection{Symmetric spectra}\label{subsec-sym-spectra}

Let $S^1$ be the simplicial circle $\Delta [1]/\partial\Delta [1]$,
obtained by identifying the two vertices of $\Delta[1]$.

\begin{definition}\label{def-spectrum}
A \emph{spectrum} is 
\begin{enumerate}
\item a sequence $X_0,X_1,\dots,X_n,\dots$ of pointed simplicial sets; and 
\item a pointed map $\sigma\mathcolon S^1\sm X_n\to X_{n+1}$ for each
$n\ge0$.
\end{enumerate}
The maps $\sigma$ are the \emph{structure maps} of the spectrum.  A
\emph{map of spectra} $f\mathcolon X\to Y$ is a sequence of pointed maps
$f_n\mathcolon X_n\to Y_n$ such that the diagram
\[
\begin{CD}
S^1\sm X_n @>\sigma>> X_{n+1}\\ 
@VS^1\sm f_nVV @Vf_{n+1}VV \\ 
S^1\sm Y_n @>\sigma>> Y_{n+1}
\end{CD}
\]
is commutative for each $n\ge0$. Let $\BF$ denote the category of
spectra.
\end{definition}
  
Replacing the sequence of pointed simplicial sets by a sequence of
pointed topological spaces in~\ref{def-spectrum} gives the original
definition of a spectrum (due to Whitehead and Lima). The categories
of simplicial spectra and of topological spectra are discussed in the work of
Bousfield and Friedlander~\cite{bousfield-friedlander}.

A symmetric spectrum is a spectrum to which symmetric group actions
have been added.  Let $\Sigma_p$ be the group of permutations of the
set $\ov p=\{1,2,\dots,p\}$, with $\ov 0=\emptyset$.  As usual, embed
$\Sigma_p\times \Sigma_q$ as the subgroup of $\Sigma_{p+q}$ with
$\Sigma_p$ acting on the first $p$ elements of $\ov{p+q}$ and
$\Sigma_q$ acting on the last $q$ elements of $\ov{p+q}$. Let
$S^p=(S^{1})^{\wedge p}$ be the $p$-fold smash power of the simplicial
circle with the left permutation action of $\Sigma_p$.

\begin{definition}\label{def-symm-spectra}
A \emph{symmetric spectrum} is 
\begin{enumerate}
\item a sequence $X_0,X_1,\dots,X_n,\dots$ of pointed simplicial sets; 
\item a pointed map $\sigma\mathcolon S^1\sm X_n\to X_{n+1}$ for each
$n\ge0$; and
\item a basepoint preserving left action of $\Sigma_n$ on $X_n$ such
that the composition
\[
\sigma^p=\sigma\circ (S^1\sm\sigma)\circ\dots\circ (S^{p-1}\sm
\sigma)\mathcolon S^p\sm X_n\to X_{n+p},
\] 
of the maps $S^i\sm S^1\sm X_{n+p-i-1}\xrightarrow{S^i\sm\sigma} S^{
i}\sm X_{n+p-i}$ is $\Sigma_p\times\Sigma_n$-equivariant for $p\ge1$
and $n\ge0$.
\end{enumerate}
A \emph{map of symmetric spectra} $f\mathcolon X\to Y$ is a sequence of
pointed maps $f_n\mathcolon X_n\to Y_n$ such that $f_n$ is
$\Sigma_n$-equivariant and the 
diagram
\[
\begin{CD}
S^1\sm X_n @>\sigma>> X_{n+1}\\ 
@VS^1\sm f_nVV @Vf_{n+1}VV \\ 
S^1\sm Y_n @>\sigma>> Y_{n+1}
\end{CD}
\] 
is commutative for each $n\ge0$. Let $\spec$ denote the category of
symmetric spectra.
\end{definition}

\begin{remark}\label{rem-transpositions}
In part three of Definition~\ref{def-symm-spectra}, one need only
assume that the maps $\sigma\mathcolon S^1\sm X_n\to X_{n+1}$ and
$\sigma^2\mathcolon S^2\sm X_n\to X_{n+2}$ are equivariant; since the
symmetric groups $\Sigma_p$ are generated by transpositions $(i,i+1)$,
if $\sigma$ and $\sigma^2$ are equivariant then all the maps
$\sigma^p$ are equivariant.
\end{remark}

\begin{example}\label{ex-symmetric-spectra}
The \emph{symmetric suspension spectrum} $\Sigma^\infty K$ of the
pointed simplicial set $K$ is the sequence of pointed simplicial sets
$S^n\sm K$ with the natural isomorphisms $\sigma\mathcolon S^1\sm S^n\sm
K\to S^{n+1}\sm K$ as the structure maps and the diagonal action of
$\Sigma_n$ on $S^n\sm K$ coming from the left permutation action on
$S^n$ and the trivial action on $K$. The composition $\sigma^p$ is the
natural isomorphism which is $\Sigma_p\times\Sigma_n$-equivariant. The
\emph{symmetric sphere spectrum} $S$ is the symmetric suspension
spectrum of the $0$-sphere; $S$ is the sequence of spheres
$S^0,S^1,S^2,\dots$ with the natural isomorphisms $S^1\sm S^n\to
S^{n+1}$ as the structure maps and the left permutation action of
$\Sigma_n$ on $S^n$.
\end{example}

\begin{example}\label{ex-homology}
The \emph{Eilenberg-Mac Lane spectrum} $H\mathbb{Z}$ is the sequence of
simplicial abelian groups $\mathbb{Z}\otimes S^{n}$, where
$(\mathbb{Z}\otimes S^{n})_{k}$ is the free abelian group on the
non-basepoint $k$-simplices of $S^{n}$.  We identify the basepoint with
$0$.  The symmetric group $\Sigma _{n}$ acts by permuting the
generators, and one can easily verify that the evident structure maps
are equivariant.  One could replace $\mathbb{Z}$ by any ring.  
\end{example}

\begin{remark}\label{more-ex}
Bordism is most easily defined as a topological symmetric spectrum, see
Example~\ref{bordism}.  As explained in~\cite[Section
6]{geisser-hesselholt}, many other examples of symmetric spectra arise
as the $K$-theory of a category with cofibrations and weak equivalences
as defined by Waldhausen~\cite[p.330]{wal}.  
\end{remark}

A symmetric spectrum with values in a simplicial category $\cc$ is
obtained by replacing the sequence of pointed simplicial sets by a
sequence of pointed objects in $\cc$. In particular, a topological
symmetric spectrum is a symmetric spectrum with values in the
simplicial category of topological spaces; see
Section~\ref{sec-topological-spectra}.

By ignoring group actions, a symmetric spectrum is a
spectrum and a map of symmetric spectra is a map of spectra.  
When no confusion can arise, the adjective ``symmetric'' may be dropped.

\begin{definition}\label{def-underlying}
Let $X$ be a symmetric spectrum. The \emph{underlying spectrum} $UX$ is
the sequence of pointed simplicial sets $(UX)_n=X_n$ with the same
structure maps $\sigma\mathcolon S^1\wedge (UX)_n \to (UX)_{n+1}$ as $X$
but ignoring the symmetric group actions.  This gives a faithful functor
$U\mathcolon\spec \to \BF$.
\end{definition} 

Since the action of $\Sigma_n$ on $S^n$ is non-trivial for $n\ge 2$, it
is usually impossible to obtain a symmetric spectrum from a spectrum
by letting $\Sigma_n$ act trivially on $X_n$.  However, many of the
usual functors to the category of spectra lift to the category
of symmetric spectra.  For example, the suspension spectrum of
a pointed simplicial set $K$ is the underlying spectrum of the
symmetric suspension spectrum of $K$.

Many examples of symmetric spectra and of functors on the category of
symmetric spectra are constructed by prolongation of simplicial
functors.

\begin{definition}
A \emph{pointed simplicial functor} or \emph{$\sset$-functor} is a
pointed functor $R\mathcolon\sset\to\sset$ and a natural transformation
$h\mathcolon RX\sm K\to R(X\sm K)$ of bifunctors such that the
composition $RX\sm S^0 \to R(X\sm S^0)\to R(X)$ is the unit isomorphism
and the diagram of natural transformations
\[
\xymatrix{ {(RX\sm K)\sm L}\ar[r]^-{h\sm L}\ar[d] & 
{R(X\sm K)\sm L}\ar[d]_{h} \\
{RX \sm (K\sm L)} \ar[r]^-{h} & {R(X\sm K\sm L)} }
\]
is commutative.  A \emph{pointed simplicial natural transformation}, or
\emph{$\sset $-natural transformation}, from the $\sset $-functors $R$
to the $\sset $-functor $R'$ is a natural transformation $\tau
\mathcolon R\to R'$ such that $\tau h=h'(\tau \sm K)$.  
\end{definition}

\begin{definition}\label{def-prolongation}
The \emph{prolongation} of a $\sset$-functor $R\mathcolon \sset \to \sset$
is the functor $R\mathcolon\spec\to\spec$ defined as follows. For $X$ a
symmetric spectrum, $RX$ is the sequence of pointed simplicial sets
$RX_n$ with the composition $\sigma\mathcolon S^1\sm R(X_n)\to R(S^1\sm
X_n)\xrightarrow{R\sigma} R(X_{n+1})$ as the structure maps and the
action of $\Sigma_n$ on $R(X_n)$ obtained by applying the functor $R$ to
the action of $\Sigma_n$ on $X_n$.  Since $R$ is a $\sset$-functor,
each map $\sigma^p$ is equivariant and so $RX$ is a symmetric spectrum.
For $f$ a map of symmetric spectra, $Rf$ is the sequence of pointed maps
$Rf_n$. Since $R$ is an $\sset$-functor, $Rf$ is a map of spectra.
Similarly, we can prolong an $\sset $-natural transformation to a
natural transformation of functors on $\spec $.  
\end{definition}

\begin{proposition}\label{prop-spec-bicomplete}
The category of symmetric spectra is bicomplete (every small diagram
has a limit and a colimit).
\end{proposition}

\begin{proof}
For any small category $I$, the limit and colimit functors $\sset^I\to
\sset$ are pointed simplicial functors; for $K\in\sset$ and $D\in\set^I$
there
is a natural isomorphism
\[
K\sm\colim D \natiso \colim (K\sm D)
\]
and a natural map
\[
K\sm \lim D \to \lim (K\sm D).
\]
A slight generalization of prolongation gives the limit and the colimit
of a diagram of symmetric spectra.
\end{proof}

In particular, the underlying sequence of
the limit is $(\lim D)_n=\lim D_n$ and the underlying sequence of the
colimit is $(\colim D)_n=\colim D_n$.

\subsection{Simplicial structure on $\spec$}
\label{subsec-simplicial-structure}

For a pointed simplicial set $K$ and  a symmetric spectrum $X$,
prolongation of the $\sset$-functor $(-)\sm K \mathcolon \sset\to\sset$
defines the \emph{smash product} $X\sm K$ and prolongation of the
$\sset$-functor $(-)^K\mathcolon\sset\to\sset$ defines the \emph{power
spectrum} $X^K$.  For symmetric spectra $X$ and $Y$, the \emph{pointed
simplicial set} of maps from $X$ to $Y$ is
$\shom_{\spec}(X,Y)=\spec(X\sm \Delta[-]_+, Y)$.

In the language of enriched category theory, the following proposition
says that the smash product $X \sm K$ is a closed action of $\sset$ on
$\spec$.  We leave the straightforward proof to the reader. 

\begin{proposition}\label{prop-closed-S-action-spec}
Let $X$ be a symmetric spectrum. Let $K$ and $L$ be pointed
simplicial sets.
\begin{enumerate}
\item  There are coherent natural isomorphisms
$X\sm (K\sm L)\natiso(X \sm K) \sm L$ and $X\sm S^0\natiso X$.
\item $(-)\sm K\mathcolon\spec \to\spec $ is the left adjoint of the
functor $(-)^K\mathcolon\spec \to\spec$.
\item $X\sm (-)\mathcolon\sset\to\spec$ is the left adjoint of the
functor $\shom_{\spec}(X,-)\mathcolon\spec\to\sset$.
\end{enumerate}
\end{proposition}

The \emph{evaluation map} $X\sm\shom_{\spec } (X,Y)\to Y$ is the adjoint
of the identity map on $\shom_{\spec }(X,Y)$. The \emph{composition
pairing}
\[
\shom_{\spec }(X,Y)\sm\shom_{\spec }(Y,Z)\to \shom_{\spec }(X,Z)
\]
is the adjoint of the composition 
\[
X\sm\shom_{\spec }(X,Y)\sm\shom_{\spec }(Y,Z)\to Y\sm\shom_{\spec
}(Y,Z)\to Z
\]
of two evaluation maps. In the language of enriched category theory, a
category with a closed action of $\sset$ is the same as a tensored and
cotensored $\sset$-category.  The following proposition, whose proof we
also leave to the reader, expresses this fact.

\begin{proposition}\label{prop-spec-is-S-category}
Let $X$, $Y$, and $Z$ be symmetric spectra and let $K$ be a pointed
simplicial set.
\begin{enumerate}
\item The composition pairing
$\shom_{\spec}(X,Y)\sm\shom_{\spec}(Y,Z)\to\shom_{\spec}(X,Z)$ is
associative.
\item The adjoint $S^0\to\shom_{\spec }(X,X)$ of the isomorphism
$X\sm S^0\to X$ is a left and a right unit of the composition pairing.
\item There are natural isomorphisms 
\[
\shom_{\spec}(X\sm K,Y)\natiso\shom_{\spec}(X,Y^K)\natiso\shom_{\spec}(X,Y)^K.
\]
\end{enumerate}
\end{proposition}

Proposition~\ref{prop-closed-S-action-spec} says that certain functors
are adjoints, whereas Proposition~\ref{prop-spec-is-S-category} says
more; they are simplicial adjoints.

The category of symmetric spectra satisfies Quillen's axiom SM7 for
simplicial model categories.

\begin{definition}\label{def-pushout-smash}
Let $f\mathcolon U\to V$ and $g\mathcolon X\to Y$ be maps of pointed
simplicial sets. The \emph{pushout smash product} $f\boxprod g$ is the
natural map on the pushout
\[
f\boxprod g\mathcolon V\sm X\amalg_{U\sm X}U\sm Y\to V\sm Y.
\]
induced by the commutative square
\[
\begin{CD}
U\sm X @>f\sm X>> V\sm X \\ @VU\sm gVV @VV{V\sm g}V \\ U\sm Y @>>f\sm
Y> V\sm Y.
\end{CD}
\]
Let $f$ be a map of symmetric spectra and $g$ be a map of pointed
simplicial sets. The \emph{pushout smash product} $f\boxprod g$ is
defined by prolongation, $(f\boxprod g)_n=f_n\boxprod g$. 
\end{definition}

Recall that a map of simplicial sets is a weak equivalence if its
geometric realization is a homotopy equivalence of CW-complexes.  One of
the basic properties of simplicial sets, proved
in~\cite[II.3]{quillen-htpy}, is:

\begin{proposition}\label{prop-SM7-sset}
Let $f$ and $g$ be monomorphisms of pointed simplicial sets.  Then
$f\boxprod g$ is a monomorphism, which is a weak equivalence if either
$f$ or $g$ is a weak equivalence.  
\end{proposition}

Prolongation gives a corollary for symmetric spectra. A map $f$ of
symmetric spectra is a monomorphism if $f_n$ is a monomorphism of
simplicial sets for each $n\ge0$.

\begin{definition}
A map $f$ of symmetric spectra is a \emph{level equivalence} if 
$f_n$ is a weak equivalence of simplicial sets for each $n\ge0$.
\end{definition}

\begin{corollary}\label{cor-SM7-spec-sset}
Let $f$ be a monomorphism of symmetric spectra and let $g$ be a
monomorphism of pointed simplicial sets. Then $f\boxprod g$ is a
monomorphism, which is a level equivalence if either $f$ is a level
equivalence or $g$ is a weak equivalence. 
\end{corollary}

By definition, a $0$-simplex of $\shom_{\spec}(X,Y)$ is a map
$X\sm\Delta[0]_+\to Y$, but $X\sm\Delta[0]_+\natiso X$ and so a
$0$-simplex of $\shom_{\spec}(X,Y)$ is a map $X\to Y$. A $1$-simplex of
$\shom_{\spec}(X,Y)$ is a \emph{simplicial homotopy} $H\mathcolon
X\sm\Delta[1]_+\to Y$ from $H\circ(X\sm i_0)$ to $H\circ(X\sm i_1)$
where $i_0$ and $i_1$ are the two inclusions
$\Delta[0]\to\Delta[1]$. Simplicial homotopy generates an equivalence
relation on $\spec (X,Y)$ and the quotient is $\pi_0\shom_{\spec
}(X,Y)$. A map $f\mathcolon X \to Y$ is a \emph{simplicial homotopy
equivalence} if it has a simplicial homotopy inverse, \ie a map
$g\mathcolon Y\to X$ such that $gf$ is simplicially homotopic to the
identity map on $X$ and $fg$ is simplicially homotopic to the identity
map on $Y$. If $f$ is a simplicial homotopy equivalence of symmetric
spectra, then each of the maps $f_n$ is a simplicial homotopy
equivalence, and so each of the maps $f_n$ is a weak equivalence. Every
simplicial homotopy equivalence is therefore a level equivalence.  The
converse is false; a map can be a level equivalence and NOT a simplicial
homotopy equivalence.

\subsection{Symmetric $\Omega$-spectra}\label{subsec-stable-category}

The stable homotopy category can be defined using $\Omega$-spectra and
level equivalences.

\begin{definition}
A Kan complex (see Example~\ref{eg-model-cats}) is a simplicial set that
satisfies the Kan extension condition.  An \emph{$\Omega$-spectrum} is a
spectrum $X$ such that for each $n\ge0$ the simplicial set $X_n$ is a
Kan complex and the adjoint $X_n\to \shom_{\sset }(S^1,X_{n+1})$ of the
structure map $S^1\sm X_n\to X_{n+1}$ is a weak equivalence of
simplicial sets.
\end{definition}

Let $\Omega\BF\subseteq\BF$ be the full subcategory of $\Omega$-spectra.
The homotopy category $\ho(\Omega\BF)$ is obtained from $\Omega\BF$ by
formally inverting the level equivalences.  By the results in
\cite{bousfield-friedlander}, the category $\ho(\Omega\BF)$ is naturally
equivalent to Boardman's stable homotopy category (or any other).
Likewise, let $\Omega\spec\subseteq\spec$ be the full subcategory of
symmetric $\Omega$-spectra (\ie symmetric spectra $X$ for which $UX$ is
an $\Omega$-spectrum). The homotopy category $\ho(\Omega\spec)$ is
obtained from $\Omega\spec$ by formally inverting the level
equivalences. Since the forgetful functor $U\mathcolon\spec\to\BF$
preserves $\Omega$-spectra and level equivalences, it induces a functor
$\ho(U)\mathcolon\ho(\Omega\spec)\to \ho(\Omega\BF)$. As a corollary of
Theorem~\ref{thm-comparison}, the functor $\ho(U)$ is a natural
equivalence of categories. Thus the category $\ho(\Omega\spec)$ is
naturally equivalent to Boardman's stable homotopy category.  To
describe an inverse of $\ho(U)$, let $\Omega^\infty \mathcolon
\BF\to\sset$ be the functor that takes a spectrum to the $0$-space of
its associated $\Omega$-spectrum.  For any spectrum $E\in\BF$, the
symmetric spectrum $VE=\Omega^\infty(E\sm S)$ is the value of the
prolongation of the $\sset$-functor $\Omega^\infty(E\sm -)$ at the
symmetric sphere spectrum $S$; the underlying sequence is
$VE_n=\Omega^\infty(E\sm S^n)$. The functor $V$ preserves
$\Omega$-spectra, preserves level equivalences, and induces a functor
$\ho(V)\mathcolon\ho(\Omega\BF)\to \ho(\Omega\spec)$ which is a natural
inverse of $\ho(U)$.

The category of symmetric $\Omega$-spectra has major defects. It is not
closed under limits and colimits, or even under pushouts and
pullbacks. The smash product, defined in
Section~\ref{sec-smash-product}, of symmetric $\Omega$-spectra is a
symmetric spectrum but not an $\Omega$-spectrum, except in trivial
cases.  For these reasons it is better to work with the category of all
symmetric spectra. But then the notion of level equivalence is no longer
adequate; the stable homotopy category is a retract of the homotopy
category obtained from $\spec$ by formally inverting the level
equivalences but many symmetric spectra are not level equivalent to an
$\Omega$-spectrum.  One must enlarge the class of equivalences.  The
stable equivalences of symmetric spectra are defined in
Section~\ref{subsec-stable-equivalence}. By
Theorem~\ref{thm-comparison}, the homotopy category obtained from
$\spec$ by inverting the stable equivalences is naturally equivalent to
the stable homotopy category.

\section{The smash product of symmetric spectra}
\label{sec-smash-product}

In this section we construct the closed symmetric monoidal product on
the category of symmetric spectra. A symmetric spectrum can be viewed as
a module over the symmetric sphere spectrum $S$, and the symmetric
sphere spectrum (unlike the ordinary sphere spectrum) is a commutative
monoid in an appropriate category.  The smash product of symmetric
spectra is the tensor product over $S$.

The closed symmetric monoidal category of symmetric sequences is
constructed in Section~\ref{subsec-sequences}. A reformulation of the
definition of a symmetric spectrum is given in
Section~\ref{subsec-spectra} where we recall the definition of monoids
and modules in a symmetric monoidal category.  In
Section~\ref{subsec-BF} we see that there is no closed symmetric
monoidal smash product on the category of (non-symmetric) spectra.  

\subsection{Symmetric sequences}\label{subsec-sequences}

Every symmetric spectrum has an underlying sequence
$X_0,X_1,\dots,X_{n},\dots$ of pointed simplicial sets with a basepoint
preserving left action of $\Sigma_n$ on $X_n$; these are
called symmetric sequences. In this section we define the closed
symmetric monoidal category of symmetric sequences of pointed
simplicial sets.

\begin{definition}\label{def-symmetric-sequences}
The \emph{category} $\Sigma=\coprod_{n\ge0}\Sigma_n $ has the finite
sets $\ov n=\{1,2,\dots,n\}$ for $n \ge 0$ ($\ov 0=\emptyset$) as its
objects and the automorphisms of the sets $\ov n$ as its maps. Let $\cc$
be a category. A \emph{symmetric sequence} of objects in $\cc$ is a
functor $\Sigma\to\cc$, and the category of symmetric sequences of
objects in $\cc$ is the functor category $\cc^\Sigma$.
\end{definition}

A symmetric sequence $X\in\symseq$ is a sequence
$X_0,X_1,\dots,X_{n},\dots$ of pointed simplicial sets with a
basepoint preserving left action of $\Sigma_n$ on $X_n$.  The
category $\cc^\Sigma$ is a product category. In
particular, $\symseq(X,Y)=\prod_p\sset ^{\Sigma _{p}}(X_p,Y_p)$.

\begin{proposition}
The category $\symseq$ of symmetric sequences in $\sset$ is bicomplete.
\end{proposition}

\begin{proof}
The category $\sset$ is bicomplete, so the functor category $\symseq$ is
bicomplete.
\end{proof}

\begin{definition}\label{def-tensor-sequences}
The \emph{tensor product} $X\otimes Y$ of the symmetric sequences
$X,Y\in\symseq$ is the
symmetric sequence
\[
(X\otimes Y)_{n}= \bigvee_{p+q=n}(\Sigma _{n})_{+}\sm_{\Sigma _{p}\times
\Sigma _{q}} (X_{p}\sm Y_{q}).
\]
The \emph{tensor product} $f\otimes g\mathcolon X\otimes
Y\xrightarrow{}X'\otimes Y'$ of the maps $f\mathcolon X\xrightarrow{}X'$
and $g\mathcolon Y\xrightarrow{}Y'$ in $\symseq$ is given by $(f\otimes
g)(\alpha,x,y )=(\alpha,f_px,g_qy)$ for $\alpha\in\Sigma_{p+q}$, $x\in
X_p$ and $y\in Y_q$.
\end{definition}

The tensor product of symmetric sequences has the universal property
for ``bilinear maps'':

\begin{proposition}\label{prop-bilinear-hom}
Let $X,Y,Z\in\symseq$ be symmetric sequences. There is a natural isomorphism
\[
\symseq(X\otimes Y,Z)\natiso\prod_{p,q}\sset ^{\Sigma _{p}\times \Sigma
_{q}} (X_p\sm Y_q,Z_{p+q})
\]
\end{proposition}

The \emph{twist isomorphism} $\tau\mathcolon X \otimes Y \to Y \otimes X$
for $X,Y\in\symseq$ is the natural map given by
$\tau(\alpha,x,y)=(\alpha\rho_{q,p},y,x)$ for $\alpha\in\Sigma_{p+q}$,
$x\in X_p$, and $y\in Y_q$, where $\rho_{q,p}\in\Sigma_{p+q}$ is the
$(q,p)$-shuffle given by $\rho_{q,p}(i)=i+p$ for $1\le i\le q$ and
$\rho_{q,p}(i)=i-q$ for $q<i\le p+q$.  The map defined without the
shuffle permutation is not a map of symmetric sequences.

\begin{remark}\label{rem-twist-isomorphism}
There is another way of describing the tensor product and the twist
isomorphism. The category $\Sigma$ is a skeleton of the category of
finite sets and isomorphisms.  Hence every symmetric sequence has an
extension, which is unique up to isomorphism, to a functor on
the category of all finite sets and isomorphisms. The tensor product of
two such functors $X$ and $Y$ is the functor defined on a finite set $C$
as 
\[
(X\otimes Y)(C)=\bigvee _{A\cup B=C,A\cap B=\emptyset}X(A)\sm Y(B).
\]
For an isomorphism $f\mathcolon C \to D$ the map $(X\otimes Y)(f)$ is the
coproduct of the isomorphisms $X(A)\sm Y(B)\to X(fA)\sm Y(fB)$. The
twist isomorphism is the map that sends the summand $X(A)\sm Y(B)$ of
$(X\otimes Y)(C)$ to the summand $Y(B)\sm X(A)$ of $(Y\otimes X)(C)$
by switching the factors.
\end{remark}

\begin{lemma}\label{lem-tensor-sequences}
The tensor product $\otimes$ is a symmetric monoidal product on the category
of symmetric sequences $\symseq$.
\end{lemma}

\begin{proof}
The unit of the tensor product is the symmetric sequence
$\Sigma(\ov 0,-)_+=(S^0,*,*,\dots)$. The unit isomorphism is obvious.  The
associativity isomorphism is induced by the associativity isomorphism
in $\sset $ and the natural isomorphism
\[
((X\otimes Y)\otimes Z)_n\natiso \bigvee
_{p+q+r=n}(\Sigma_n)_+\sm_{\Sigma_p\times\Sigma_q\times\Sigma_r}(X_p\sm
Y_q\sm Z_r).
\] 
The twist isomorphism is described in
Remark~\ref{rem-twist-isomorphism}. The coherence of the natural
isomorphisms follows from coherence of the natural isomorphisms for the
smash product in $\sset$.
\end{proof}

We now introduce several functors on the category of symmetric sequences.

\begin{definition}
\label{def-sequence-evaluation}
The \emph{evaluation functor} $\Ev_n\mathcolon\symseq \to\sset$ is given
by $\Ev_nX=X_n$ and $\Ev_nf=f_n$. The \emph{free functor}
$G_n\mathcolon\sset\to\symseq$ is the left adjoint of the evaluation
functor $\Ev_n$.  The \emph{smash product} $X\sm K$ of $X\in\symseq$
and $K\in\sset$ is the symmetric sequence $(X\sm K)_n=X_n\sm K$ with
the diagonal action of $\Sigma_n$ that is trivial on $K$. The
\emph{pointed simplicial set} $\shom_{\symseq}(X,Y)$ of maps from
$X$ to $Y$ is the pointed simplicial set
$\symseq(X\sm\Delta[-]_+, Y)$.
\end{definition}

For each $n\ge0$, the free symmetric sequence is $\Sigma[n]=\Sigma(\ov
n,-)$ and the free functor is $G_n=\Sigma[n]_+ \sm -\mathcolon \sset \to
\symseq$.  So, for a pointed simplicial set $K$,
$(G_nK)_n=(\Sigma_n)_+\sm K$ and $(G_nK)_k=*$ for $k\neq n$. In
particular, $G_nS^0=\Sigma[n]_+$, $G_0K=(K,*,*,\dots)$ and $G_0S^0$ is
the unit of the tensor product $\otimes$.

We leave the proof of the following basic proposition to the reader.

\begin{proposition}\label{prop-symseq-nat-isos}
There are natural isomorphisms\textup{:}
\begin{enumerate}
\item $G_pK\otimes G_qL\natiso G_{p+q}(K\sm L)$ for $K,L\in\sset$.
\item $X\otimes G_0K\natiso X\sm K$ for $K\in\sset$ and $X\in\symseq$.
\item $\shom_{\symseq }(G_nK,X)\natiso \shom_{\sset }(K,X_n)$ for
$K\in\sset$ and $X\in\symseq$.
\item $\shom_{\symseq }(X\otimes Y,Z)\natiso
\prod_{p,q}\shom_{\Sigma_p\times\Sigma_q}(X_p\sm Y_q, Z_{p+q})$ for
$X,Y,Z\in\spec$.
\end{enumerate}
\end{proposition}

A map $f$ of symmetric sequences is a level equivalence if each of the
maps $f_n$ is a weak equivalence. Since $\symseq$ is a product
category, a map $f$ of symmetric sequences is a monomorphism if and
only if each of the maps $f_n$ is a monomorphism.

\begin{proposition}\label{prop-otimes-preserves}
Let $X$ be a symmetric sequence, $f$ be a map of symmetric
sequences and $g$ be a map of pointed simplicial sets.
\begin{enumerate}
\item $X\otimes(-)$ preserves colimits. 
\item If $f$ is a monomorphism then $X\otimes f$ is a  monomorphism.
\item If $f$ is a level equivalence then $X\otimes f$ is a level equivalence.
\item If $g$ is a monomorphism then $G_ng$ is a
monomorphism for $n\ge0$. 
\item If $g$ is a weak equivalence then $G_ng$ is a
level equivalence for $n\ge0$.
\end{enumerate}
\end{proposition}

\begin{proof}
Parts (1), (2) and (3) follow from the definition of $\otimes$
and the corresponding properties for the smash product of pointed
simplicial sets. For Parts (4) and (5) use the isomorphism
$G_nK=\Sigma[n]_+\sm K$.
\end{proof}

By part three of Proposition~\ref{prop-symseq-nat-isos},
$\shom(\Sigma[n]_+,X)\natiso X_n$.  As $n$ varies, $\Sigma[-]_+$ is a
functor $\Sigma^{\textup{op}}\to\symseq$, and for $X\in\symseq$, the
symmetric sequence $\shom_{\symseq }(\Sigma[-]_+,X)$ is naturally
isomorphic to $X$.

\begin{definition}\label{def-sequence-hom} 
Let $X$ and $Y$ be symmetric sequences. The \emph{symmetric sequence of
maps} from $X$ to $Y$ is
\[
\ihom_{\Sigma}(X,Y)=\shom_{\symseq }(X\otimes \Sigma[-]_+,Y).
\]
\end{definition}

\begin{theorem}\label{thm-tensor-closed-sym-monoidal}
The tensor product is a closed symmetric monoidal product on the
category of symmetric sequences.
\end{theorem}

\begin{proof}
The tensor product is a symmetric monoidal product by
Lemma~\ref{lem-tensor-sequences}. The product is closed if there is a
natural isomorphism
\[
\symseq(X \otimes Y, Z) \natiso \symseq\left(X,\ihom_{\Sigma}(Y,Z)\right).
\]
for symmetric sequences $X,Y$ and $Z$.

By Proposition~\ref{prop-bilinear-hom}, a map of symmetric sequences
$f\mathcolon X\otimes Y\to Z$ is a collection of
$\Sigma_p\times\Sigma_q$-equivariant maps $f_{p,q}\mathcolon X_p\sm
Y_q\to Z_{p+q}$. This is adjoint to a collection of
$\Sigma_p$-equivariant maps $g_{p,q}\mathcolon
X_p\to\shom_{\Sigma_q}(Y_q,Z_{p+q})$. So there is a natural isomorphism
\[
\symseq(X\otimes Y,Z)\natiso
\prod_p\sset ^{\Sigma _{p}}(X_p,\prod_q\shom_{\Sigma_q}(Y_q,Z_{p+q})) 
\]
By Proposition~\ref{prop-symseq-nat-isos}, the functor sending $\ov p$
to $\prod_q\shom_{\Sigma_q}(Y_q,Z_{p+q})$ is the functor sending $\ov p$
to $\shom(Y\otimes \Sigma[p]_+,Z)$ which by definition is
$\ihom_{\Sigma}(Y,Z)$. Combining the isomorphisms gives the natural
isomorphism that finishes the proof.
\end{proof}

\subsection{Symmetric spectra}\label{subsec-spectra}

In this section we recall the language of ``monoids'' and ``modules'' in
a symmetric monoidal category and apply it to the category of symmetric
sequences. In this language, the symmetric sequence of spheres
$S=(S^0,S^1,\dots,S^n,\dots)$ is a commutative monoid in the
category of symmetric sequences and a symmetric spectrum is a (left)
$S$-module.

Consider the symmetric sphere spectrum $S$.  By
Proposition~\ref{prop-bilinear-hom}, the natural
$\Sigma_p\times\Sigma_q$-equivariant maps $m_{p,q}\mathcolon S^p\sm S^q\to
S^{p+q}$ give a pairing $m\mathcolon S\otimes S\to S$. The adjoint $G_0S^0\to
S$ of the identity map $S^0\to \Ev_0S=S^0$ is a two sided
unit of the pairing.  The diagram of natural isomorphisms
\[
\begin{CD}
S^p\sm S^q\sm S^r @>>> S^p\sm S^{q+r}\\
@VVV    @VVV\\
S^{p+q}\sm S^r @>>> S^{p+q+r}
\end{CD}
\]
commutes, showing that $m$ is an associative pairing of symmetric
sequences.

A symmetric spectrum $X$ has an underlying symmetric sequence of pointed
simplicial sets $\{X_{n} \}$ and a collection of pointed
$\Sigma_p\times\Sigma_q$-equivariant maps $\sigma^p\mathcolon S^p\sm
X_q\to X_{p+q}$ for $p\ge1$ and $q\ge0$. Let $\sigma^0\mathcolon S^0\sm
X_n\to X_n$ be the unit isomorphism. Then the diagram
\[
\xymatrix{ S^p\sm S^q\sm X_r \ar[r]^-{S^p\sm \sigma^q}\ar[d]_{\natiso}
& S^p\sm X_{q+r} \ar[d]^{\sigma^p}\\ S^{p+q}\sm X_r
\ar[r]^-{\sigma^{p+q}} & X_{p+q+r} }
\]
commutes for $p,q,r\ge0$. By Proposition~\ref{prop-bilinear-hom},
the equivariant maps $\sigma^p$ give a pairing $\sigma\mathcolon S\otimes
X\to X$ such that the composition $G_0(S^0)\otimes X\to S\otimes X\to X$
is the unit isomorphism and the diagram
\[
\xymatrix{
S\otimes S\otimes X \ar[r]^-{S\otimes\sigma}\ar[d]_{m\otimes X} & S\otimes X
\ar[d]^{\sigma}\\ 
S\otimes X \ar[r]^-{\sigma} & X
}
\]
commutes. 

In the language of monoidal categories, $S$ is a monoid in the
category of symmetric sequences and a symmetric spectrum is a left
$S$-module. Moreover, $S$ is a commutative monoid, \ie the diagram
\[
\xymatrix{{S\otimes S}\ar [rr]^{\tau}\ar [dr]_m&&{S\otimes S}\ar[dl]^m\\
&S&},
\]
commutes, where $\tau$ is the twist isomorphism.  To see this, one
can use either the definition of the twist isomorphism or the
description given in Remark~\ref{rem-twist-isomorphism}.  Then,
as is the case for commutative monoids in the category of sets and for
commutative monoids in the category of abelian groups
(\ie commutative rings), there is a tensor product $\otimes_S$, having
$S$ as the unit.  This gives a symmetric monoidal product on the category
of $S$-modules.  The smash product $X\sm Y$ of $X,Y\in\spec$ is the
symmetric spectrum $X\otimes_S Y$.

We review the necessary background on monoidal categories.  Monoids and
modules can be defined in any symmetric monoidal category, see
\cite{maclane-categories}. Let $\otimes$ be a symmetric monoidal product
on a category $\cc$ with unit $e\in\cc$.  A \emph{monoid} in $\cc$ is an
object $R\in\cc$, a multiplication $\mu \mathcolon R\otimes R \to R$,
and a unit map $\eta\mathcolon e\to R$ such that the diagram
\[
\begin{CD}
R\otimes R \otimes R @>{m\otimes R}>> R \otimes R\\
@V{R\otimes m}VV  @VVmV\\
R\otimes R @>>m>  R
\end{CD}
\]
commutes (\ie $m$ is associative) and such that the compositions
$e\otimes R \xrightarrow{\eta \otimes 1} R\otimes R\xrightarrow{\mu} R$
and $R\otimes e \xrightarrow{\eta \otimes 1} R\otimes R\xrightarrow{\mu}
R$ are the unit isomorphisms of the product $\otimes$.  The monoid $R$
is \emph{commutative} if $\mu = \mu \circ \tau$ where $\tau$ is the twist
isomorphism of $\otimes$. A left \emph{R-module} is an object $M$ with an
associative multiplication $\alpha\mathcolon R\otimes M \to M$ that
respects the unit; right $R$-modules are left $R^{\textup{op}}$-modules,
where $R^{\textup{op}}$ is $R$ with the multiplication $\mu \circ \tau
$. A map of left $R$-modules $M$ and $N$ is a map $f\mathcolon M\to N$
in $\cc$ that commutes with the left $R$-actions. The category of left
$R$-modules is denoted $R\text{-}Mod$.  If $\cc$ is complete, the module
category $R\text{-}Mod$ is complete and the forgetful functor
$R\text{-}Mod\to\cc$ preserves limits. If $\cc$ is cocomplete and the
functor $R\otimes(-)$ preserves coequalizers (in fact it suffices that
$R\otimes(-)$ preserve reflexive coequalizers) then the module category
$R\text{-}Mod$ is cocomplete and the forgetful functor
$R\text{-}Mod\to\cc$ preserves colimits.

The symmetric sequence of spheres $S=(S^0,S^1,S^2,\dots)$ is a
commutative monoid in the category of symmetric sequences.
 
\begin{proposition}\label{prop-S-mod=sym-spectra}
The category of symmetric spectra is naturally equivalent to
the category of left $S$-modules. 
\end{proposition}

\begin{proof}
A pairing $m\mathcolon S\otimes X\to X$ is the same as a collection of
$\Sigma_p\times\Sigma_q$-equivariant maps $m_{p,q}\mathcolon S^p\sm
X_q\to X_{p+q}$. If $X$ is a left $S$-module, there is a spectrum for
which $X$ is the underlying symmetric sequence and the structure maps
are the maps $\sigma=m_{1,n}\mathcolon S^1 \wedge X_n \to X_{n+1}$. The
compositions $\sigma^p$ are the $\Sigma_p\times\Sigma_q$-equivariant
maps $m_{p,q}$. Conversely, for $X$ a symmetric spectrum, the map of
symmetric sequences $m\mathcolon S\otimes X\to X$ corresponding to the
collection of $\Sigma_p \times \Sigma_q$-equivariant maps
$m_{p,q}=\sigma^p \mathcolon S^p \wedge X_q \to X_{p+q}$, where
$\sigma^0$ is the natural isomorphism $S^0\sm X_n\to X_n$, makes $X$ a
left $S$-module. These are inverse constructions and give a natural
equivalence of categories.
\end{proof}

The smash product on the category of symmetric spectra is a special
case of the following lemma.

\begin{lemma} \label{lem-R-mod}
Let $\cc$ be a symmetric monoidal category that is cocomplete and let
$R$ be a commutative monoid in $\cc$ such that the functor
$R\otimes(-)\mathcolon \cc\to\cc$ preserves coequalizers. Then there is
a symmetric monoidal product $\otimes_R$ on the category of $R$-modules
with $R$ as the unit.
\end{lemma}

\begin{proof} 
Let $Q$ be a monoid in $\cc$, $M$ be a right $Q$-module, and $N$ be a
left $Q$-module.  The tensor product $M\otimes_Q N$ is the colimit in
$\cc$ of the diagram
\[
\xymatrix{ M\otimes Q\otimes N \ar@<\sep>[r]^-{m\otimes 1}
\ar@<-\sep>[r]_-{ 1\otimes m} & M\otimes N }.
\]
If $M$ is an $(R,Q)$-bimodule (\ie the right action of $Q$ commutes with
the left action of $R$) then this is a diagram of left $R$-modules.
Since the functor $R\otimes(-)$ preserves coequalizers, the colimit
$M\otimes _Q N$ is a left $R$-module.  If $M$ is a left module over the
commutative monoid $R$, the composition $M\otimes R\xrightarrow{\tau}
R\otimes M \xrightarrow{\alpha} M$ is a right action of $R$; since $R$
is commutative, the two actions commute and $M$ is an
$(R,R)$-bimodule. Hence, the tensor product $M\otimes_R N$ is a left
$R$-module.  The unit of the product $\otimes_R$ is $R$ as a left
$R$-module. The associativity, unit, and twist isomorphisms of the
product on $\cc$ induce corresponding isomorphisms for the product
$\otimes_R$ on the category of left $R$-modules. Thus $\otimes_R$ is a
symmetric monoidal product.
\end{proof}

\begin{definition}
The \emph{smash product} $X\sm Y$ of symmetric spectra $X$ and $Y$
is the symmetric spectrum $X\otimes_S Y$.
\end{definition}

Apply Lemma~\ref{lem-R-mod} to the commutative monoid $S$ in the
bicomplete category of symmetric sequences $\symseq$ to obtain the
following corollary.

\begin{corollary}\label{cor-S-mod}
The smash product $X\sm Y$ is a
symmetric monoidal product on the category of symmetric spectra.
\end{corollary}

Next, some important functors on the category of symmetric spectra.

\begin{definition} \label{def-spectra-evaluation}
The functor $S\otimes(-)\mathcolon\symseq\to\spec$ gives the \emph{free
$S$-module} $S\otimes X$ generated by the symmetric sequence $X$.  For
each $n\ge 0$, the \emph{evaluation functor} $\Ev_{n}\mathcolon
\spec\to\sset$ is given by $\Ev_nX=X_n$ and $\Ev_nf=f_n$. The
\emph{free functor} $F_{n}\mathcolon\sset \to\spec $ is the left adjoint
of the evaluation functor $\Ev_n$. The functor
$R_n\mathcolon\sset\to\spec$ is the right adjoint of the evaluation
functor $\Ev_n\mathcolon\spec\to\sset$.  
\end{definition}

The functor $S\otimes(-)$ is left adjoint to the forgetful functor
$\spec\to\symseq$. The free functor $F_n$ is the composition $S\otimes
G_n$ of the left adjoints $G_n\mathcolon\sset\to\symseq$
(Definition~\ref{def-sequence-evaluation}) and
$S\otimes(-)\mathcolon\symseq\to\spec$. Thus, for $X\in\spec$ and
$K\in\sset$, the left $S$-module $X\sm F_nK$ is naturally isomorphic to
the left $S$-module $X\otimes G_nK$. In particular, $X\sm F_0K$ is
naturally isomorphic to the symmetric spectrum $X\sm K$ defined by
prolongation in Section~\ref{subsec-simplicial-structure}.  Furthermore
$F_0K=S\sm K $ is the symmetric suspension spectrum $\Sigma^\infty K$ of
$K$, and $F_0S^0$ is the symmetric sphere spectrum $S$.  For a pointed
simplicial set $K$, $R_nK$ is the symmetric sequence
$\ihom_{\symseq}(S,K^{\Sigma(-,\ov n)_+})$, which is a left $S$-module
since $S$ is a right $S$-module.

We leave the proof of the following proposition to the reader.

\begin{proposition}\label{prop-adjunctions-for-spectra}
There are natural isomorphisms\textup{:}
\begin{enumerate}
\item $F_{m}(K)\sm F_{n}(L)\natiso F_{m+n}(K\sm L)$ for
$K,L\in\sset$.
\item $\shom_{\spec}(S\otimes X,Y)\natiso\shom_{\symseq} (X,Y)$ for
$X\in\symseq$ and $Y\in\spec$.
\item $\shom_{\spec}(F_nK,X)\natiso \shom_{\sset}(K,\Ev_nX)$ for
$K\in\sset$ and $X\in\spec$.
\end{enumerate}
\end{proposition}

\begin{proposition}
\label{prop-F-preserves}
Let $f$ be a map of pointed simplicial sets.
\begin{enumerate}
\item $F_n\mathcolon\sset\to\spec$ preserves colimits. 
\item If $f$ is a monomorphism then $F_nf$ is a  monomorphism.
\item If $f$ is a weak equivalence then $F_nf$ is a level equivalence.
\end{enumerate}
\end{proposition}

\begin{proof}
Use the isomorphism $F_nf=S\otimes G_nf$ and
Proposition~\ref{prop-otimes-preserves}. 
\end{proof}

The internal Hom on the category of symmetric spectra is a special
case of the following lemma.

\begin{lemma} \label{lem-R-Hom}
Let $\cc$ be a closed symmetric monoidal category that is bicomplete
and let $R$ be a commutative monoid in $\cc$. Then there is a function
$R$-module $\ihom_R(M,N)$, natural for $M,N\in\cc$, such that the
functor $(-)\otimes_R M$ is left adjoint to the functor
$\ihom_R(M,-)$
\end{lemma}

\begin{proof} 
Let $R$ be a monoid in $\cc$ and let $M$ and $N$ be left
$R$-modules.  Then $\ihom_R(M,N)$ is the limit in $\cc$  of the
diagram
\[
\xymatrix{
{\ihom_{\cc}(M,N)}\ar@<\sep>[r]^-{m^*}\ar@<-\sep>[r]_-{m_*} &
{\ihom_{\cc}(R\otimes M,N)} 
}
\]
where $m^*$ is pullback along the multiplication $m\mathcolon R\otimes
M\to M$ and $m_*$ is the composition
\[
\ihom_{\cc}(M,N)\xrightarrow{R\otimes -}\ihom_{\cc}(R\otimes M,R\otimes
N)\xrightarrow{m_*}\ihom_{\cc}(R\otimes M,N)
\]
If $Q$ is another monoid and $N$ is an $(R,Q)$-bimodule then this is a
diagram of right $Q$-modules and the limit $\ihom_R(M,N)$ is a right
$Q$-module.  If $R$ is a commutative monoid and $N$ is a left $R$-module
then the left action is also a right action and $N$ is an
$(R,R)$-bimodule. So $\ihom_R(M,N)$ is a right $R$-module and hence a
left $R$-module.  It follows from the properties of the internal Hom in
$\cc$ and the definition of $\otimes_R$ that $(-)\otimes_R M$ is left
adjoint to $\ihom_R(M,-)$.
\end{proof}

\begin{definition}\label{def-function-spectrum}
Let $X$ and $Y$ be symmetric spectra.  The \emph{function spectrum}
$\Hom_S(X,Y)$ is the limit of the diagram in $\spec$
\[
\xymatrix{
{\ihom_{\Sigma }(X,Y)}\ar@<\sep>[r]^-{m^*}
\ar@<-\sep>[r]_-{m_*} & {\ihom_{\Sigma }(S\otimes X,Y)}
}.
\]
\end{definition}

Combining Lemmas \ref{lem-R-mod} and \ref{lem-R-Hom}:

\begin{theorem}\label{thm-monoidal-smash} 
The smash product is a closed symmetric monoidal product on the category
of symmetric spectra.  In particular, there is a natural adjunction
isomorphism
\[
\spec(X \sm Y, Z) \natiso\spec(X, \ihom_S(Y,Z)).
\]
\end{theorem}

\begin{proof}
The smash product $\sm$ is a symmetric monoidal product by
Corollary \ref{cor-S-mod}.  The adjunction isomorphism follows from
Lemma~\ref{lem-R-Hom}
\end{proof}

The adjunction is also a simplicial adjunction and an internal
adjunction.

\begin{corollary}\label{cor-spectra-adjunctions}
There are natural isomorphisms
\[
\shom_{\spec}(X \sm Y, Z) \natiso\shom_{\spec}(X, \ihom_S(Y,Z)).
\]
and
\[
\ihom_S(X \sm Y, Z) \natiso\ihom_S(X, \ihom_S(Y,Z))
\]
\end{corollary}

\begin{remark}\label{rem-describe-internal-hom}
We use Proposition~\ref{prop-adjunctions-for-spectra} to give another
description of the function spectrum $\ihom_S(X,Y)$. For a symmetric
spectrum $X$, the pointed simplicial set of maps $\shom_{\spec
}(F_nS^0,X)$ is naturally isomorphic to $\shom_{\sset
}(S^0,\Ev_nX)=X_n$. The symmetric spectrum $F_nS^0$ is the $S$-module
$S\otimes \Sigma[n]_+$ and as $n$ varies, $S\otimes \Sigma[-]_+$ is a
functor $\Sigma^{\textup{op}}\to\spec$. The symmetric sequence
$\shom_{\spec }(S\otimes\Sigma[-]_+ ,X)$ is the underlying symmetric
sequence of $X$.  In particular, the natural isomorphism
$X_n=\shom_{\spec }(F_nS^0,X)$ is $\Sigma_n$-equivariant.  Applying this
to $\ihom _{S}(X,Y)$ and using Corollary~\ref{cor-spectra-adjunctions},
we find that the underlying symmetric sequence of $\ihom _{S}(X,Y)$ is
the symmetric sequence $\shom _{\spec }(X\sm (S\otimes \Sigma
[-]_{+}),Y)$.  

We must also describe the structure maps of $X$ from this point of view.
Recall that $\shom_{\spec }(F_nS^0,X)=X_n$, $\shom_{\spec
}(F_nS^1,X)=\shom_{\sset }(S^1,X_n)$.  Let
$\badmap \mathcolon F_1S^1\to F_0S^0$ be the adjoint of the identity map
$S^1\to\Ev_1 F_0S^0=S^1$. The induced map $\shom_{\spec
}(\badmap,X)\mathcolon X_0\to\shom_{\sset }(S^1,X_1)$ is adjoint to the
structure map $S^1\sm X_0\to X_1$. The map
\[
\badmap \sm F_nS^0\mathcolon F_1S^1\sm F_nS^0=F_{n+1}S^1\to F_0S^0\sm
F_nS^0=F_nS^0
\]
is $\Sigma_1\times\Sigma_n$-equivariant; the induced map 
\[
\shom_{\spec }(\badmap\sm F_nS^0,X)\mathcolon X_n\to\shom_{\sset
}(S^1,X_{n+1})
\]
is $\Sigma_1\times\Sigma_n$-equivariant and is adjoint to the structure
map $\sigma\mathcolon S^1\sm X_n\to X_{n+1}$.  In order to apply this to
$\ihom _{S}(X,Y)$, use Proposition~\ref{prop-adjunctions-for-spectra}
and Corollary~\ref{cor-spectra-adjunctions} to find a natural isomorphism
\[
\shom_{\spec }(X\sm F_{n+1}S^1, Y)\cong\shom_{\sset }(S^1,\shom_{\spec
}(X\sm F_{n+1}S^0,Y).
\]
Using this natural isomorphism, we find that the structure maps of
$\ihom_S(X,Y)$ are the adjoints of the maps
\[
\shom_{\spec }(X\sm F_nS^0,Y)\to\shom_{\spec }(X\sm F_{n+1}S^1,Y)
\]
induced by $\badmap \sm F_nS^0$.   

For example, $\Hom_S(F_kS^0,X)$ is the \emph{$k$-shifted} spectrum; its
underlying symmetric sequence is the sequence of pointed simplicial sets
\[
X_k,X_{1+k},\dots,X_{n+k}\dots
\]
with $\Sigma_n$ acting on $X_{n+k}$ by restricting the action of
$\Sigma_{n+k}$ to the copy of $\Sigma_n$ that permutes the first $n$
elements of $\ov{n+k}$.  The structure maps of the $k$-shifted spectrum
are the structure maps $\sigma\mathcolon S^1\sm X_{n+k}\to X_{n+k+1}$ of
$X$.  More generally, $\Hom_S(F_kK,X)$ is the $k$-shifted spectrum of
$X^K$.
\end{remark}

\subsection{The ordinary category of spectra}\label{subsec-BF}

An approach similar to the last two sections can be used to describe
(non-symmetric) spectra as modules over the sphere spectrum in a
symmetric monoidal category.  But in this case the sphere spectrum is
not a commutative monoid, which is why there is no closed
symmetric monoidal smash product of spectra.

We begin as in Section~\ref{subsec-sequences} by considering a
category of sequences.

\begin{definition}\label{def-sequences}
The category $\cN$ is the category with the non-negative
integers as its objects and with the identity maps of the objects as
its only maps. 
The category of \emph{sequences} $\seq $ is the category of functors
from $\cN$ to $\sset$.  An object of $\seq $ is a sequence
$X_0,X_1,\dots, X_n,\dots$ of pointed simplicial sets and a map
$f\mathcolon X\xrightarrow{}Y$ is a sequence of pointed simplicial maps
$f_{n}\mathcolon X_{n}\to Y_{n}$.
\end{definition}

\begin{definition}\label{def-graded-smash}
The \emph{graded smash product} of sequences $X$ and $Y$ is the
sequence $X\otimes Y$ given in degree $n$ by
\[
(X\otimes Y)_{n}=\bigvee_{p+q=n} X_{p}\wedge Y_{q}
\]
\end{definition}

\begin{lemma}\label{lem-sequences-closed-symmetric-monoidal} 
The category of sequences is a bicomplete category and the graded
smash product is a symmetric monoidal product on $\seq$.
\end{lemma}

\begin{proof}
Limits and colimits are defined at each object, so $\seq $ has
arbitrary limits and colimits since $\sset$ does.  The graded tensor
product is a symmetric monoidal product on $\sset^\cN$. The unit is
the sequence $S^0,*,*,\cdots$\,.  The coherence isomorphisms for the
graded smash product follow from the corresponding coherence
isomorphisms in $\sset$.  The twist isomorphism takes $(x,y)$ to
$(y,x)$.
\end{proof}

\begin{proposition}\label{prop-spectra-are-modules}
The sequence $S$ whose $n$th level is $S^n$ is a monoid in the
category of sequences.  The category of left $S$-modules is isomorphic to the
ordinary category of spectra, $\BF $.
\end{proposition}

\begin{proof}
By definition, a pairing of sequences $X\otimes Y\to Z$ is the same as
a collection of maps $X_p\sm Y_q\to Z_{p+q}$ for $p,q\ge0$. The
associative pairing $S^p\sm S^q\to S^{p+q}$ gives an associative
pairing of sequences $S\otimes S\to S$ with an obvious unit map. Thus
$S$ is a monoid. Since $S^p=S^1\sm \dots \sm S^1$, a left action of
the monoid $S$ on the sequence $X$ is the same as a collection of maps
$S^1\sm X_q\to X_{q+1}$ for $q\ge0$.
\end{proof}

Suppose there were a closed symmetric monoidal structure on $\BF$ with
$S$ as the unit; since $\BF$ is the category of left $S$-modules it
would follow that $S$ is a commutative monoid in the category of
sequences.  But the twist map on $S^1\sm S^1$ is not the identity map and
so $S$ is not a commutative monoid in $\seq$.  Therefore there is no
closed symmetric monoidal smash product of spectra.  Note that the
category of right $S$-modules is isomorphic to the category of left
$S$-modules, even though $S$ is not commutative.  On the other hand, $S$
is a graded commutative monoid up to homotopy, in fact up to
$E_\infty$-homotopy.  This observation underlies Boardman's construction
of handicrafted smash products \cite{adams-blue}. 

\section{Stable homotopy theory of symmetric spectra}\label{sec-homotopy}

To use symmetric spectra for the study of stable homotopy theory, one
should have a stable model category of symmetric spectra such that the
category obtained by inverting the stable equivalences is naturally
equivalent to Boardman's stable homotopy category of spectra (or to any
other known to be equivalent to Boardman's). In this section we define
the stable model category of symmetric spectra. In
Section~\ref{sec-comparison} we show that it is Quillen equivalent to
the stable model category of spectra discussed in
\cite{bousfield-friedlander}.

In Section~\ref{subsec-stable-equivalence} we define the class of stable
equivalences of symmetric spectra and discuss its non-trivial
relationship to the class of stable equivalences of (non-symmetric)
spectra. In Section~\ref{subsec-model-cat} we recall the axioms and
basic theory of model categories.  In Section~\ref{subsec-level} we
discuss the level structure in $\spec $, and in
Section~\ref{subsec-stable-model-cat} we define the stable model
structure on the category of symmetric spectra which has the stable
equivalences as the class of weak equivalences.  The rest of the section
is devoted to checking that the stable model structure satisfies the
axioms of a model category.

\subsection{Stable equivalence}\label{subsec-stable-equivalence}

One's first inclination is to define stable equivalence using the
forgetful functor $U\mathcolon \spec \to \BF$; one would like a map $f$
of symmetric spectra to be a stable equivalence if the underlying map
$Uf$ of spectra is a stable equivalence, \ie if $Uf$ induces an
isomorphism of stable homotopy groups. The reader is warned:
\textbf{THIS WILL NOT WORK}. Instead, stable equivalence is defined
using cohomology; a map $f$ of symmetric spectra is a stable equivalence
if the induced map $E^*f$ of cohomology groups is an isomorphism for
every generalized cohomology theory $E$.  The two alternatives, using
stable homotopy groups or using cohomology groups, give equivalent
definitions on the category of (non-symmetric) spectra but not on the
category of symmetric spectra.

It would be nice if the $0$th cohomology group of the symmetric spectrum
$X$ with coefficients in the symmetric $\Omega$-spectrum $E$ could be
defined as $\pi_0\shom_{\spec }(X,E)$, the set of simplicial homotopy
classes of maps from $X$ to $E$. But, even though the contravariant
functor $E^0=\pi_0\shom_{\spec }(-,E)$ takes simplicial homotopy
equivalences to isomorphisms, $E^0$ may not take level equivalences to
isomorphisms. This is a common occurrence in simplicial categories, but
is a problem as every level equivalence should induce an isomorphism of
cohomology groups; a level equivalence certainly induces an isomorphism
of stable homotopy groups. We introduce injective spectra as a class of
spectra $E$ for which the functor $E^0$ behaves correctly.

\begin{definition}\label{def-inj-spec}
An \emph{injective spectrum} is a symmetric spectrum $E$ that has the
extension property with respect to every monomorphism $f$ of symmetric
spectra that is a level equivalence.  That is, for every diagram in $\spec$
\[
\begin{CD}
X @>g>> E \\
@VfVV @. \\
Y @.
\end{CD}
\]
where $f$ is a monomorphism and a level equivalence there is a map
$h\mathcolon Y\to E$ such that $g=hf$.
\end{definition}

Some examples of injective spectra follow. Recall that
$R_n\mathcolon\sset\to\spec$ is the right adjoint of the evaluation functor
$\Ev_n\mathcolon\spec\to\sset$. Also recall that a Kan complex has the
extension property with respect to every map of pointed simplicial sets
that is a monomorphism and a weak equivalence.

\begin{lemma}\label{lem-function-injective-spectra}
If the pointed simplicial set $K$ is a Kan complex then $R_nK$ is an
injective spectrum. If $X$ is a symmetric sequence and $E$
is an injective spectrum then $\ihom_S(S\otimes X,E)$ is an
injective spectrum.
\end{lemma}

\begin{proof}
Since $\Ev_n$ is left adjoint to $R_n$, the spectrum $R_nK$ has
the extension property with respect to the 
monomorphism and level equivalence $f$ if and only if $K$ has the
extension property with respect to the monomorphism and weak
equivalence $\Ev_nf$.  Since $K$ is a Kan complex, it does have the extension
property with respect to $\Ev_nf$. Hence $R_nK$ is injective.

Since $(S\otimes X)\sm (-)$ is the left adjoint of $\ihom_S(S\otimes
X,-)$, the spectrum $\ihom_S(S\otimes X,E)$ has the extension property
with respect to the monomorphism and level equivalence $f$ if and only
$E$ has the extension property with respect to the map $(S\otimes X)\sm
f$.  There is a natural isomorphism of maps of symmetric sequences
$(S\otimes X)\sm f\natiso X\otimes f$. Since $f$ is a monomorphism and a
level equivalence, $X\otimes f$ is a monomorphism and a level
equivalence by Proposition~\ref{prop-otimes-preserves}. Thus $(S\otimes
X)\sm f$ is also a monomorphism and level equivalence of symmetric
spectra.  So $\ihom_S(S\otimes X, E)$ is injective.
\end{proof}

In fact, injective spectra are the fibrant objects of a model structure
on $\spec$ for which every object is cofibrant
(Section~\ref{subsec-level-model}).  In particular, we will see in
Corollary~\ref{cor-enough-injectives}, there are enough injectives;
every symmetric spectrum embeds in an injective spectrum by a map that
is a level equivalence.

\begin{definition}\label{def-stable-equivalence}
A map $f\mathcolon X\to Y$ of symmetric spectra is a \emph{stable
equivalence} if $E^0f\mathcolon E^0Y\to E^0X$ is an isomorphism
for every injective $\Omega$-spectrum $E$.
\end{definition}

There are two other ways to define stable equivalence.

\begin{proposition}\label{prop-def-stable-equivalence}
Let $f\mathcolon X \to Y$ be a map of symmetric spectra. The following
conditions are equivalent\textup{:}
\begin{itemize}
\item $E^0f$ is an isomorphism for every injective $\Omega$-spectrum
$E$;
\item $\shom_{\spec }(f,E)$ is a weak equivalence for every injective
$\Omega$-spectrum $E$;
\item $\ihom_S(f,E)$ is a level equivalence for every
injective $\Omega$-spectrum $E$.
\end{itemize}
\end{proposition}

\begin{proof}
Let $K$ be a pointed simplicial set and $E$ be a symmetric
$\Omega$-spectrum.  The adjoints of the structure maps of $E$ are weak
equivalences of Kan complexes. From
Remark~\ref{rem-describe-internal-hom}, for $k,n\ge0$
$\Ev_k\ihom_S(F_nK,E)=E_{n+k}^K$. The adjoints of the structure maps
of $\ihom_S(F_nK,E)$ are the weak equivalences of Kan complexes
$E_{n+k}^K\to E_{n+k+1}^{S^1\sm K}$ induced by the adjoints of the
structure maps of $E$.  Therefore, $\ihom_S(F_nK,E)$ is an
$\Omega$-spectrum.

Now let $E$ be an injective $\Omega$-spectrum. Using the natural
isomorphism $F_nK\natiso S\otimes G_nK$ and
Lemma~\ref{lem-function-injective-spectra}, $\ihom_S(F_nK,E)$ is an
injective spectrum. By the preceding paragraph, $\ihom_S(F_nK,E)$ is an
$\Omega$-spectrum.  Hence $E^{S^n}= \ihom_S(F_0(S^n),E)$ and the
$k$-shifted spectrum $\ihom_S(F_kS^0,E)$ are injective $\Omega$-spectra.
Given a stable equivalence $f\mathcolon X\xrightarrow{}Y$, we want to
show that $\shom _{\spec }(f,E)$ is a weak equivalence.  Since
\[
\pi_n\shom_{\spec }(f,E)=\pi_0\shom_{\spec }(f,E^{S^n})
\]
and the simplicial sets $\shom _{\spec }(Y,E)$ and $\shom _{\spec
}(X,E)$ are Kan complexes by Lemma~\ref{lem-kan-complex}, $\shom _{\spec
}(f,E)$ is a weak equivalence on the basepoint components.  We must extend
this to all components.  To do so, note that $\shom _{\spec
}(f,E)^{S^{1}}$ is a weak equivalence for any injective $\Omega
$-spectrum $E$, since the loop space only depends on the basepoint
component.  Consider the commutative diagram 
\[
\begin{CD}
\shom _{\spec }(Y,E) @>>> \shom _{\spec }(Y, \ihom
_{S}(F_{1}S^{0},E))^{S^{1}} \\
@V\shom _{\spec }(f,E)VV @VVV \\
\shom _{\spec }(X,E) @>>> \shom _{\spec }(X,\ihom
_{S}(F_{1}S^{0},E))^{S^{1}} 
\end{CD}
\]
where the horizontal maps are induced by the map $E\xrightarrow{}\ihom
_{S}(F_{1}S^{0},E)^{S^{1}}$ adjoint to the structure map of $E$.  Since
$E$ is an injective $\Omega $-spectrum, this map is a level equivalence
of injective spectra.  By Lemma~\ref{lem-injective-spectra}, it is a
simplicial homotopy equivalence.  Hence the horizontal maps in the
diagram above are weak equivalences.  Since the right-hand vertical map
is a weak equivalence, so is the left-hand vertical map $\shom _{\spec
}(f,E)$.  Thus the first two conditions in the proposition are
equivalent.  Since $\Ev_k\ihom_S(f,E)=\shom_{\spec
}(f,\ihom_S(F_kS^0,E))$, the second two conditions are equivalent.
\end{proof}

\begin{lemma}\label{lem-kan-complex}
Suppose $X$ is a symmetric spectrum and $E$ is an injective spectrum.
Then the pointed simplicial set $\shom _{\spec }(X,E)$ is a Kan
complex.  In particular, each pointed simplicial set $E_{n}$ is a Kan
complex.  
\end{lemma}

\begin{proof}
Suppose that $f\mathcolon K\xrightarrow{}L$ is a monomorphism and weak
equivalence of simplicial sets.  We must show that $\shom _{\spec
}(X,E)$ has the extension property with respect to $f$.  By adjointness,
this is equivalent to showing that $E$ has the extension property with
respect to $X\sm f$.  But $X\sm f$ is a monomorphism and level
equivalence by Corollary~\ref{cor-SM7-spec-sset} applied to the
monomorphism $*\xrightarrow{}X$ and $f$, so $E$ does have the
required extension property.
\end{proof}

The basic properties of injective spectra which are needed in the rest
of this section are stated in the following lemma.

\begin{lemma}\label{lem-injective-spectra}
Let $f\mathcolon X\to Y$ be a map of symmetric spectra.
\begin{enumerate}
\item If $E\in\spec$ is an injective spectrum and $f$ is a level
equivalence then $E^0f$ is an isomorphism of sets.
\item If $f\mathcolon X\to Y$ is a map of injective spectra, $f$ is a level
equivalence if and only if $f$ is a simplicial homotopy equivalence.
\end{enumerate}
\end{lemma}

The proof uses the following construction.

\begin{construction}[Mapping cylinder construction]\label{mapping-cylinder}
The \emph{mapping cylinder} construction for maps of symmetric spectra
is the prolongation of the reduced mapping cylinder construction for
maps of pointed simplicial sets.  Let $i_0$ and $i_1$ be the two
inclusions $\Delta[0]\to\Delta[1]$. The cylinder spectrum $Mf$ of a map
$f\mathcolon X\to Y\in\spec$ is the corner of the pushout square
\[
\begin{CD}
X\sm \Delta[0]_+=X @>f>> Y\\
@V{X\sm i_0}VV   @VVsV\\
X \sm \Delta[1]_+ @>>g> Mf 
\end{CD}
\]
Let $i=g\circ ( X \sm i_1) \mathcolon X \to Mf$. Let $r\mathcolon Mf \to
Y$ be the map on the pushout induced by the identity map on $Y$ and the
composition $X \sm \Delta[1]_+\to X \to Y$.  Then $f=ri$, $i$ is a
monomorphism, $rs=id_Y$, and there is a simplicial homotopy from $sr$ to
the identity map of $Mf$.
\end{construction}

\begin{proof}[Proof of Lemma \ref{lem-injective-spectra}] 
For part one of the lemma, let $f\mathcolon X\to Y$ be a level equivalence
and let $Mf$ be the mapping cylinder of $f$.  As above $f=ri$,
$i\mathcolon X \to Mf$ is a monomorphism, and $r\mathcolon Mf\to Y$ is a
simplicial homotopy equivalence. Then $E^0r$ is an isomorphism and, if
$E^0i$ is an isomorphism, the composition $E^0f$ is an isomorphism.
The map $i$ is a monomorphism which, by the 2-out-of-3 property, is a
level equivalence. By the extension property of $E$ with respect to
$i$, the map $E^0i$ is surjective . The inclusion of the boundary
$j\mathcolon\partial\Delta[1]\to\Delta[1]$ is a monomorphism. By
Corollary~\ref{cor-SM7-spec-sset}, the map
\[
i\boxprod j\mathcolon X\sm\Delta[1]_+\amalg_{X\sm\partial\Delta[1]_+}
Mf\sm\partial\Delta[1]_+\to Mf\sm\Delta[1]_+
\]
is a monomorphism and a level equivalence.  The extension property of
$E$ with respect to $i\boxprod j$ implies that if $g,h\mathcolon Mf\to
E$ are maps such that $gi$ and $hi$ are simplicially homotopic, then $g$
and $h$ are simplicially homotopic. So $E^0i$ is a monomorphism and
hence $E^0i$ is an isomorphism.

For the second part of the lemma, if $f$ is a simplicial homotopy
equivalence, each $f_n$ is a simplicial homotopy equivalence of
simplicial sets and so each $f_n$ is a weak equivalence. Conversely,
suppose $f\mathcolon X\to Y$ is a level equivalence of injective spectra.
By part one, $X^0f\mathcolon X^0Y\to X^0X$ is an isomorphism. The inverse
image of the equivalence class of the identity map $X\to X$ is an
equivalence class of maps $Y\to X$. Since $Y$ is injective, $Y^0f\mathcolon
Y^0Y\to Y^0X$ is an isomorphism.  Hence each map in the equivalence
class is a simplicial homotopy inverse of $f$.
\end{proof}

Restricting part one of Lemma~\ref{lem-injective-spectra} to injective 
$\Omega$-spectra gives:

\begin{corollary}\label{cor-level-is-stable}
Every level equivalence of symmetric spectra is a stable equivalence.
\end{corollary}

Next recall the definition of stable homotopy equivalence in the
category of (non-symmetric) spectra $\BF$.

\begin{definition}\label{def-stable-groups}
For each integer $k$ the \emph{$k$th homotopy group} of the
spectrum (or symmetric spectrum) $X$ is the colimit
\[
\pi_k X=\colim_n\pi_{k+n}X_n.
\]
of the directed system given by the compositions
\[
\pi_{k+n}X_n\xrightarrow{E}\pi_{k+n+1}(S^1\sm
X_n)\xrightarrow{\pi_{k+n+1}\sigma} \pi_{k+n+1}X_{n+1}
\]
of the suspension homomorphism and the map induced by $\sigma$.
\end{definition}

A map of spectra $f\in\BF$ is a \emph{stable homotopy equivalence} if
$\pi_*f$ is an isomorphism.  For example, every level equivalence in
$\BF$ is a stable homotopy equivalence as it induces an isomorphism of
homotopy groups. We do not define stable equivalence of symmetric
spectra in this way; as the following example shows, a stable
equivalence of symmetric spectra need not induce an isomorphism of
homotopy groups.

\begin{example}\label{eg-stable-equiv}
The map $\badmap \mathcolon F_{1}S^{1}\to F_{0}S^{0}$ (see
\ref{rem-describe-internal-hom}) is the adjoint of the identity map
$S^1\to \Ev_1S=S^1$. The $n$th space of $F_{1}S^{1}$ is
$(\Sigma_{n})_{+}\wedge_{\Sigma _{n-1}}S^n$, which is a wedge of $n$
copies of $S^{n}$.  One can calculate that $\pi_0F_{1}S^{1}$ is an
infinite direct sum of copies of the integers $\mathbb{Z}$, whereas
$\pi_0F_{0}S^{0}$ is $\mathbb{Z}$. So $\pi_*\badmap $ is not an
isomorphism of homotopy groups and thus $U\badmap $ is not a stable
homotopy equivalence of (non-symmetric) spectra.  But $\badmap $ is a
stable equivalence of symmetric spectra. For a symmetric
$\Omega$-spectrum $E$, $\shom_{\spec}(F_1S^1,E)=\shom_{\sset}(S^1,E_1)$,
$\shom_{\spec}(F_0S^0,E)=E_0$, and the induced map
$\shom_{\spec}(\badmap,E)=E_0 \to E_1^{S^1}$ is adjoint to the structure
map $S^1\sm E_0\to E_1$. So $\shom_{\spec }(\badmap,E)$ is a weak
equivalence for every $\Omega$-spectrum $E$, including the injective
ones, and so $\badmap$ is a stable equivalence. By the same argument,
the maps $\badmap\sm F_nS^0$ are stable equivalences.
\end{example}

The forgetful functor $U\mathcolon\spec\to\BF$ does not preserve stable
equivalences.  On the other hand, the functor $U$ does reflect stable
equivalences.

\begin{theorem}\label{thm-U-detects-stable}
Let $f$ be a map of symmetric spectra such that $\pi_*f$ is an
isomorphism of homotopy groups.  Then $f$ is a stable equivalence.
\end{theorem}

\begin{proof}
To ease notation, let $RX=\Hom_S(F_1S^1, X)$, so that
$R^nX=\Hom_S(F_nS^n,X)$ for $n\ge0$ and $\Ev_kR^nX=\shom_{\sset
}(S^n,X_{n+k})$. In particular $R^0X=X$ and there is a natural map
$\badmap^* \mathcolon X\to RX$ induced by the map $\badmap \mathcolon
F_{1}S^{1}\xrightarrow{} F_{0}S^{0}$ discussed in
Example~\ref{eg-stable-equiv} and
Remark~\ref{rem-describe-internal-hom}.  The maps $R^n(\badmap
^{*})\mathcolon R^nX\to R^{n+1}X$ give a directed system. Let
\[
R^{\infty}X=\colim_{n\geq 0} R^n X= \colim_{n\geq 0} \Hom_S(F_nS^n,X)
\]
and let $r_X\mathcolon X\to R^{\infty}X$ be the natural map from $X$ to the
colimit.  

Let $E$ be an injective $\Omega$-spectrum.  Since $E$ is an
$\Omega$-spectrum, the map 
\[
R^n(\badmap ^{*}) \mathcolon R^nE\to R^{n+1}E
\]
is a level equivalence for each $n$ and the map $r_E\mathcolon E\to
R^{\infty}E$ is a level equivalence.  Since $E$ is injective, the proof of 
part one of Lemma~\ref{lem-injective-spectra} applied to $r_E$ shows
that there is a map $g\mathcolon R^\infty E\to E$ such that the
composition $gr_E$ is simplicially homotopic to the identity map on $E$,
though the other composite $r_E g$ may not be simplicially homotopic to
the identity map on $R^{\infty }E$. There is a natural transformation
$E^0X\to E^0(R^\infty X)$ sending the map $f\mathcolon X\to E$ to the
map $g R^\infty f$; there is a natural transformation $E^0(R^\infty
X)\to E^0X$ induced by composition with the map $r_X\mathcolon X\to
R^{\infty}X$.  Since $gr_E$ is simplicially homotopic to the identity
map on $E$, the composition of the natural transformations is the
identity natural transformation of the functor $E^0$. In the diagram
\[
\xymatrix{
E^0Y \ar[d]_{E^0f}\ar[r] & {E^0R^\infty Y} \ar[d]_{E^0R^\infty f}\ar[r]
& E^0Y\ar[d]_{E^0f}\\
E^0X \ar[r] & {E^0R^\infty X} \ar[r] & E^0X
}
\]
the composition of the horizontal maps is the identity, showing that
$E^0f$ is a retract of $E^0R^\infty f$. 

Let $X$ be a symmetric spectrum that is a level Kan complex, \ie each
$X_n$ is a Kan complex. Since the functor $\pi_k$ commutes with filtered
colimits, the homotopy group $\pi_k\Ev_nR^{\infty}X$ of the pointed
simplicial set $\Ev_nR^{\infty}X$ is naturally isomorphic to the
homotopy group $\pi_{k-n}X$ of the symmetric spectrum $X$.  A warning:
even though the groups $\pi_k\Ev_nR^{\infty}X$ and
$\pi_{k+1}\Ev_{n+1}R^{\infty}X$ are abstractly isomorphic, the structure
map of the symmetric spectrum $R^{\infty}X$ need not induce an
isomorphism between them. In particular, despite its similarity to the
standard construction of the $\Omega$-spectrum associated to a
(non-symmetric) spectrum, $R^{\infty}X$ need not be an $\Omega$-spectrum
and $r_X\mathcolon X \to R^{\infty}X$ need not induce an isomorphism of
homotopy groups.

Now, let $f$ be a map of symmetric spectra such that $\pi_*f$ is an
isomorphism. Assume as well that $X$ and $Y$ are level Kan complexes.
Then $R^\infty f$ is a level equivalence.  By Proposition
\ref{cor-level-is-stable}, $E^0R^\infty f$ is an isomorphism for every
injective $\Omega$-spectrum $E$.  Thus $E^0f$, which is a
retract of $E^0R^\infty f$, is an isomorphism for every injective
$\Omega$-spectrum $E$, and so $f$ is a stable equivalence.

To finish the proof, let $f\mathcolon X\to Y$ be a map of arbitrary
symmetric spectra such that $\pi_*f$ is an isomorphism.  For every
simplicial set $X$ there is a natural weak equivalence $X\to KX$ where
$KX$ is a Kan complex. There are several such functors: $K$ can be Kan's
$\text{Ex}^\infty$ functor; $K$ can be the total singular complex of the
geometric realization; or $K$ can be constructed using a simplicial
small object argument. In each case, $K$ is a $\sset$-functor and $X\to
KX$ is a $\sset$-natural transformation.  By prolongation, for every
symmetric spectrum $X$ there is a natural level equivalence $X\to KX$
where $KX$ is a level Kan complex. In the commutative diagram
\[
\begin{CD}
X @>>> KX\\
@VfVV  @VVKfV\\
Y @>>> KY
\end{CD}
\]
the horizontal maps are stable equivalences by
\ref{cor-level-is-stable}; the map $Kf$ is a stable equivalence by the
preceding paragraph; hence $f$ is a stable equivalence.
\end{proof}

As a corollary some of the standard results about spectra translate
into results about symmetric spectra.

\begin{definition}
Let $(X,A)$ be a pair of symmetric spectra where $A$ is a subspectrum of $X$. 
The \emph{$k$th relative homotopy group} of the pair $(X,A)$ is the
colimit
\[
\pi_k(X,A)=\colim_n\pi_{k+n}(X_n,A_n;*).
\]
of the relative homotopy groups of the pointed pairs of simplicial
sets $(X_n,A_n;*)$.
\end{definition}

\begin{lemma}[Stable excision]\label{lem-stable-excision}
Let $(X,A)$ be a pair of symmetric spectra with $A$ a subspectrum of
$X$.  The map of homotopy groups $\pi_k(X,A)\to \pi_k(X/A)$ is an
isomorphism.
\end{lemma}

\begin{proof}
Consider the diagram:
\[
\xymatrix{
{\pi_{k+q}(X_q,A_q)}\ar[r]\ar[d] & {\pi_{k+q+p}(S^p\sm X_q,S^p\sm A_q)}
\ar[r]\ar[d] & {\pi_{k+q+p}(X_{p+q},A_{p+q})}\ar[d]\\
{\pi_{k+q}(X_q/A_q)}\ar[r] & {\pi_{k+q+p}(S^p\sm X_q/A_q)}
\ar[r] & {\pi_{k+q+p}(X_{p+q}/A_{p+q})}
}.
\]
Let $(K,L)$ be a pair of pointed simplicial sets. By the homotopy
excision theorem the map $\pi_n(S^p\sm K,S^p\sm L)\to\pi_n(S^p\sm K/L)$
is an isomorphism when $n<2p$. So the middle vertical arrow in the
diagram is an isomorphism when $p>k+q$ and hence the map of colimits
$\pi_k(X,A)\to \pi_k(X/A)$ is an isomorphism.
\end{proof}

\begin{theorem} \label{thm-traditional-results}
\begin{enumerate}
\item Let $f\mathcolon X\to B$ be a map of symmetric spectra such that
$f_n\mathcolon X_n\to B_n$ is a Kan fibration for each $n\ge0$ and let
$F$ be the fiber over the basepoint. Then the map $X/F\to B$ is a stable
equivalence.
\item A map $f\in\spec$ of symmetric spectra is a stable equivalence
if and only if its suspension $f\sm S^1$ is a stable equivalence.
\item For symmetric spectra $X$ and $Y$ such that $Y$ is a level Kan
complex, a map $X\to Y^{S^1}$ is a stable equivalence if and only its
adjoint $X\sm S^1\to Y$ is a stable equivalence.
\end{enumerate}
\end{theorem}

\begin{proof}
By stable excision, the map $X/F\to B$ induces an 
isomorphism of homotopy groups and hence is a stable equivalence by
Theorem \ref{thm-U-detects-stable}. 

For part two, let $E$ be an injective $\Omega$-spectrum. By
Lemma~\ref{lem-function-injective-spectra}, the spectra $E^{S^1}$,
$\ihom_S(F_1S^0,E)$ and $\ihom_S(F_1S^1,E)$ are injective
$\Omega$-spectra.  If $f$ is a stable equivalence of symmetric spectra
then the map $\shom_{\spec }(f,E^{S^1})=\shom_{\spec }(f\sm S^1,E)$ is a
weak equivalence of simplicial sets and so $f\sm S^1$ is a stable
equivalence. Conversely if $f\sm S^1$ is a stable equivalence then the
map
\[
\shom_{\spec }(f\sm S^1,\ihom_S(F_1S^0,E))=\shom_{\spec
}(f,\ihom_S(F_1S^1,E))
\]
is a weak equivalence.  The map $E\to\ihom_S(F_1S^1,E)$ is a level
equivalence of injective spectra and thus a simplicial homotopy
equivalence, by Lemma~\ref{lem-injective-spectra}. So, for every
symmetric spectrum $X$, the induced map 
\[
\shom_{\spec }(X,E)\to\shom_{\spec }(X,\ihom_S(F_1S^1,E))
\]
is a simplicial homotopy equivalence.  Therefore $\shom_{\spec }(f,E)$
is a weak equivalence, and so $f$ is a stable equivalence.

For part three, let $f\mathcolon X\to Y^{S^1}$ be a map and
$f^a\mathcolon X\sm S^1 \to Y$ be the adjoint of $f$.  The diagram
\[
\xymatrix{ X\sm S^1 \ar[rr]^{f\sm S^1} \ar[rd]_{f^a} && Y^{S^1}\sm S^1
\ar[ld]^{\textup{ev}}\\ & Y
}
\]
is commutative where the map $\textup{ev}$ is the evaluation map. By
part one applied to the prolongation of the path fibration, the map
$\textup{ev}$ is a stable equivalence; by part two, $f$ is a stable
equivalence if and only if $f\sm S^1$ is a stable equivalence. Therefore
$f$ is a stable equivalence if and only if $f^a$ is a a stable
equivalence.
\end{proof}

Once we have the stable model category of symmetric spectra, part
three of this theorem tells us that it really is stable; \ie that the
suspension functor $-\sm S^1$ is an equivalence of model categories.

\subsection{Model categories}\label{subsec-model-cat}

In this section we recall the definition and the basic properties of
model categories; see~\cite{dwyer-spalinski}, \cite{hovey-model},
or~\cite{kan-model} for a more detailed introduction.

\begin{definition}\label{def-lifting-property}
An ordered pair $(f,g)$ of maps in the category $\cc$ has the
\emph{lifting property} if every commutative square 
\[
\xymatrix{
A \ar[r] \ar[d]_f & X \ar[d]^g \\  B \ar[r]  & Y
}
\]
in $\cc$ extends to a commutative diagram 
\[
\xymatrix{
A \ar[r] \ar[d]_f & X \ar[d]^g  \\
B \ar[r] \ar[ru]  &    Y
}.
\]
We also say that $f$ has the \emph{left lifting property} with respect
to $g$ and that $g$ has the \emph{right lifting property} with respect
to $f$. More generally, if $I$ and $J$ are classes of maps in $\cc$,
the pair $(I,J)$ has the \emph{lifting property} if every
pair $(f,g)$ with $f\in I$ and $g\in J$ has the lifting
property. We also say that $I$ has the \llp $J$ and that $J$ has the
\rlp $I$.
\end{definition}

It would be more accurate to say that the pair $(f,g)$ has the
lifting-extension property but we prefer the shorter term.

\begin{definition}
Let $f$ and $g$ be maps in a category $\cc$. The map $f$ is a
\emph{retract} of $g$ if it is a retract in the category of arrows,
\ie if there is a commutative diagram
\[
\begin{CD}
A @>>> B @>>> A \\
@VfVV  @VgVV @VfVV \\
X @>>> Y @>>> X.
\end{CD}
\]
such that the horizontal compositions are the identity maps. A class
of maps is closed under retracts if whenever $f$ is a retract of
$g$ and $g$ is in the class then $f$ is in the class.
\end{definition}

\begin{definition}\label{def-model-category}
A \emph{model category} is a category $\cm$ with three distinguished
classes of maps; the class of \emph{weak equivalences}, the class of
\emph{cofibrations}, and the class of \emph{fibrations}; that satisfy the
model category axioms below.  We call a map that is both a cofibration
and a weak equivalence a \emph{trivial cofibration}, and we call a map
that is both a fibration and a weak equivalence a \emph{trivial
fibration}. 
\begin{enumerate}
\item[M1] \label{bicomplete}
\emph{Limit axiom} The category $\cm$ is bicomplete (closed under
arbitrary small limits and colimits).
\item[M2] \label{2-out-of-3}
\emph{Two out of three axiom.} Let $f$ and $g$ be maps in
$\cm$ such that $gf$ is defined. If two of $f$, $g$ and $gf$
are weak equivalences then the third is a weak equivalence.
\item[M3] \label{retract} \emph{Retract axiom.} The class of weak
equivalences, the class of cofibrations, and the class of fibrations
are closed under retracts.
\item[M4] 
\label{lifting} \emph{Lifting axiom.}  A cofibration has the
\llp every trivial fibration.  A fibration has the \rlp every trivial
cofibration.
\item[M5] 
\label{factorization} \emph{Factorization axiom.} Every map
$f\in\cm$ has a factorization $f=pi$ where $i$ is a cofibration and
$p$ is a trivial fibration and a factorization $f=qj$ where $j$ is a
trivial cofibration and $q$ is a fibration.
\end{enumerate}
\end{definition}

Three classes of maps that satisfy axioms M2, M3, M4 and M5 are a
\emph{model structure} on the category. One should keep in mind that a
category can have more than one model structure; there can even be
distinct model structures with the same class of weak equivalences.

A bicomplete category has an initial object $\emptyset$ and a terminal
object $*$. In a model category, an object $X$ is \emph{cofibrant} if
the unique map $\emptyset\to X$ is a cofibration and an object $X$ is
\emph{fibrant} if the unique map $X\to *$ is a fibration.  A model
category is \emph{pointed} if the unique map $\emptyset\to *$ is an
isomorphism.

\begin{proposition}[The Retract Argument]\label{prop-retract}
Let $\cc$ be a category and let $f=pi$ be a factorization in $\cc$.
\begin{enumerate}
\item If $p$ has the right lifting property
with respect to $f$ then $f$ is a retract of $i$.
\item If $i$ has the left lifting property with
respect to $f$ then $f$ is a retract of $p$.
\end{enumerate}
\end{proposition}

\begin{proof}
We only prove the first part, as the second is similar.  Since $p$
has the right lifting property with respect to $f$, we have a lift
$g\mathcolon Y\xrightarrow{}Z$ in the diagram
\[
\begin{CD}
X @>i>> Z \\
@VfVV @VpVV \\
Y @= Y
\end{CD}
\]
This gives a diagram
\[
\begin{CD}
X @= X @= X\\
@VfVV @ViVV @VfVV \\
Y @>g>> Z @>p>> Y
\end{CD}
\]
where the horizontal compositions are identity maps, showing that $f$ is
a retract of $i$.
\end{proof}

The following proposition is a converse to the lifting axiom.

\begin{proposition}[Closure property]\label{prop-closure}
In a model category\textup{:} 
\begin{enumerate}
\item The cofibrations are the maps having the \llp
every trivial fibration.
\item The trivial cofibrations are the the maps having the \llp every
fibration.
\item The fibrations are the maps having the \rlp
every trivial cofibration. 
\item The trivial fibrations are the maps having the \rlp every
cofibration.
\end{enumerate}
\end{proposition}

\begin{proof}
Use the factorization axiom and the retract argument.
\end{proof}

In particular, any two of the three classes of maps in a model
category determine the third. For example, a weak equivalence is a map
that factors as a trivial cofibration composed with a trivial
fibration.

\begin{example}\label{eg-model-cats}
We recall the standard model structure on the category of simplicial
sets~\cite[II.3]{quillen-htpy}. A weak equivalence is a
map whose geometric realization is a homotopy equivalence of
CW-complexes. The cofibrations are the monomorphisms and every simplicial
set is cofibrant.  Recall, the standard $n$-simplex is
$\Delta[n]=\Delta(-,\ov n)$.  The boundary of $\Delta[n]$ is the
subfunctor $\partial\Delta[n]\subseteq\Delta[n]$ of non-surjective
maps. For $0\le i\le n$, the $i$th horn of $\Delta[n]$ is the
subfunctor $\Lambda^i[n]\subseteq\partial\Delta[n]$ of maps for which
$i$ is not in the image. Geometrically, $\Lambda^i[n]$ is obtained
from the boundary of $\Delta[n]$ by removing the $i$th face.  The
fibrations are the Kan fibrations, the maps that have the \rlp the
maps $\Lambda^i[n]\to\Delta[n]$ for $n>0$ and $0\le i\le n$; the
fibrant simplicial sets are the Kan complexes, the simplicial sets
that satisfy the Kan extension condition. A map is a trivial fibration
(a fibration and a weak equivalence) if and only if it has the \rlp
the maps $\partial\Delta[n] \to \Delta[n]$.  It follows that the
pointed weak equivalences, the pointed monomorphisms, and the pointed
(Kan) fibrations are a model structure on the category of pointed
simplicial sets.
\end{example}

When constructing a model category, the factorization axiom can be the
hardest to verify.  After some preliminary definitions,
Lemma~\ref{lem-factorization} constructs functorial factorizations in the
category of symmetric spectra. 

\begin{definition}\label{def-injectives}
Let $I$ be a class of maps in a category $\cc$.
\begin{enumerate}
\item
A map is \emph{$I$-injective} if it has the \rlp every map in $I$. The
class of $I$-injective maps is denoted $I\inj$.
\item
A map is \emph{$I$-projective} if it has the \llp every map in $I$.
The class of $I$-projective maps is denoted $I\proj$.
\item
A map is an \emph{$I$-cofibration} if it has the \llp every
$I$-injective map.  The class of $I$-cofibrations is the class
$(I\inj)\proj$ and is denoted $I\cof$.
\item 
A map is an \emph{$I$-fibration} if it has the \rlp every
$I$-projective map.  The class of $I$-fibrations is the class
$(I\proj)\inj$ and is denoted $I\fib$.
\end{enumerate}
\end{definition}

Injective and projective are dual notions; an $I$-injective map in $\cc$
is an $I$-projective map in $\cc^{\textup{op}}$; an $I$-fibration in
$\cc$ is an $I$-cofibration in $\cc^{\textup{op}}$. The class $I\inj$
and the class $I\proj$ are analogous to the orthogonal complement of a
set vectors. That analogy helps explain the following proposition, whose
proof we leave to the reader.

\begin{proposition}\label{prop-lifting-properties}
Let $I$ and $J$ be classes of maps in a category $\cc$.
\begin{enumerate}
\item
If $I\subseteq J$ then $J\inj\subseteq I\inj$ and $J\proj\subseteq
I\proj$.
\item Repeating the operations\textup{:} $I\subseteq I\cof$, $I\subseteq
I\fib$, $I\proj=(I\proj)\cof=(I\fib)\proj$, and
$I\inj=(I\inj)\fib=(I\cof)\inj$.
\item The following conditions are equivalent\textup{:}
\begin{itemize}
\item The pair $(I,J)$ has the lifting property.
\item $J\subseteq I\inj$.
\item $I\subseteq J\proj$.
\item The pair $(I\cof,J)$ has the lifting property.
\item The pair $(I,J\fib)$ has the lifting property.
\end{itemize}
\item The classes $I\inj$ and $I\proj$ are subcategories of
$\cc$ and are closed under retracts.
\item
The class $I\inj$ is closed under base change. That is, if 
\[
\begin{CD}
A @>>> B \\
@VfVV    @VgVV\\
X @>>> Y
\end{CD}
\]
is a pullback square and $g$ is an $I$-injective map then $f$ is an
$I$-injective map.
\item
The class $I\proj$ is closed under cobase change.  That is, if
\[
\begin{CD}
A @>>> B \\
@VfVV    @VgVV\\
X @>>> Y
\end{CD}
\]
is a pushout square and $f$ is an $I$-projective map then $g$ is an
$I$-projective map.
\end{enumerate}
\end{proposition}

\begin{corollary}
Let $I$ be a class of maps in a category $\cc$. The class $I\cof$ is a
subcategory of $\cc$ that is closed under retracts and cobase
change. The class $I\fib$ is a subcategory of $\cc$ that is closed
under retracts and base change.
\end{corollary}

Another useful elementary lemma about the lifting property is the
following.

\begin{lemma}\label{lem-lifting-adjoints}
Let $L\mathcolon\cc\to\cd$ be a functor that is left adjoint to the
functor $R\mathcolon\cd \to\cc$. If $I$ is a class of maps in $\cc$ and
$J$ is a class of maps in $\cd$, the pair $(I,RJ)$ has the lifting
property if and only if the pair $(LI,J)$ has the lifting property.
\end{lemma}

The next lemma is used repeatedly to construct factorizations.

\begin{lemma}[Factorization Lemma]\label{lem-factorization}
Let $I$ be a set of maps in the category $\spec$.  There is a
functorial factorization of every map of symmetric spectra as an
$I$-cofibration followed by an $I$-injective map.
\end{lemma}

The factorization lemma is proved using the transfinite small
object argument.  We begin by showing that every symmetric spectrum is
suitably small.

Recall that an ordinal is, by recursive definition, the well-ordered set
of all smaller ordinals.  In particular, we can regard an ordinal as a
category.  A cardinal is an ordinal of larger cardinality than all
smaller ordinals.  

\begin{definition}\label{defn-filtered}
Let $\gamma$ be an infinite cardinal. An ordinal $\alpha$ is
\emph{$\gamma$-filtered} if every set $A$ consisting of ordinals less
than $\alpha $ such that $\sup A=\alpha$ has cardinality greater than
$\gamma$.
\end{definition}

Every $\gamma$-filtered ordinal is a limit ordinal.  In fact, since
$\gamma$ is infinite, every $\gamma$-filtered ordinal is a limit
ordinal $\alpha$ for which there is no countable set $A$ of ordinals
less than $\alpha$ such that $\sup A=\alpha$.  The smallest
$\gamma$-filtered ordinal is the first ordinal of cardinality greater
than $\gamma$.  For example, $\omega_1$ is the smallest
$\aleph_0$-filtered ordinal. If $\gamma<\ov\gamma$ and $\alpha$ is
$\ov\gamma$-filtered, then $\alpha$ is $\gamma$-filtered.

Define the \emph{cardinality} of a spectrum $X$ to be the cardinality of
its underlying set $\amalg_n\amalg_k (Ev_nX)_k$. Then the cardinality of
$X$ is always infinite, which is convenient for the following lemma.

\begin{proposition}\label{prop-spec-smallness}
Let $X$ be a simplicial spectrum of cardinality $\gamma$.  Let $\alpha$
be a $\gamma$-filtered ordinal and let $D\mathcolon\alpha\to\spec$ be an
$\alpha$-indexed diagram of symmetric spectra. Then the natural map
\[
\colim_{\alpha}\spec(X,D)\to\spec(X,\colim_{\alpha} D)
\]
is an isomorphism.
\end{proposition}

\begin{proof}
Every symmetric spectrum has a presentation as a coequalizer
\[
\xymatrix{S\otimes S\otimes X\ar@<\sep>[r]\ar@<-\sep>[r] & S\otimes
X\ar[r] & X}
\]
of free symmetric spectra in the category $\spec$.  The symmetric
spectra $X$ and $S\otimes X$ have the same cardinality.  So 
the proposition follows once it is proved for
free symmetric spectra.  There is a natural isomorphism $\spec(S\otimes
X,Y)=\prod_p\sset ^{\Sigma_{p}}(X_p,Y_p)$.  The functors $\sset ^{\Sigma
_{p}}(X_n,-)$ have the property claimed for $\spec(X,-)$. 

This fact is the heart of the proposition.  To begin the proof of it,
suppose we have a map $f\mathcolon X_{n}\to \colim _{\alpha }D$, where
$D$ is an $\alpha $-indexed diagram of $\Sigma _{p}$-simplicial sets.
For each simplex $x$ of $X_{n}$, we can choose an ordinal
$\beta_{x}<\alpha $ and a simplex $y_{x}\in D_{\beta_x }$ such that
$f(x)$ is the image of $y_{x}$.  Because $\alpha $ is $\gamma
$-filtered, we can then find one ordinal $\beta <\alpha $ and a map
$g\mathcolon X\xrightarrow{}D_{\beta }$ factoring $f$.  The map $g$ may
not be simplicial or equivariant, but, again using the fact that $\alpha
$ is $\gamma $-filtered, we can go out far enough in the colimit so that
$g$ will be both simplicial and equivariant.  We leave the details to
the reader.

Since $\gamma$ is infinite, for every countable set $A$ of ordinals that
are strictly less than $\alpha$, the ordinal $\sup A$ is strictly less
than $\alpha$. Therefore, the countable product of functors
$\symseq(X,-)$ has the property claimed for $\spec(X,-)$ and the
proposition is proved.
\end{proof}

\begin{proof}[Proof of Lemma~\ref{lem-factorization}]
We begin by constructing a functorial factorization    
\[
X\xrightarrow{Ig}Eg\xrightarrow{Pg}Y
\]
of $g\mathcolon X\to Y$ such that $Ig$ is an $I$-cofibration
and $g=Pg\circ Ig$.  For a map $f\mathcolon b_f\to c_f$ in $I$, let
$Df$ be the set of commutative squares
\[
\begin{CD}
b_f  @>>> X\\
@VfVV   @VVgV\\
c_f @>>> Y.
\end{CD}
\]
Let 
\[
B=\amalg_{f\in I}\amalg_{Df}b_f,\ C=\amalg_{f\in I}\amalg_{Df}c_f,
\text{ and }F=\amalg_{f\in I}\amalg_{Df}f\mathcolon B\to C.
\]
By the definition of $Df$, there is a commutative square
\[
\begin{CD}
B @>>> X\\ @VFVV @VVgV\\ C @>>>Y.
\end{CD}
\]
Let $Eg$ be the pushout $X\amalg_B C$, $Ig$ be the map $X\to
Eg=X\amalg_B C$, and $Pg\mathcolon Eg\to Y$ be the natural map on the
pushout.  By construction, the map $Ig\mathcolon X\to X\amalg_B C$ is an
$I$-cofibration and $g=Pg\circ Ig$. However, the map $Pg$ need not
be an $I$-injective map.

Use transfinite induction to define functorial factorizations of $g$
\[
X\xrightarrow{I^{\alpha}g}E^\alpha g\xrightarrow{P^\alpha g}Y
\]
for every ordinal $\alpha$.  The induction starts at $0$ with $E^0g=X$,
$I^0g=id_X$, and $P^0g=g$. For a successor ordinal $\alpha +1$,
$E^{\alpha+1}g=E(P^\alpha g)$, $I^{\alpha+1}g=I(P^\alpha g)\circ
I^{\alpha}g$, and $P^{\alpha+1}g=P(P^{\alpha}g)$.  For $\beta$ a limit
ordinal, $E^\beta g=\colim_{\alpha<\beta}E^\alpha g$, $I^\beta
g=\colim_{\alpha<\beta} I^\alpha g$ and $P^\beta
g=\colim_{\alpha<\beta}P^\alpha g$.

The map $I^\alpha g\mathcolon X\to E^\alpha g$ is an $I$-cofibration for
each $\alpha$; the required lift is constructed by transfinite
induction. The proof of the lemma is completed by finding an ordinal
$\beta$ for which $P_{\beta}g$ is an $I$-injective map. Let 
\[
\begin{CD}
b_f @>>> {E^\beta g} \\ @VfVV @VV{P^\beta g}V \\ c_f @>>> Y
\end{CD}
\]
be a commutative square with $f\in I$. If the map $b_f\to E^\beta g$
factors as $b_f\to E^\alpha g\to E^\beta g$ for $\alpha<\beta$ then by
construction there is a lift $c_f\to E^{\alpha +1}g$ and, since
$\alpha+1\le\beta$, a lift $c_f\to E^{\beta}g$. Let $\gamma_f$ be the
cardinality of $b_f$ and let $\gamma=\sup_I\gamma_f$. Let $\beta$ be a
$\gamma$-filtered ordinal. Then $\spec(b_f,P^\beta g)=
\colim_{\alpha<\beta}\spec(b_f,P^\alpha g)$ for every $f\in I$ and
hence $P^\beta g$ is an $I$-injective map.
\end{proof}

\subsection{Level structure}\label{subsec-level} 

Prolongation of the model structure on $\sset$ (see \ref{eg-model-cats})
gives the level structure on the category of symmetric spectra. It is
not a model structure but it is a basic tool in the construction of the
stable model structure. Its use is already implicit in
Sections~\ref{subsec-simplicial-structure}
and~\ref{subsec-stable-equivalence}.

\begin{definition}\label{def-level-maps}
Let $f\mathcolon  X\xrightarrow{}Y$ be a map of symmetric spectra.
\begin{enumerate}
\item The map $f$ is a \emph{level equivalence} if each map 
$f_n\mathcolon X_n\to Y_n$ is a weak equivalence of simplicial sets.
\item The map $f$ is a \emph{level (trivial) cofibration} if each map
$f_{n}\mathcolon X_n\to Y_n$ is a (trivial) cofibration of simplicial sets.
\item The map $f$ is a \emph{level (trivial) fibration} if each map
$f_{n}\mathcolon X_n\to Y_n$ is a (trivial) fibration of simplicial sets.
\end{enumerate}
\end{definition}

The level cofibrations are the monomorphisms of symmetric spectra. 
Next, we characterize the level fibrations and trivial
fibrations.

\begin{definition}\label{def-level-generators}
\begin{enumerate}
\item \label{lambda-cof} Let $I_\Lambda$ denote the set of maps
$\Lambda^k[r]_+ \xrightarrow{}\Delta [r]_+$ for $r>0$  
and $0\le k\le r$. Let $FI_\Lambda= \cup_{n\geq 0}
F_n(I_\Lambda)$.
\item \label{boundary-cof} Let $I_{\partial}$ denote the set of maps
$\partial\Delta [r]_+\xrightarrow{}\Delta [r]_+$ for $r\ge0$. Let
$FI_\partial= \cup_{n\geq 0} F_n(I_\partial)$.
\end{enumerate}
\end{definition}

\begin{proposition}\label{prop-level-fib}
The level fibrations are the $FI_\Lambda$-injective maps. The level
trivial fibrations are the $FI_\partial$-injective maps.
\end{proposition}

\begin{proof}
A map $g$ is a level (trivial) fibration if and only if $\Ev_ng=g_n$
is a (trivial) Kan fibration for each $n\ge0$.  But $\Ev_ng$ is a
(trivial) Kan fibration if and only if it has the \rlp the class
($I_\partial$) $I_\Lambda$. Then by adjunction, $g$ is a level (trivial)
fibration if and only if $g$ has the \rlp the class ($FI_\partial$)
$FI_\Lambda$. 
\end{proof}

The level structure is not a model structure; it satisfies the
two-out-of-three axiom, the retract axiom, and the factorization axiom
but not the lifting axiom.  A model structure is determined by any two
of its three classes and so the level structure is over determined. In
Section~\ref{subsec-level-model} we prove there are two ``level''
model structures with the level equivalences as the weak equivalences:
one that is generated by the level equivalences and the level
cofibrations and one that is generated by the level equivalences and
the level fibrations. In any case, the level homotopy category
obtained by inverting the level equivalences is not the stable
homotopy category of spectra.

The pushout smash product (Definition~\ref{def-pushout-smash}) has an
adjoint construction.

\begin{definition}
Let $f\mathcolon U\to V$ and $g\mathcolon X\to Y$ be maps of pointed
simplicial sets. The map
\[
\shom_{\square}(f,g)\mathcolon\shom(V,X)\to
\shom(U,X)\times_{\shom(U,Y)}\shom(V,Y). 
\]
is the map to the fiber product induced by the commutative square
\[
\begin{CD}
\shom(V,X) @>f^*>> \shom(U,X)\\ @Vg_*VV @VVg_*V\\ \shom(V,Y) @>>f^*> 
\shom(U,Y).
\end{CD}
\]
Let $f$ be a map of pointed simplicial sets and $g$ be a map of
symmetric spectra. Then $\ihom_\square(f,g)$ is the map of
symmetric spectra that is defined by prolongation,
$\Ev_n\ihom_\square(f,g)=\shom_\square(f,g_n)$.
\end{definition}

\begin{proposition}\label{prop-SM7-boxhom}
\begin{enumerate}
\item If $f\in\sset$ is a monomorphism and $g\in\sset$ is a Kan
fibration then $\shom_\square(f,g)$ is a Kan fibration.  If, in
addition, either $f$ or $g$ is a weak equivalence, then $\shom
_{\square}(f,g)$ is a weak equivalence. 
\item If $f\in\sset$ is a monomorphism and $g\in\spec$ is a level
fibration then $\ihom_\square(f,g)$ is a level fibration. If, in
addition, either $f$ is a weak equivalence or $g$ is a level
equivalence, then $\ihom _{\square}(f,g)$ is a level equivalence.
\end{enumerate}
\end{proposition}

\begin{proof}
Part one is a standard property of simplicial sets, proved
in~\cite[II.3]{quillen-htpy}. Part two follows from part one by
prolongation.
\end{proof}

\begin{definition}\label{def-boxHom-category}
Let $f\mathcolon U\to V$ and $g\mathcolon X\to Y$ be maps in a category
$\cc$. Then $\cc_\square(f,g)$ is the natural map of sets to the fiber
product
\[
\cc(V,X)\to \cc(U,X)\times_{\cc(U,Y)}\cc(V,Y)
\]
coming from the commutative square
\[
\begin{CD}
\cc(V,X) @>f^*>> \cc(U,X)\\ @Vg_*VV @VVg_*V\\ \cc(V,Y) @>>f^*> \cc(U,Y)
\end{CD}
\]
\end{definition}

A pair $(f,g)$ has the lifting property if and only if
$\cc_\square(f,g)$ is surjective.

\begin{definition}
Let $f\mathcolon U\to V$ and $g\mathcolon X\to Y$ be maps of symmetric
spectra. Then $\shom_\square(f,g)$ is the natural map to the fiber product
\[
\spec_\square(f\sm\Delta[-]_+,g)\mathcolon\shom(V,X)\to 
\shom(U,X)\times_{\shom(U,Y)}\shom(V,Y)
\]
\end{definition}

\begin{proposition}\label{prop-pushout-smash-adjunctions}
Let $f$ and $h$ be maps of symmetric spectra and $g$ be a map of
pointed simplicial sets. There are natural isomorphisms
\[
\spec_\square(f\boxprod g,h)\natiso(\sset)_\square(g,\shom_\square(f,h))
\natiso\spec_\square(f,\ihom_\square(g,h))
\]
\end{proposition}

In fact this proposition holds in any simplicial model category.

\begin{proof}
Let $f\mathcolon U\to V$ and $h\mathcolon X\to Y$ be maps in $\spec$ and
$g\mathcolon K \to L$ be a map in $\sset$. Using adjunction and the
defining property of pushouts and of pullbacks, each of the three
maps in the proposition is naturally isomorphic to the map from
$\spec(V\sm L, X)$ to the limit of the diagram
\[
\xymatrix{
{\spec(U\sm L,X)}\ar[d]\ar[dr] & {\spec(V\sm K, X)}\ar[dl]\ar[dr] &
{\spec(V \sm L, Y)} \ar[d] \ar[dl] \\
{\spec(U\sm K,X)}\ar[dr] & {\spec(U\sm L, Y)}\ar[d] &
{\spec(V\sm K ,Y)}\ar[dl]\\
 & {\spec(U\sm K,Y)} 
}
\]
\end{proof}

\begin{corollary}\label{cor-boxhom-lifting}
Let $f$ and $h$ be maps of symmetric spectra and $g$ be a map of
pointed simplicial sets. 
The following are equivalent\textup{:}
\begin{itemize}
\item $(f\boxprod g,h)$ has the lifting property.
\item $(g, \shom_{\square}(f,h))$ has the lifting property.
\item $(f, \ihom_\square(g,h))$ has the lifting property.
\end{itemize}
\end{corollary}

\subsection{Stable model category}\label{subsec-stable-model-cat}

In this section we define the stable cofibrations and the stable
fibrations of symmetric spectra. The main result is
that the class of stable equivalences, the class of stable cofibrations,
and the class of stable fibrations are a model structure on $\spec$.

Recall that $f$ is a level trivial fibration if $f_n$ is a trivial Kan
fibration for each $n\ge0$.

\begin{definition}\label{def-stable-cof}
A map of symmetric spectra is a \emph{stable cofibration} if it has
the left lifting property with respect to every level trivial
fibration.  A map of symmetric spectra is a \emph{stable trivial
cofibration} if it is a stable cofibration and a stable equivalence. A
symmetric spectrum $X$ is \emph{stably cofibrant} if $*\to X$ is a
stable cofibration.
\end{definition}

The basic properties of the class of stable cofibrations are next.

\begin{proposition}\label{prop-stable-cof}
\begin{enumerate}
\item The class of stable cofibrations is the class $FI_\partial\cof$.
\item The class of stable cofibrations is a subcategory that is closed
under retracts and closed under cobase change.
\item If $f$ is a cofibration of pointed simplicial sets and $n\ge0$
then $F_nf$ is a stable cofibration. In particular, 
$F_nK$ is stably cofibrant for $K\in\sset$.
\item If $f\in\spec$ is a stable cofibration and $g\in\sset$ is a
cofibration then the pushout smash product $f\boxprod g$ is a stable
cofibration.
\item If $f\in\spec$ is a stable cofibration and $h\in\spec$ is a level
fibration then $\shom_{\square}(f,h)$ is a Kan fibration.
\end{enumerate}
\end{proposition}

\begin{proof}
The stable cofibrations are the maps having the \llp the level trivial
fibrations, which by Proposition~\ref{prop-level-fib} are the
$FI_{\partial}$-injective maps. So the stable cofibrations are the
$FI_\partial$-cofibrations.

Every class $I\cof$ has the properties stated in part two.

Suppose $g\in \spec $ is a level trivial fibration, and $f\in \sset $ is
a cofibration.  Then $f$ has the \llp the trivial Kan fibration
$\Ev_ng$. By adjunction, $F_nf$ has the \llp $g$.  Hence $F_nf$ is a
stable cofibration.  In particular, for every pointed simplicial set
$K$, the map $*\to F_nK$ is a stable cofibration, and so $F_nK$ is stably
cofibrant.

Now suppose $f\in \spec$ is a stable cofibration and $g\in \sset $ is a
cofibration.  Then, given a level trivial fibration $h\in\spec$, the map
$\ihom_\square(g,h)$ is a level trivial fibration by
Proposition~\ref{prop-SM7-boxhom}.  Therefore the pair
$(f,\ihom_\square(g,h))$ has the lifting property. Then by
Corollary~\ref{cor-boxhom-lifting}, the pair $(f\boxprod g, h)$ has the
lifting property, and so $f\boxprod g$ is a stable cofibration.

Finally, suppose $f\in \spec $ is a stable cofibration and $h\in \spec $
is a level fibration.  Given a trivial cofibration $g\in\sset$,
$\ihom_\square(g,h)$ is a level trivial fibration by
Proposition~\ref{prop-SM7-boxhom}. Therefore, the pair
$(f,\ihom_\square(g,h))$ has the lifting property.  Then, by
Corollary~\ref{cor-boxhom-lifting}, the pair $(g,\shom_\square(f, h))$
has the lifting property.  Therefore $\shom_\square(f,h)$ is a Kan
fibration.
\end{proof}

The next definition is natural in view of the closure properties in a
model category, see Proposition~\ref{prop-closure}.

\begin{definition}\label{def-stable-fib}
A map of symmetric spectra is a \emph{stable fibration} if it has the
right lifting property with respect to every map that is a stable
trivial cofibration. A map of symmetric spectra is a \emph{stable
trivial fibration} if it is a stable fibration and a stable equivalence.
A spectrum $X$ is \emph{stably fibrant} if the map $X\to *$ is a stable
fibration.
\end{definition}

\begin{theorem}\label{thm-stable-model-cat}
The category of symmetric spectra with the class of stable equivalences,
the class of stable cofibrations, and the class of stable fibrations is
a model category.
\end{theorem}

\begin{proof}
The category $\spec$ is bicomplete by
Proposition~\ref{prop-spec-bicomplete}. The two out of three axiom and
the retract axiom are immediate consequences of the definitions.  By
definition, $(i,p)$ has the lifting property when $i$ is a stable
trivial cofibration and $p$ is a stable fibration. The lifting axiom for
$i$ a stable cofibration and $p$ a stable trivial fibration is verified
in Corollary~\ref{cor-lift-cof-tfib}.  The two parts of the
factorization axiom are verified in Corollary~\ref{cor-factor-cof-tfib}
and Lemma~\ref{lem-factor-tcof-fib}.
\end{proof}

\begin{lemma}\label{lem-stfib=ltfib}
A map is a stable trivial fibration if and only if it is a level trivial
fibration.
\end{lemma}

\begin{proof}
Suppose $g$ is a level trivial fibration. By definition, every stable
cofibration has the \llp $g$ and in particular every stable trivial
cofibration has the \llp $g$. So $g$ is a stable fibration which is a
level equivalence and hence a stable equivalence. So $g$ is a stable
trivial fibration.

Conversely, suppose $g$ is a stable trivial fibration. 
Recall that at this point we do not know that $g$ has the \rlp stable
cofibrations.  By Lemma~\ref{lem-factorization}, $g$ can be factored as
$g=pi$ with $i$ an $FI_{\partial}$-cofibration and $p$ an
$FI_{\partial}$-injective map. Since $p$ is a level equivalence, it is a
stable equivalence.  By the two out of three property, $i$ is a stable
equivalence. Therefore, $i$ is a stable trivial cofibration and has the
\llp $g$.  By the Retract Argument~\ref{prop-retract}, $g$ is a retract
of $p$ and so $g$ is a level trivial fibration.
\end{proof}

\begin{corollary}\label{cor-factor-cof-tfib}
Every map $f$ of symmetric spectra has a factorization
$f=pi$ as a stable cofibration $i$ followed by a stable trivial
fibration $p$.
\end{corollary}
 
\begin{proof}
By the Factorization Lemma~\ref{lem-factorization}, every map $f$ in
$\spec$ can be factored as $f=pi$ with $i$ an
$FI_{\partial}$-cofibration and $p$ an $FI_{\partial}$-injective map. Then
$i$ is a stable cofibration and $p$ is level trivial fibration, which,
by Lemma~\ref{lem-stfib=ltfib}, means that $p$ is a stable
trivial fibration. 
\end{proof}

\begin{corollary}\label{cor-lift-cof-tfib}
A stable cofibration has the \llp every stable trivial fibration.
\end{corollary}

\begin{proof}
By Lemma~\ref{lem-stfib=ltfib} every stable trivial fibration is a
level trivial fibration. By definition,
stable cofibrations have the \llp every level trivial fibration.
\end{proof}

The following lemma will finish the proof of
Theorem~\ref{thm-stable-model-cat}.

\begin{lemma}\label{lem-factor-tcof-fib}
Every map $f$ of symmetric spectra has a factorization $f=pi$ as a
stable trivial cofibration $i$ followed by a stable fibration $p$.
\end{lemma}

To prove the lemma we need a set of maps $J$ such that a
$J$-cofibration is a stable trivial cofibration and a $J$-injective
map is a stable fibration. Using the Factorization Lemma with the set
$J$ will prove Lemma~\ref{lem-factor-tcof-fib}.  The set $J$ is
defined in~\ref{def-J} and Corollary~\ref{cor-J-works} verifies its
properties. This takes up the rest of the section.

The maps $\badmap\sm F_nS^0$ used in the definition below appeared in the
description of the function spectrum in
Remark~\ref{rem-describe-internal-hom}. They are stable equivalences
(see Example~\ref{eg-stable-equiv}) but are not stable cofibrations or
even level cofibrations.  We modify them to get the set $J$.

\begin{definition}\label{def-J}
Let $\badmap \mathcolon F_1S^1 \to F_0S^0$ be the adjoint of the
identity map $S^1\to\Ev_1F_0S^0=S^1$ and let $\badmap_n$ be the map
$\badmap\sm F_nS^0\mathcolon F_{n+1}S^1\to F_nS^0$, so that
$\badmap_0=\badmap$. The mapping cylinder
construction~\ref{mapping-cylinder} gives a factorization
$\badmap_n=r_nc_n$ where $r_n\mathcolon C\badmap_n\to F_nS^0$ is a
simplicial homotopy equivalence and $c_n\mathcolon F_{n+1}S^1\to
C\badmap_n$ is a level cofibration. For $n\ge0$, let $K_n=c_n\boxprod
I_{\partial}$, \ie $K_n$ is the set of maps $c_n\boxprod j$ for $j\in
I_{\partial}$.  Let $K=\cup_n K_n$ and finally, let $J=FI_\Lambda\cup
K$.
\end{definition}

\begin{lemma}\label{lem-stable-generators}
For each $n\ge0$, the map $c_n\mathcolon F_{n+1}S^{1}\to C\badmap_n$ is
a stable trivial cofibration.
\end{lemma}

\begin{proof}
The map $\badmap_n$ is a stable equivalence
(Example~\ref{eg-stable-equiv}) and the simplicial homotopy equivalence
$r_n$ is a stable equivalence. Using the factorization $\badmap_n=r_nc_n$
and the two out of three property of stable equivalences, $c_n$ is a
stable equivalence.

Next we show that $c_n$ is a stable cofibration.  The mapping cylinder
$C\badmap_n$ can also be defined as the corner in the pushout square
\[
\begin{CD}
 @. F_{n+1}S^1\vee F_{n+1}S^1  @>{i_0 \vee i_1}>> F_{n+1}S^1 \sm \Delta[1]_+\\
 @.    @V{\badmap_n \vee 1}VV         @VVV  \\
F_{n+1}S^1 @>k_n>> F_nS^0 \vee F_{n+1}S^1 @>g_n>> C\badmap_n,
\end{CD}
\]
where $k_n$ is the inclusion on the second factor, $i_0$ and $i_1$ come
from the two inclusions $\Delta[0]\to\Delta[1]$, and $g_n$ is the
natural map to the pushout.  Using the properties of stable cofibrations
in Proposition~\ref{prop-stable-cof}, we find that the map $*\to F_nS^0$
is a stable cofibration and, by cobase change, that $k_n$ is a stable
cofibration.  Let $j$ be the cofibration
$\partial\Delta[1]_+\to\Delta[1]_+$.  Then $(* \to F_{n+1}S^{1})\boxprod
j=i_0 \vee i_1$ is a stable cofibration and, by cobase change, $g_n$ is
a stable cofibration.  Thus the composition $c_n=g_nk_n$ is a stable
cofibration.
\end{proof}

Next we characterize the $J$-injective maps.

\begin{definition}\label{def-homotopy-pullback}
A commutative square of simplicial sets
\[
\begin{CD}
X @>>> Z \\
@VpVV @VqVV \\
Y @>>f> W
\end{CD}
\]
where $p$ and $q$ are fibrations is  a \emph{homotopy pullback
square} if the following equivalent conditions hold:
\begin{itemize}
\item The induced map $X\xrightarrow{}Y\times_{W}Z$ is a weak
equivalence.
\item
For every $0$-simplex $v\in Y_0$, the map of fibers $p^{-1}v\to q^{-1}fv$
is a weak equivalence.
\end{itemize}
\end{definition}

\begin{lemma}\label{lem-J-injective}
A map of symmetric spectra $p\mathcolon E\to B$ is $J$-injective if and
only if $p$ is a level fibration and the diagram
\begin{equation}
\begin{CD}
E_n @>{\sigma^a}>> E_{n+1}^{S^1} \\
@Vp_nVV         @VVp_{n+1}V\\
B_n @>{\sigma^a}>> B_{n+1}^{S^1},
\end{CD}\tag{$*$}
\end{equation}
is a homotopy pullback square for each $n\ge0$, where the horizontal maps
are the adjoints of the structure maps.
\end{lemma}

\begin{proof}
Since $J=FI_\Lambda\cup K$ and $K=\cup_n K_n$, a map is $J$-injective if
and only if it is $FI_\Lambda$-injective and $K_n$-injective for each
$n\ge0$. By Proposition~\ref{prop-level-fib}, the $FI_\Lambda$-injective
maps are the level fibrations. By definition, $p\in\spec$ is a
$K_n$-injective map if and only if $p$ has the \rlp the class
$c_n\boxprod I_\partial$.  Then, by Corollary~\ref{cor-boxhom-lifting},
$p$ is $K_n$-injective if and only if $\shom_\square(c_n,p)$ has the
\rlp the class $I_\partial$. Hence, $p$ is $K_n$-injective if and only
if $\shom_\square(c_n,p)$ is a trivial Kan fibration. If $p$ is a level
fibration, $\shom_\square(c_n,p)$ is a Kan fibration by
Proposition~\ref{prop-stable-cof}. So, a level fibration $p$ is
$K_n$-injective if and only if $\shom_\square(c_n,p)$ is a weak
equivalence.  Taken together, $p$ is $J$-injective if and only if $p$ is
a level fibration and $\shom_\square(c_n,p)$ is a weak equivalence for
each $n\ge0$.

For each $n\ge0$, the map $r_n\mathcolon C\badmap_n\to F_nS^0$ has a
simplicial homotopy inverse $s_n\mathcolon F_nS^0 \to C\badmap_n$ for
which $r_ns_n$ is the identity map on $F_nS^0$
(see~\ref{mapping-cylinder}). Then $\shom_\square(c_n,p)$ is
simplicially homotopic to $\shom_\square(s_n\badmap_n,p)$. Since $F_nS^0$
is a simplicial deformation retract of $C\badmap_n$, $\badmap_n$ is a
simplicial deformation retract of $s_n\badmap_n$ and
$\shom_\square(\badmap_n,p)$ is a simplicial deformation retract of
$\shom_\square(s_n\badmap_n,p)$. Therefore, $\shom_\square(c_n,p)$ is a
weak equivalence if and only if $\shom_\square(\badmap_n,p)$ is a weak
equivalence.

The map
\[
\shom_\square(\badmap_n,p)\mathcolon\shom(F_nS^0,E)\to\shom(F_nS^0,B)
\times_{\shom(F_{n+1}S^1, B)} \shom(F_{n+1}S^1, E)
\]
is naturally isomorphic to the map
\[
E_n\to B_n\times_{B_{n+1}^{S^1}} E_{n+1}^{S^1}.
\]
induced by the diagram $(*)$.  If $p$ is a level fibration then by
definition the diagram $(*)$ is a homotopy pullback square if and only
if the map $\shom_\square(\badmap_n,p)$ is a weak equivalence.

Combining the conclusions of the three paragraphs completes the proof.
\end{proof}

\begin{corollary}\label{cor-fibrant-is-omega-spectra}
The map $F \to *$ is $J$-injective if and only if $F$ is an
$\Omega$-spectrum.
\end{corollary}

We also get the following corollary, which is not needed in the sequel.
Its proof uses properness (see Section~\ref{subsec-proper}).  

\begin{corollary}\label{cor-fibration-between-omega-spectra}
A level fibration between two $\Omega$-spectra is $J$-injective.
\end{corollary}

\begin{lemma}\label{stable-fib-equiv-level}
Let $p\mathcolon X\to Y$ be a map of symmetric spectra.  If $p$ is
$J$-injective and $p$ is a stable equivalence then $p$ is a level
equivalence.
\end{lemma}

\begin{proof}
Suppose $p\mathcolon X\to Y$ is a $J$-injective stable equivalence. In
particular, $p$ is a level fibration. Let $F$ be the fiber over the
basepoint. Since the class $J\inj$ is closed under base change, the
map $F\to *$ is $J$-injective and $F$ is an $\Omega$-spectrum.  The
map $p$ factors as $X\to X/F \to Y$.  The map $X/F \to Y$ is a stable
equivalence by Theorem \ref{thm-traditional-results}. Since $p\mathcolon X\to
Y$ is a stable equivalence, $q\mathcolon X\to X/F$ is a stable
equivalence.

A Barratt-Puppe type sequence for symmetric
spectra is constructed by prolongation to give the diagram
\[
X\to X\amalg_F (F\sm \Delta[1]_+) \to F\sm S^1 \to X \sm S^1 \to
(X\amalg_F (F\sm \Delta[1]_+))\sm S^1.
\]
Let $E$ be an injective $\Omega$-spectrum. 
Since the  map $X\amalg_F (F\sm \Delta[1]_+)\to X/F$ is a level
equivalence, after applying $E^0(-)$ to this sequence we can rewrite the
terms involving the homotopy cofiber as $E^0(X/F)$.
This gives an exact sequence
\[
E^0X \xleftarrow{E^0q} E^0(X/F) \from E^0(F\sm S^1) \from E^0(X \sm
S^1) \xleftarrow{E^0(q\sm S^1)}
E^0(X/F\sm S^1).
\]
Since $q\mathcolon X\to X/F$ is a stable equivalence, $E^0q$ is an
isomorphism by definition, and $E^0(q\sm S^1)$ is an isomorphism by part
two of Theorem~\ref{thm-traditional-results}.  Hence, $E^0(F\sm S^1)=*$
for every injective $\Omega$-spectrum $E$, and so, by part two of
Theorem~\ref{thm-traditional-results}, $E^0F=*$ for every injective
$\Omega$-spectrum $E$.  By Corollary~\ref{cor-enough-injectives}, there
is a level equivalence $F\to E$ where $E$ is an injective spectrum;
since $F$ is an $\Omega$-spectrum, $E$ is an injective
$\Omega$-spectrum.  By Lemma~\ref{lem-injective-spectra},
$E^0E=E^0F=*$. So $E$ is simplicially homotopic to $*$ and $F$ is level
equivalent to $*$.

This does not finish the argument as the base of the fibration $X_n\to
Y_n$ need not be connected. Since $p$ is $J$-injective  
\[
\begin{CD}
X_n @>{\sigma^a}>> X_{n+1}^{S^1}\\
@Vp_nVV         @VVp_{n+1}^{S^1}V\\
Y_n @>{\sigma^a}>> Y_{n+1}^{S^1},
\end{CD}
\]
is a homotopy pullback square for each $n\ge0$. The proof is completed
by showing that $p_{n+1}^{S^1}$ is a trivial Kan fibration for every
$n\ge0$ which implies that $p_n$ is a trivial Kan fibration for every
$n\ge 0$.  For a pointed simplicial set $K$, let $cK$ denote the
connected component of the basepoint.  If $E\to B$ is a pointed Kan
fibration, then $cE\to cB$ is a Kan fibration; if the fiber over the
basepoint $*\in B$ is contractible then $cE\to cB$ is a trivial Kan
fibration. In particular, $cX_n\to cY_n$ is a trivial Kan fibration and
therefore, $(cX_n)^{S^1}\to (cY_n)^{S^1}$ is a trivial fibration.  Since
$K^{S^1}=(cK)^{S^1}$ for any pointed simplicial set $K$,
$p_n^{S^1}\mathcolon X_n^{S^1}\to Y_n^{S^1}$ is a trivial Kan fibration
for every $n\ge0$.
\end{proof}

The next corollary finishes the proof of
Lemma~\ref{lem-factor-tcof-fib}.

\begin{corollary}\label{cor-J-works}
The $J$-cofibrations are the stable trivial cofibrations and the 
$J$-injective maps are the stable fibrations.
\end{corollary}

\begin{proof}
Every level trivial fibration is $J$-injective since it satisfies the
condition in Lemma~\ref{lem-J-injective}.  Thus, a $J$-cofibration has
the \llp every level trivial fibration, and hence a $J$-cofibration is a
stable cofibration. Let $E$ be an $\Omega$-spectrum. The maps
$p\mathcolon E\to *$ and $q=\ihom(j,E)\mathcolon E^{\Delta[1]}\to
E\times E$, where $j\mathcolon\partial\Delta[1]_+\to\Delta[1]_+$ is the
inclusion, are $J$-injective by Lemma \ref{lem-J-injective}. Let $E$ be
an injective $\Omega$-spectrum and $f$ be a $J$-cofibration. Since $f$
has the \llp $p\mathcolon E\to *$, $E^0f$ is surjective. Since $f$ has
the \llp $q\mathcolon E^{\Delta[1]}\to E\times E$, $E^0f$ is injective.
So $E^0f$ is an isomorphism and every $J$-cofibration is a stable
trivial cofibration.

Conversely, let $f$ be a stable trivial cofibration.  By the
Factorization Lemma \ref{lem-factorization}, $f$ factors as $f=pi$ where
$i$ is a $J$-cofibration and $p$ is a $J$-injective map.  We have just
seen that $i$ is a stable equivalence.  So, the $J$-injective map $p$ is
a stable equivalence and by Lemma~\ref{stable-fib-equiv-level}, $p$ is a
level equivalence.  Therefore the stable cofibration $f$ has the \llp
the map $p$.  By the Retract Argument~\ref{prop-retract}, $f$ is a
$J$-cofibration.

Let $F$ be the class of stable fibrations.  Since $J\cof$ is the class
of stable trivial cofibrations, one has by the definition of stable
fibrations that $F=(J\cof)\inj$.  But $(J\cof)\inj=J\inj$ by Proposition
\ref{prop-lifting-properties} (2).  In other words, the stable
fibrations are the $J$-injective maps.
\end{proof}

In particular, Lemma~\ref{lem-J-injective} characterizes the stable
fibrations. The stably fibrant objects are the $\Omega$-spectra by
Corollary \ref{cor-fibrant-is-omega-spectra}.
Corollary~\ref{cor-J-works} finishes the proof of
Lemma~\ref{lem-factor-tcof-fib} and the verification of the axioms for
the stable model category of symmetric spectra.

\section{Comparison with the Bousfield-Friedlander
category}\label{sec-comparison}

The goal of this section is to show that the stable homotopy theory of
symmetric spectra and the stable homotopy theory of spectra are
equivalent.  We begin in Section~\ref{subsec-Quillen-equiv} by recalling
the general theory of Quillen equivalences of model categories.  In
Section~\ref{subsec-stable-BF} we provide a brief recap of the stable
homotopy theory of (non-symmetric) spectra.  In
Section~\ref{subsec-equivalence} we show that the forgetful functor $U$
from symmetric spectra to spectra is part of a Quillen equivalence.  The
left adjoint $V$ of $U$ plays very little role in this proof, beyond its
existence, so we postpone its construction to Section~\ref{subsec-V}.

\subsection{Quillen equivalences}\label{subsec-Quillen-equiv}

In this section, we briefly recall Quillen functors and Quillen
equivalences between model categories.  

\begin{definition}\label{def-Quillen-functor}
Let $\cc$ and $\cd$ be model categories. Let $L\mathcolon\cc\to\cd$ and
$R\mathcolon\cd\to\cc$ be functors such that $L$ is left adjoint to $R$.
The adjoint pair of functors $L$ and $R$ is a \emph{Quillen adjoint
pair} if $L$ preserves cofibrations and $R$ preserves fibrations.  We
refer to the functors in such a pair as \emph{left and right Quillen
functors}.  A Quillen adjoint pair is a \emph{Quillen equivalence} if
for every cofibrant object $X\in \cc$ and every fibrant object $Y\in
\cd$, a map $LX\xrightarrow{}Y$ is a weak equivalence if and only if its
adjoint $X\xrightarrow{}RY$ is a weak equivalence.
\end{definition}

The definition of a Quillen adjoint pair can be reformulated. 

\begin{lemma}\label{lem-quillen-adjoints}
Let $L$ and $R$ be a pair of functors between model categories such that
$L$ is left adjoint to $R$.
\begin{enumerate}
\item $L$ preserves cofibrations if and only if $R$ preserves trivial
fibrations.
\item $L$ preserves trivial cofibrations if and only if $R$ preserves
fibrations.
\end{enumerate}
\end{lemma}

This lemma is an immediate corollary of
Lemma~\ref{lem-lifting-adjoints}; see also~\cite[9.8]{dwyer-spalinski}.
A useful lemma associated to these questions is Ken Brown's lemma.

\begin{lemma}[Ken Brown's Lemma]\label{lem-ken-brown}
Let $F$ be a functor between model categories.
\begin{enumerate}
\item If $F$ takes trivial cofibrations between cofibrant objects to
weak equivalences, then $F$ preserves all weak equivalences between
cofibrant objects.
\item If $F$ takes trivial fibrations between fibrant objects to weak
equivalences, then $F$ preserves all weak equivalences between fibrant
objects.
\end{enumerate}
\end{lemma}

For the proof of this lemma see~\cite[9.9]{dwyer-spalinski}.

In particular, a left Quillen functor $L$ preserves weak equivalences
between cofibrant objects, and a right Quillen functor $R$ preserves
weak equivalences between fibrant objects.

The following proposition is the reason Quillen equivalences are
important.  

\begin{proposition}\label{prop-induced-quillen}
A Quillen adjoint pair of functors between model categories induces an
adjoint pair of functors on the homotopy categories which is an adjoint
equivalence if and only if the adjoint pair of functors is a Quillen
equivalence.
\end{proposition}

For the proof of this proposition, see~\cite[Theorem 9.7]{dwyer-spalinski}.

We now describe a useful sufficient condition for a Quillen adjoint pair
to be a Quillen equivalence.  

\begin{definition}\label{def-reflect}
Suppose $F\mathcolon\cc\to\cd$ is a functor between model categories.  For
any full subcategory $\mathcal{C}'$ of $\cc$, we say that $F$
\emph{detects and preserves weak equivalences of $\mathcal{C}'$} if a
map $f$ in $\mathcal{C}'$ is a weak equivalence if and only if $Ff$ is.
\end{definition}

In practice, very few functors detect and preserve weak equivalences on
the whole category. However, many functors detect and preserve weak 
equivalences between cofibrant objects or fibrant objects, so the next
lemma is often useful.  Before stating it, we need a definition.  

\begin{definition}\label{def-fibrant-replacement}
Suppose $\cc $ is a model category.  A \emph{fibrant replacement
functor} on $\cc $ is a functor $K\mathcolon \cc \to \cc $ whose image
lies in the full subcategory of fibrant objects, together with a natural
weak equivalence $i \mathcolon X\to KX$.
\end{definition}

There is a dual notion of a cofibrant replacement functor, but we do not
use it.  Fibrant replacement functors are usually obtained by using a
version of the Factorization Lemma~\ref{lem-factorization} appropriate
for $\cc $ to functorially factor the map $X\xrightarrow{}1$ into a
trivial cofibration followed by a fibration.  We have already used
fibrant replacement functors in $\sset $ in the proof of
Theorem~\ref{thm-U-detects-stable}.  

\begin{lemma}\label{lem-Quillen-equivalences}
Suppose $L\mathcolon\cc\to\cd$ is a left Quillen functor with right
adjoint $R$, and suppose $K$ is a fibrant replacement functor on $\cd$.
Suppose $R$ detects and preserves weak equivalences between fibrant
objects and the composition $X\to RLX\xrightarrow{Ri}RKLX$ is a weak
equivalence for all cofibrant objects $X$ of $\cc$.  Then the pair $(L,
R)$ is a Quillen equivalence.
\end{lemma}

There is also a dual statement, but this is the criterion we use.

\begin{proof}
Suppose $f\mathcolon LX\xrightarrow{}Y$ is a map, where $X$ is cofibrant
and $Y$ is fibrant.  Consider the commutative diagram below.
\[
\begin{CD}
X @>>> RLX @>Rf>> RY \\
@. @VRi_{LX}VV @VVRi_{Y}V \\
\ @. RKLX @>>RKf>> RKY
\end{CD}
\]
The top composite is the adjoint $g\mathcolon X\xrightarrow{}RY$ of $f$.
The map $i_{Y}$ is a weak equivalence between fibrant objects, so
$Ri_{Y}$ is a weak equivalence.  The composite $X\xrightarrow{} RLX
\xrightarrow{Ri_{LX}} RKLX$ is a weak equivalence by hypothesis.  Thus
$g$ is a weak equivalence if and only if $RKf$ is a weak equivalence.
But $R$ detects and preserves weak equivalences between fibrant objects,
so $RKf$ is a weak equivalence if and only if $Kf$ is a weak
equivalence.  Since $i$ is a natural weak equivalence, $Kf$ is a weak
equivalence if and only if $f$ is a weak equivalence.
\end{proof}

\subsection{The stable Bousfield-Friedlander
category}\label{subsec-stable-BF} 

In this section we describe the stable homotopy theory of (non-symmetric)
spectra.  The basic results are proved in~\cite{bousfield-friedlander},
but we use the approach of Section~\ref{sec-homotopy}.  In particular,
we can define injective spectra and $\Omega $-spectra just as we did for
symmetric spectra.  

\begin{definition}\label{defn-BF}
Suppose $f\mathcolon X\xrightarrow{}Y$ is a map of spectra.  
\begin{enumerate}
\item The map $f$ is a \emph{stable equivalence} if $E^{0}f$ is an
isomorphism for every injective $\Omega $-spectrum $E$. 
\item The map $f$ is a \emph{stable cofibration} if $f$ has the \llp
every level trivial fibration.  
\item The map $f$ is a \emph{stable fibration} if $f$ has the \rlp every
map which is both a stable cofibration and a stable equivalence.  
\end{enumerate}
\end{definition}

We can then attempt to carry out the program of
Section~\ref{sec-homotopy} with these definitions.  The category $\BF $
of non-symmetric spectra is neither symmetric monoidal nor closed, so we
must be careful in some places.  However, if $X$ is a sequence and $Y$
is a spectrum, we can form the sequences $X\otimes Y$ and $\ihom
_{\cN}(X,Y)$, and these are $S$-modules in a canonical way.
With this minor change, all of the results of Section~\ref{sec-homotopy}
go through without difficulty.  There is one important change: in the
proof of Theorem~\ref{thm-U-detects-stable}, we considered the spectrum
$R^{\infty }X$, which looks like an $\Omega $-spectrum but is not in
general. The analogous spectrum in $\BF $ is an $\Omega $-spectrum, so
we obtain the following result.

\begin{theorem}\label{thm-BF-U-detects-stable}
A map in $\BF $ is a stable equivalence if and only if it is a stable
homotopy equivalence.
\end{theorem}

We also obtain the usual characterizations of the stable trivial
fibrations and the stable fibrations.  

\begin{proposition}\label{prop-BF-stable-triv-fib}
A map in $\BF$ is a stable trivial fibration if and only if it is a level 
trivial fibration.
\end{proposition}

\begin{proposition}\label{prop-stable-BF-fibration}
A map $f\mathcolon E\xrightarrow{}B$ in $\BF $ is a stable fibration if
and only if it is a level fibration and, for all $n\geq 0$, the diagram 
\[
\begin{CD}
E_n @>>> E_{n+1}^{S^1}\\
@VVV    @VVV\\
B_n @>>> B_{n+1}^{S^1}
\end{CD}
\]
is a homotopy pullback square.
\end{proposition}

\begin{corollary}\label{cor-BF-fibrant-is-omega}
A spectrum $X$ in $\BF $ is fibrant in the
stable model category if and only if $X$ is an $\Omega $-spectrum.
\end{corollary}

Finally, we recover the theorem of Bousfield and Friedlander.

\begin{theorem}[~\cite{bousfield-friedlander}]\label{thm-BF}
The stable equivalences, stable cofibrations, and stable fibrations
define a model structure on $\BF $.  
\end{theorem}

Note that the model structure we have just described is the same as the
model structure in ~\cite{bousfield-friedlander}, since the stable
equivalences and stable cofibrations are the same.  

In particular, we have the following characterization of stable
cofibrations.

\begin{proposition}\label{prop-BF-cofib}
A map $f\mathcolon X\xrightarrow{}Y$ in $\BF $ is a stable cofibration
if and only if $f_{0}\mathcolon X_{0}\xrightarrow{}Y_{0}$ is a
monomorphism and the induced map $X_{n}\amalg _{X_{n-1}\sm
S^{1}}(Y_{n-1}\sm S^{1})\xrightarrow{}Y_{n}$ is a monomorphism for all
$n>0$.  In particular, $X$ is cofibrant if and only if the structure
maps are monomorphisms.  
\end{proposition}

This proposition is proved in~\cite{bousfield-friedlander}; the proof is
very similar to the proof of the analogous fact for symmetric
spectra, proved in Proposition~\ref{prop-stable-cofib}.  

\subsection{The Quillen equivalence}\label{subsec-equivalence}   

The goal of this section is to show that the forgetful functor $U$ from
symmetric spectra to spectra is part of a Quillen equivalence.
Obviously this requires that $U$ have a left adjoint $V\mathcolon \BF
\xrightarrow{}\spec $.  We will assume the existence of $V$ in this
section, and construct $V$ in Section~\ref{subsec-V}.  

\begin{proposition}\label{prop-comparison}
The functors $U\mathcolon \spec \xrightarrow{}\BF$ and its left adjoint
$V$ are a Quillen adjoint pair.
\end{proposition}

\begin{proof}
Proposition~\ref{prop-stable-BF-fibration} implies that $U$ preserves
stable fibrations.  The stable trivial fibrations in $\spec $ and in
$\BF $ are the level trivial fibrations, so $U$ preserves stable trivial
fibrations as well.
\end{proof} 

\begin{theorem}\label{thm-comparison}
The functor $U\mathcolon \spec \to \BF $ and its left adjoint
$V$ form a Quillen equivalence of the stable model categories.
\end{theorem}

We prove this theorem by using
Lemma~\ref{lem-Quillen-equivalences}.  In particular, we need to
understand stable equivalences between stably fibrant objects. 

\begin{lemma}\label{lem-stably-fibrant-equiv}
Suppose $f\mathcolon X\xrightarrow{}Y$ is a stable equivalence between
stably fibrant objects in either $\spec $ or $\BF $.  Then $f$ is a
level equivalence.  
\end{lemma}

\begin{proof}
Factor $f$ as a stably trivial cofibration, $i$, followed by a stable
fibration, $p$.  Since $f$ is a stable equivalence, $p$ is also.  Hence,
$p$ is a level trivial fibration by Lemma~\ref{lem-stfib=ltfib}.  Also,
$i$ is a stably trivial cofibration between stably fibrant objects,
hence it is a strong deformation retract, see~\cite[II
p. 2.5]{quillen-htpy}.  To see this, note that $i$ has the left lifting
property with respect to $X \to *$, so the lift constructs a homotopy
inverse to $i$.  Because the simplicial structure is given on levels, a
strong deformation retract here is a level equivalence.  So both $i$ and
$p$ are level equivalences, hence so is $f$.
\end{proof}

\begin{corollary}\label{lem-U-reflects-equiv}
$U\mathcolon \spec \xrightarrow{}\BF $ detects and preserves 
stable equivalences between stably fibrant objects.
\end{corollary}

Let $L$ denote a fibrant replacement functor in $\spec $, obtained by
factoring $X\xrightarrow{}*$ into a stable trivial cofibration followed
by a stable fibration.  By Lemma~\ref{lem-Quillen-equivalences} and
Corollary~\ref{lem-U-reflects-equiv}, to prove
Theorem~\ref{thm-comparison} it suffices to show that $X \to ULVX$ is a
stable equivalence for all cofibrant (non-symmetric) spectra $X$.  We
prove this in several steps.

\begin{definition}
Given a simplicial set $X$, define $\wt{F}_{n}(X)$ to be the
(non-symmetric) spectrum whose $m$th level is $S^{m-n}\sm X$ for $m\geq
n$ and the basepoint otherwise, with the obvious structure maps.  This
defines a functor $\wt{F}_{n}\mathcolon \sset \xrightarrow{}\BF $ left
adjoint to the evaluation functor $\Ev _{n}$.  
\end{definition}

Note that $\wt{F}_{0}X=\Sigma ^{\infty }X$.  Also, since $U\circ \Ev
_{n}=\Ev _{n}$, the left adjoints satisfy $V\circ \wt{F}_{n}=F_{n}$.  

\begin{lemma}\label{lem-comparison-free}
The map $X\xrightarrow{}ULVX$ is a stable equivalence when
$X=\Sigma^{\infty}Y=\wt{F}_{0}Y$ for any $Y\in \sset $.
\end{lemma}

\begin{proof}
Consider the functor on simplicial sets $QZ=\colim\Omega^n K \Sigma^n
Z$, where $K$ is a simplicial fibrant replacement functor.  Because $Q$
is simplicial we can prolong it to a functor on $\spec$.  The map
$F_{0}Y\xrightarrow{}QF_{0}Y$ induces an isomorphism on stable homotopy.
Also $QF_0 Y$ is an $\Omega $-spectrum since $QZ \to \Omega Q\Sigma Z$ is a
weak equivalence for any $Z \in \sset$.  Hence $QF_0Y$ is level
equivalent to $LF_0Y$, so $F_0Y \to LF_0Y$ induces an isomorphism in
stable homotopy.  Since $\widetilde{F}_0 Y \to UF_0Y$ is a level
equivalence and $UF_0 Y \to ULF_0 Y$ is a stable homotopy equivalence,
the lemma follows.
\end{proof}

Because both $\spec $ and $\BF $ are stable model categories, the
following lemma is expected.  

\begin{lemma}\label{lem-comparison-suspension-free}
Suppose $X$ is a cofibrant spectrum in $\BF $.  Then the map
$X\xrightarrow{}ULVX$ is a stable equivalence if and only if $X\sm
S^{1}\xrightarrow{}ULV(X\sm S^{1})$ is a stable equivalence.
\end{lemma}

\begin{proof}
For notational convenience, we write $\Sigma X$ for $X\sm S^{1}$ and
$\Omega X$ for $X^{S^{1}}$ in this proof, for $X$ a (possibly symmetric)
spectrum.  Consider the stable trivial cofibration $\Sigma
VX\xrightarrow{}L\Sigma VX$ in $\spec $.  By
Theorem~\ref{thm-traditional-results} part three, $VX \to \Omega L\Sigma
VX$ is also a stable equivalence.  By the lifting property of the stable
trivial cofibration $VX \to LVX$ and the 2-out-of-3 property, there is a
stable equivalence $LVX \to \Omega L\Sigma VX$.  This map is a stable
equivalence between stably fibrant objects, so by
Corollary~\ref{lem-U-reflects-equiv}, $f\mathcolon ULVX \to U\Omega L
\Sigma VX$ is a stable equivalence.

So $g \mathcolon X \to ULVX$ is a stable equivalence if and only if
$fg\mathcolon X \to U\Omega L \Sigma VX$ is a stable equivalence.  Since
$\Omega$ and $U$ commute, $gf$ is a stable equivalence if and only if
$\Sigma X \to U L \Sigma VX$ is a stable equivalence by part three of
Theorem~\ref{thm-traditional-results} for (non-symmetric) spectra.  But,
since $U$ commutes with $\Omega$, the left adjoint $V$ commutes with
$\Sigma$, so we have a natural isomorphism $U L \Sigma VX \to ULV \Sigma
X$.  This completes the proof.
\end{proof}

\begin{lemma}\label{lem-good-properties}
Let $f\mathcolon X \to Y$ be a stable equivalence between cofibrant
spectra in $\BF$.  Then $X \to ULVX$ is a stable equivalence if and only
if $Y \to ULVY$ is a stable equivalence.
\end{lemma}

\begin{proof}
Consider the following commuting square.
\[
\begin{CD}
X @>>> Y \\
@VVV   @VVV\\
ULVX @>>> ULVY
\end{CD}
\]
Since $V$ is a left Quillen functor by
Proposition~\ref{prop-comparison}, it preserves trivial cofibrations.
Hence, by Ken Brown's Lemma \ref{lem-ken-brown}, $V$ preserves stable
equivalences between cofibrant objects.  Hence $VX \to VY$ is a stable
equivalence.  $L$ takes stable equivalences to level equivalences, by
Lemma~\ref{lem-stably-fibrant-equiv}.  So $ULVX \to ULVY$ is a level
equivalence since $U$ preserves level equivalences.  Hence the top and
bottom maps are stable equivalences, so the right map is a stable
equivalence if and only if the left map is.
\end{proof}

Using the preceding three lemmas we can extend
Lemma~\ref{lem-comparison-free} to any cofibrant strictly bounded below
spectrum.
  
\begin{definition}\label{def-bounded-above}
Define a spectrum $X\in \BF $ to be \emph{strictly bounded below} if
there is an $n$ such that for all $m\geq n$ the structure map $S^{1}\sm
X_{m}\xrightarrow{}X_{m+1}$ is an isomorphism.
\end{definition}

\begin{lemma}\label{lem-comparison-bounded-above}
Suppose $X\in \BF $ is cofibrant and
strictly bounded below.  Then the map $X\xrightarrow{}ULVX$ is a stable
equivalence.
\end{lemma}

\begin{proof}
Suppose $X$ is strictly bounded below at $n$.  Then we have a map
$\wt{F}_nX_{n}\xrightarrow{g} X$ which is the identity on all levels
$\geq n$.  In particular, $g$ is a stable homotopy equivalence.
Applying Lemma~\ref{lem-good-properties}, this shows that to prove the
lemma it is enough to show that $\wt{F}_nX_n \to ULF_nX_n$ is a stable
equivalence.  But there is an evident map $\wt{F}_{n}X_{n}\sm
S^{n}\natiso \wt{F}_{n}(X_{n}\sm S^{n})\to \Sigma^{\infty} X_n$ which is
the identity map above level $n-1$, and so is a stable equivalence.
Lemmas~\ref{lem-comparison-free}, \ref{lem-comparison-suspension-free},
and~\ref{lem-good-properties} complete the proof.  
\end{proof}

We now extend this lemma to all cofibrant objects, completing the proof
of Theorem~\ref{thm-comparison}.  First, we need to recall a basic fact
about simplicial sets.  Recall that the homotopy group $\pi _{n}X$ of a
pointed simplicial set $X$ is defined to be $\pi _{0}\shom _{\sset
}(S^{n},KX)$, where $K$ is a fibrant replacement functor.  This ensures
that weak equivalences are homotopy isomorphisms.  If $X$ is already a
Kan complex, $X$ is simplicially homotopy equivalent to $KX$, and so
$\pi _{n}X\natiso \pi _{0}\shom _{\sset }(S^{n},X)$.  Since the
simplicial sets $\partial \Delta [n]_{+}$ and $\Delta [n]_{+}$ are
finite, the colimit of a sequence of Kan complexes is again a Kan
complex.  Since the simplicial sets $S^{n}$ and $S^{n}\sm \Delta
[1]_{+}$ are finite, homotopy commutes with
filtered colimits of Kan complexes, and in particular with transfinite
compositions of maps of Kan complexes.

In fact, homotopy commutes with transfinite compositions of arbitrary
monomorphisms of simplicial sets.  To see this, apply the geometric
realization to get a sequence of cofibrations of CW complexes.
Since homotopy commutes with such transfinite compositions, the result
follows.  

\begin{lemma}\label{lem-comparison-arbitrary}
Suppose $X$ is a cofibrant object of $\BF $.  Then the map
$X\xrightarrow{}ULVX$ is a stable equivalence.
\end{lemma}

\begin{proof}
Let $X^{i}$ denote the truncation of $X$ at $i$.  That is, we have
$X^{i}_{n}=X_{n}$ for $n\leq i$ and $X^{i}_{n}=X_{i}\wedge S^{n-i}$
for $n\geq i$.  Then the $X^{i}$ are strictly bounded below and
cofibrant, and there are monomorphisms $X^{i}\xrightarrow{}X^{i+1}$ with
$\colim_i X^{i}=X$.  Thus each map $X^{i}\xrightarrow{}ULVX^{i}$ is a
stable equivalence.  

We claim that the induced map $X\xrightarrow{}\colim_i ULVX^{i}$ is a
stable equivalence.  To see this, note that 
\[
\pi _{n}X=\colim_{i} \pi _{n+i}X_{i} = \colim _{i} \pi _{n+i}\colim
_{j}X_{i}^{j} 
\]
Since homotopy commutes with transfinite compositions of monomorphisms, 
we find that $\pi _{n}X\natiso \colim _{j}\pi _{n}X^{j}$.  Similarly,
since homotopy of Kan complexes commutes with arbitrary filtered
colimits, we find $\pi _{n}\colim_i ULVX^{i}\natiso \colim_i \pi
_{n}ULVX^{i}$.  It follows that $X\xrightarrow{}\colim_i ULVX^{i}$ is a
stable homotopy equivalence, as required. 

We now examine the relationship between $ULVX$ and $\colim_i ULVX^{i}$.
Since $V$ is a left adjoint, $VX\natiso \colim_i VX^{i}$.  Each map
$VX^{i}\xrightarrow{}LVX^{i}$ is a stable trivial cofibration; we claim
that the induced map $\colim_i VX^{i}\xrightarrow{}\colim_i LVX^{i}$ is a
stable equivalence.  To see this, we define a new sequence $Y^{i}$ and
maps of sequences $VX^{i}\xrightarrow{}Y^{i}$ and
$Y^{i}\xrightarrow{}LVX^{i}$ inductively.  Define $Y^{0}=LVX^{0}$.
Having defined $Y^{i}$, define $Y^{i+1}$ by factoring the map
$Y^{i}\amalg _{VX^{i}}VX^{i+1}\xrightarrow{}LVX^{i+1}$ into a stable
trivial cofibration $Y^{i}\amalg _{VX^{i}}VX^{i+1}\xrightarrow{}Y^{i+1}$
followed by a stable fibration $Y^{i+1}\xrightarrow{}LVX^{i+1}$.  Then
the induced map $\colim_i VX^{i}\xrightarrow{}\colim_i Y^{i}$ is a stable
trivial cofibration, by a lifting argument.  On the other hand, each map
$Y^{i}\xrightarrow{}LVX^{i}$ is a stable equivalence, by the two out of
three axiom.  Since $LVX^{i}$ and hence $Y^{i}$ are stably fibrant, the
maps $Y^{i}\xrightarrow{}LVX^{i}$ are level equivalences, by
Lemma~\ref{lem-stably-fibrant-equiv}.  Since homotopy of level Kan
complexes commutes with filtered colimits, we find that $\colim
Y^{i}\xrightarrow{}\colim_i LVX^{i}$ is a level equivalence, and therefore
that $VX \natiso \colim_i VX^{i}\xrightarrow{}\colim_i LVX^{i}$ is a stable
equivalence.  

We now claim that $\colim_i LVX^{i}$ is an $\Omega $-spectrum, and thus is
stably fibrant.  Indeed, $\colim_i LVX^{i}$ is a level Kan complex by the
comments preceding this lemma.  Similarly, $(\colim
LVX^{i})_{n+1}^{S^{1}}=\colim_i ((LVX^{i})_{n+1}^{S^{1}})$.  Since homotopy
of Kan complexes commutes with filtered colimits, it follows that
$\colim_i LVX^{i}$ is an $\Omega $-spectrum.  

Hence the stable equivalence $VX \xrightarrow{}\colim_i LVX^{i}$ extends
to a stable equivalence $LVX \xrightarrow{}\colim_i LVX^{i}$.  By
Lemma~\ref{lem-stably-fibrant-equiv}, this map is actually a level
equivalence.  Since $U$ preserves level equivalences and colimits, the
map $ULVX\xrightarrow{}\colim_i ULVX^{i}$ is also a level equivalence.
We have seen above that the map $X\xrightarrow{}\colim_i ULVX^{i}$ is a
stable equivalence, so $X\xrightarrow{}ULVX$ must also be a stable
equivalence.
\end{proof}

\begin{remark}\label{rem-smash-comparison}
It follows from the results of Section~\ref{subsec-pushout-smash} that
the smash product on $\spec $ induces a smash product on $\ho \spec $.
The handicrafted smash products of~\cite{adams-blue} induce a smash
product on $\ho \BF $.  We now consider to what extent the equivalence
$RU\mathcolon \ho \spec \xrightarrow{}\ho \BF $ induced by $U$ preserves
these smash products.  Since $U$ is a simplicial functor, there is a
natural isomorphism $RU(X\sm Y)\natiso (RU)(X)\sm (RU)(Y)$ for all
(arbitrary desuspensions of) suspension spectra $X$.  On the other hand,
in either $\ho \spec $ or $\ho \BF $, $X\sm Y$ is determined by the
collection of $F\sm Y$ for all finite spectra $F$ mapping to $X$.  To be
precise, $X\sm Y$ is the minimal weak colimit~\cite{hovey-axiomatic} of
the $F\sm Y$.  As an equivalence of categories, $RU$ preserves minimal
weak colimits, so there is an isomorphism $RU(X\sm Y)\natiso (RU)(X)\sm
(RU)(Y)$, However, we do not know if this is natural, as the minimal
weak colimit is only a weak colimit.  This isomorphism is natural up to
phantom maps, however.
\end{remark}

\subsection{Description of $V$}\label{subsec-V}

This short section is devoted to the construction of the left adjoint
$V\mathcolon \BF \xrightarrow{}\spec $ to the forgetful functor
$U\mathcolon \spec \xrightarrow{}\BF $.  

Recall that, in any cocomplete symmetric monoidal category $\cc$, the
\emph{free monoid} or \emph{tensor algebra} generated by an object $X$ is
$T(X)= e \vee X \vee X^{\otimes 2} \vee \cdots \vee X^{\otimes n} \vee
\cdots$, where $e$ is the unit and $\vee$ is the coproduct. The
multiplication on $T(X)$ is the concatenation $X^{\otimes n} \otimes
X^{\otimes m} \to X^{\otimes (n+m)}$.  Similarly, the \emph{free
commutative monoid} on an object $X$ is $\sym(X)= e \vee X \vee
(X^{\otimes 2}/ \Sigma_2) \vee \cdots \vee (X^{\otimes n}/\Sigma_n) \vee
\cdots$.

Recall that the evaluation functor $\Ev _{n}\mathcolon \symseq
\xrightarrow{}\sset $ has a left adjoint $G_{n}$, where $G_{n}X$ is
$(\Sigma _{n})_{+}\sm X$ at level $n$ and the basepoint everywhere
else.  Similarly, the evaluation functor $\Ev _{n}\mathcolon \seq
\xrightarrow{}\sset $ has a left adjoint $\wt{G}_{n}$, where
$\wt{G}_{n}$ is $X$ at level $n$ and the basepoint everywhere else. 

\begin{lemma}\label{lem-comm-monoid}
In the category $\seq $ of sequences, the sphere spectrum $S$ is the
tensor algebra on the sequence $\wt{G}_{1}S^{1}=(*,S^{1},*,\dots
,*,\dots )$.  In the category $\symseq $ of symmetric sequences, the
sphere symmetric spectrum $S$ is the free commutative monoid on the
symmetric sequence $G_{1}S^{1}=(*,S^{1},*,\dots ,*,\dots )$.
\end{lemma}

\begin{proof}
The first statement follows directly from the definitions.  In the
category of symmetric sequences, $(G_{1}S^{1})^{\otimes n}=G_{n}S^{n}$,
so $T(G_{1}S^{1})$ is $(\Sigma _{n})_{+}\sm S^{n}$ in degree $n$.
Therefore $\sym(G_{1}S^{1})$ is $S^{n}$ in degree $n$.  Since we already
know $S$ is a commutative monoid, the map $G_{1}S^{1}\xrightarrow{}S$
induces a map $\sym(G_{1}S^{1})\xrightarrow{}S$ which is an isomorphism.
\end{proof}

This lemma explains why left $S$-modules and right $S$-modules are
equivalent in the category of sequences, since this is true for any
tensor algebra.  This lemma also explains why
Remark~\ref{rem-transpositions} holds, since an analogous statement
holds for any free commutative monoid.  

Now, the forgetful functor $U\mathcolon \symseq \xrightarrow{}\seq $ has
a left adjoint $G$, defined by $GX=\bigvee G_{n}X_{n}$, so that the
$n$th level of $GX$ is just $(\Sigma _{n})_{+}\sm X_{n}$.  The functor
$G$ is monoidal; that is, there is a natural isomorphism $G(X)\otimes
G(Y)\xrightarrow{}G(X\otimes Y)$ compatible with the associativity and
unit isomorphisms.  However, $G$ is definitely not a \emph{symmetric}
monoidal functor; this natural isomorphism is not compatible with the
commutativity isomorphisms.  This explains how $S$ can be commutative in
$\symseq $ yet $US=S$ is not commutative in $\seq $.  

Since $G$ is a monoidal functor, $G$ preserves monoids and modules, and
so defines a functor $G\mathcolon \BF \xrightarrow{}T(G_{1}S^{1})$-mod,
left adjoint to the forgetful functor
$T(G_{1}S^{1})\text{-mod}\xrightarrow{}\BF $.  On the other hand, the
map of monoids $T(G_{1}S^{1})\xrightarrow{p}\sym (G_{1}S^{1})=S$ defines
the usual adjoint pair of induction and restriction.  Induction takes a
(left) $T(G_{1}S^{1})$-module $X$ to $S\otimes _{T(G_{1}S^{1})}X$, where
the tensor product uses the right action of $T(G_{1}S^{1})$ on $S$
determined by $p$.  It follows that the left adjoint $V\mathcolon \BF
\xrightarrow{}\spec $ of the forgetful functor $U\mathcolon \spec
\xrightarrow{}\BF $ is $V(X)=S\otimes _{T(G_{1}S^{1})} GX$.

\section{Additional properties of symmetric spectra}\label{sec-misc}

In this section we discuss some properties of the category of symmetric
spectra.  In Section~\ref{subsec-level-model}, we consider the level
model structures on $\spec $.  In particular, we show that every
symmetric spectrum embeds in an injective spectrum by a level
equivalence, completing the proof that the stable structures define a
model structure on $\spec $.  In Section~\ref{subsec-stable-cofib} we
characterize the stable cofibrations.  In
Sections~\ref{subsec-pushout-smash} and~\ref{subsec-monoids}, we study
the relationship between the stable model structure on $\spec $ and the
smash product.  This is necessary for constructing model categories of
monoids, algebras, and modules, as is done
in~\cite{schwede-shipley-monoids}.  In Section~\ref{subsec-proper}, we
show that the stable model structure on $\spec $ is proper.  Finally,
in Section~\ref{sec-semistable} we define semistable spectra and investigate
their relationship to stable homotopy equivalences.

\subsection{Level model structure}
\label{subsec-level-model}

In this section we construct the two level model structures on the
category of symmetric spectra. 

\begin{definition}
A \emph{projective cofibration} of symmetric spectra is a map that has
the \llp every level trivial fibration. The projective cofibrations are
the stable cofibrations from Section~\ref{subsec-stable-model-cat}. The
\emph{projective level structure} on $\spec$ is the class of level
equivalences, the class of projective cofibrations, and the class of
level fibrations.  An \emph{injective fibration} of symmetric spectra is
a map that has the \rlp every level trivial cofibration (the adjective
``injective'' refers to the lifting properties of the map and not to its
being a monomorphism). The \emph{injective level structure} is the class
of level equivalences, the class of level cofibrations, and the class of
injective fibrations.
\end{definition}

\begin{theorem}\label{thm-level-model}
The projective level structure and the injective level
structure are model structures on the category of symmetric spectra. 
\end{theorem}

\begin{proof}
The category of symmetric spectra is bicomplete. The class of level
equivalences has the two-out-of-three property. The retract axiom holds
by construction in both the projective and injective level structures.

We now prove the lifting and factorization axioms, beginning with the
projective level structure.  We use the sets of maps $FI_{\partial }$
and $FI_{\Lambda }$ defined in Definition~\ref{def-level-generators}.
The lifting axiom for a projective cofibration and a level trivial
fibration holds by definition.  The other lifting and factorization
axioms follow by identifying the respective classes in terms of
$FI_{\partial }$ and $FI_{\Lambda }$.  By part~\ref{FILambda} of
Lemma~\ref{lem-level-properties}, an $FI_\Lambda$-cofibration is a
projective cofibration which is a level equivalence and an
$FI_\Lambda$-injective map is a level fibration.  Since $(J\cof,J\inj)$
has the lifting property for any class $J$, the lifting axiom for a map
that is both a level equivalence and a projective cofibration and a map
that is a level fibration follows by setting $J=FI_{\Lambda}$.
Moreover, every map can be factored as the composition of an
$FI_\Lambda$-cofibration and an $FI_\Lambda$-injective map.  Similarly,
every map can be factored as the composition of an
$FI_\partial$-cofibration and an $FI_\partial$-injective map, by
Lemma~\ref{lem-factorization}.  By part~\ref{FIpartial} of
Lemma~\ref{lem-level-properties}, an $FI_\partial$-cofibration is a
projective cofibration and an $FI_\partial$-injective map is a level
trivial fibration.

Now consider lifting and factorization for the injective level model
structure.  Here we use a set $C$ containing a map from each
isomorphism class of monomorphisms $i\mathcolon X\to Y$ with $Y$ a
countable symmetric spectrum, and a set $tC$ containing a map from each
isomorphism class of level trivial cofibrations $i\mathcolon X\to Y$
with $Y$ a countable symmetric spectrum.  The lifting axiom for a level
trivial cofibration $i$ and an injective fibration $p$ holds by
definition. By part~\ref{C} of
Lemma~\ref{lem-level-properties}, a $C$-cofibration is a level
cofibration and a $C$-injective map is an injective fibration that is a
level equivalence.  Since $(J\cof,J\inj)$ has the lifting property for
any class $J$, the lifting axiom for a level cofibration and a map that
is both an injective fibration and a level equivalence follows with
$J=C$.  Also, every map can be factored as the composition of a
$C$-cofibration followed by a $C$-injective map, by
Lemma~\ref{lem-factorization}.  Similarly, every map can be factored as
the composition of a $tC$-cofibration and a $tC$-injective map. By
part~\ref{tC} of Lemma~\ref{lem-level-properties}, a $tC$-cofibration is
a level trivial cofibration and a $tC$-injective map is an injective
fibration.
\end{proof}

\begin{corollary}\label{cor-enough-injectives}
Every symmetric spectrum embeds in an injective spectrum by a map that
is a level equivalence.
\end{corollary}

\begin{proof} 
For a symmetric spectrum $X$, the map $X\to *$ is the composition of a
level trivial cofibration $X\to E$ and an injective fibration $E\to
*$. The fibrant object $E$ is an injective spectrum.
\end{proof}

Some parts of the next lemma have already been proved. They are repeated
for easy reference.  Recall that $R_n \mathcolon \sset \to \spec $ is
the right adjoint of the evaluation functor $\Ev_n\mathcolon \spec\to
\sset$.

\begin{lemma}\label{lem-level-properties} 
\begin{enumerate}
\item[]
\item \label{RK} Let $K\subseteq \sset$ be the class of Kan fibrations
and let $RK=\cup_n R_nK$. Then a map is $RK$-projective if and only if
it is a level trivial cofibration.
\item \label{RtK} Let $tK\subseteq \sset$ be the class of trivial Kan
fibrations and let $R(tK)=\cup_n R_n(tK)$.  Then a map is
$R(tK)$-projective if and only if it is a level cofibration.
\item \label{FIpartial} Let $FI_\partial$ be the set defined in
\ref{def-level-generators}. Then a map is $FI_\partial$-injective if and
only if it is a level trivial fibration.  A map is an
$FI_\partial$-cofibration if and only if it is a projective cofibration.
\item \label{FILambda} Let $FI_\Lambda$ be the set defined in
\ref{def-level-generators}. Then a map is $FI_\Lambda$-injective if and
only if it is a level fibration.  A map is an $FI_\Lambda$-cofibration
if and only if it is a projective cofibration and a level equivalence.
\item \label{C} Let $C$ be a set containing a map from each
isomorphism class of monomorphisms $i\mathcolon X\to Y$ with $Y$ a
countable symmetric spectrum.  Then a map is $C$-injective if and only if
it is an injective fibration and a level equivalence.  A map is a
$C$-cofibration if and only if it is a level cofibration.  
\item \label{tC} Let $tC$ be a set containing a map from each
isomorphism class of level trivial cofibrations $i\mathcolon X\to Y$ with
$Y$ a countable symmetric spectrum.  Then a map is $tC$-injective if and
only if it is a injective fibration.  A map is a $tC$-cofibration if and
only if it is a level trivial cofibration.
\end{enumerate}
\end{lemma}

\begin{proof}
Parts~\ref{RK} and~\ref{RtK}: By adjunction, a map $g$ has the \llp
the class $RK$ ($R(tK)$) if and only if for each $n\ge0$ the map
$\Ev_ng$ has the \llp $K$ ($tK$). But $\Ev_ng$ has the \llp $K$ ($tK$)
if and only if $\Ev_ng$ is a trivial cofibration (arbitrary
cofibration), \ie if and only if $g$ is a level trivial cofibration
(arbitrary level cofibration).

Part~\ref{FIpartial} is proved in 
Propositions~\ref{prop-level-fib} and ~\ref{prop-stable-cof}. 

Part~\ref{FILambda}: The first claim is proved in
Proposition~\ref{prop-level-fib}. Every $FI_{\Lambda }$-cofibration has
the \llp level fibrations, so is in particular a projective cofibration
by Part~\ref{FIpartial}.  Every map in $FI_\Lambda$ is a level trivial
cofibration by Proposition~\ref{prop-F-preserves}, so is $RK$-projective
by Part~\ref{RK}.  So every $FI_{\Lambda}$-cofibration
is also $RK$-projective and hence is a level trivial cofibration by 
Part~\ref{RK} again.  So in particular it is a level equivalence.  

Conversely, suppose $f$ is a projective cofibration and a level
equivalence.  We can factor $f$ as the composition of an $FI_\Lambda
$-cofibration $i$ and an $FI_\Lambda $-injective map $p$, by
Lemma~\ref{lem-factorization}.  By the two-out-of-three property, $p$ is
a level equivalence. Therefore the projective cofibration $f$ has the
\llp the level trivial fibration $p$. By the Retract
Argument~\ref{prop-retract}, $f$ is a retract of $i$, and so is an
$FI_\Lambda$-cofibration.

For part~\ref{C}, first note that, by part~\ref{RtK}, every
$C$-cofibration is a level cofibration.  Conversely, suppose
$f\mathcolon X\xrightarrow{}Y$ is a level cofibration.  Then $f$ is a
$C$-cofibration if, for every $C$-injective map $g$ and commutative
square 
\[
\begin{CD}
X @>>> E \\
@VfVV  @VVgV \\
Y @>>> Z
\end{CD}
\]
there is a lift $h\mathcolon Y\to E$ making the diagram commute. Let
$\cp$ be the partially ordered set of partial lifts: an object of $\cp$
is a pair $(U,h_U)$ such that $X\subseteq U \subseteq Y$ and the diagram
\[
\xymatrix{
X \ar[d]_{i_U} \ar[r]&  E \ar[d]^g \\
  U \ar[r] \ar[ur]^{h_U} & Z,
}
\]
is commutative.  We define $(U,h_U)\le (V,h_V)$ if $U\subseteq V$ and
$h_V$ extends $h_U$.  Every chain in $\cp$ has an upper bound and so
Zorn's lemma gives a maximum $(M,h_M)$.  Suppose $M$ is strictly
contained in $Y$.  Then, by taking the subspectrum generated by a
simplex not in $M$, we find a countable subspectrum $D$ (by
Lemma~\ref{lem-countable-subspectra} below) such that the level
cofibration $D\cap M\to D$ is not an isomorphism.  By construction, the
map $D\cap M\to D$ is isomorphic to a map in $C$.  By cobase change,
$M\to D\cup M$ is a $C$-cofibration.  Thus $h_{M}$ extends to a partial
lift on $D\cup M$, contradicting the maximality of $(M,h_{M})$.
Therefore $M=Y$, and so $f$ is a $C$-cofibration.

We now identify $C\inj$. Since $(C\cof)\inj=C\inj$, every
$C$-injective map has the \rlp every monomorphism. In particular,
every $C$-injective map is an injective fibration. Let $f\mathcolon E\to
B$ be a map having the \rlp every monomorphism. Let $s\mathcolon
B\to E$ be a lift in the diagram 
\[
\begin{CD}
* @>>> E \\
@VVV @VVfV \\
B @= B
\end{CD}
\] 
Then $fs$ is the identity map on $B$.  To study the composite $sf$, let
$j$ be the monomorphism $\partial\Delta[1]\to\Delta[1]$. The diagram
\[
\begin{CD}
E\vee E @>sf\vee 1>> E\\
@VE\sm j_+VV    @VVfV\\
E\sm\Delta[1]_+ @>>> B
\end{CD}
\]
is commutative since $fsf=f$ and has a lift since $E\sm j_+$ is a
monomorphism. The lift is a simplicial homotopy from $sf$ to the
identity on $E$. Therefore $f$ is a simplicial homotopy equivalence and
in particular $f$ is a level equivalence.  Conversely suppose $f$ is
both an injective fibration and a level equivalence.  We can factor $f$
as the composition of a $C$-cofibration $i$ and a $C$-injective map
$p$. By the two-out-of-three property, $i$ is a level equivalence. The level
trivial cofibration $i$ has the \llp the injective fibration $f$.  By the
Retract Argument~\ref{prop-retract}, $f$ is a retract of $p$ and so is a
$C$-injective map.

The proof of part~\ref{tC} is similar, though slightly more complex.  By
part~\ref{RK}, every $tC$-cofibration is a level trivial
cofibration. Conversely, suppose $f\mathcolon X\to Y$ is a level trivial
cofibration. Then $f$ is a $tC$-cofibration if, for every $tC$-injective
map $g$ and commutative square
\[
\begin{CD}
X @>>> E \\
@VfVV  @VVgV \\
Y @>>> Z
\end{CD}
\]
there is a lift $h\mathcolon Y\to E$ making the diagram commute.  We
again let $\cp$ be the partially ordered set of partial lifts: an object
of $\cp$ is a pair $(U,h_U)$ such that $X\subseteq U \subseteq Y$, and
the diagram
\[
\xymatrix{
X \ar[d]_{i_U} \ar[r]&  E \ar[d]^g \\
  U \ar[r] \ar[ur]^{h_U} & Z,
}
\]
is commutative, but we also require that the inclusion $i_U\mathcolon
X\to U$ is a weak equivalence.  We define $(U,h_U)\le (V,h_V)$ as
before.  Every chain in $\cp$ has an upper bound (using the fact that a
transfinite composition of level trivial cofibrations is a level trivial
cofibration) and so Zorn's lemma gives a maximum $(M,h_M)$.  The inclusion
$X\to M$ is a level trivial cofibration, so, by the two-out-of-three
property, the inclusion $M\to Y$ is a weak equivalence.  If $M$ is
strictly contained in $Y$, Lemma~\ref{lem-countable-equivalence}, proved
below, applied to the countable subspectrum of $Y$ generated by a
simplex not in $M$, gives a countable subspectrum $D$ of $Y$ such that the
monomorphism $D\cap M\to D$ is a weak equivalence but is not an
isomorphism. By construction, $D\cap M\to D$ is isomorphic to a map in
$tC$.  By cobase change, $M\to D\cup M$ is a $tC$-cofibration. So $h_M$
extends to a partial lift on $D\cup M$.  This is a contradiction since
$(M,h_M)$ is maximal.  Thus $M=Y$, and so $f$ is a $tC$-cofibration.

Since $(tC\cof)\inj=tC\inj$, the $tC$-injective maps are the injective
fibrations.
\end{proof}

\begin{corollary}
Every injective fibration is a level fibration and every projective
cofibration is a level cofibration.
\end{corollary}

\begin{proof}
By Proposition~\ref{prop-F-preserves}, every map in $FI_\partial $ is a
level cofibration. Therefore, by part~\ref{RtK} of
Lemma~\ref{lem-level-properties}, every projective cofibration is a
level cofibration. By Proposition~\ref{prop-F-preserves}, every map in
$FI_\Lambda$ is a level trivial cofibration. Therefore every injective
fibration is a level fibration.
\end{proof}

The following lemmas are used in the proof of
Lemma~\ref{lem-level-properties}.

\begin{lemma}\label{lem-countable-subspectra}
Let $X$ be a spectrum, and suppose $x$ is a simplex of $X_{n}$ for some
$n\geq 0$.  Then the smallest subspectrum of $X$ containing $x$ is
countable.  
\end{lemma}

\begin{proof}
First note that if $L$ is a countable collection of simplices in a
simplicial set $K$, then the smallest subsimplicial set of $K$
containing $L$ is also countable.  Indeed, we need only include all
degeneracies of all faces of simplices in $L$, of which there are only
countably many.  Similarly, if $L$ is a countable collection of
simplices in a $\Sigma _{n}$-simplicial set $K$, then the smallest
sub-$\Sigma _{n}$-simplicial set containing $L$ is countable.  Indeed,
we only need to include the orbits of all degeneracies of all faces of
simplices in $L$.  

Now, let $Y_{n}$ denote the sub-$\Sigma _{n}$-simplicial set of $X_{n}$
generated by $x$.  We have just seen that $Y_{n}$ is countable.  We
then inductively define $Y_{n+k}$ to be the smallest sub-$\Sigma
_{n+k}$-simplicial set of $X_{n+k}$ containing the image of $S^{1}\sm
Y_{n+k-1}$.  Then each $Y_{n+k}$ is countable, and the $Y_{n+k}$ define
a subspectrum of $X$ containing $x$.  
\end{proof}

It follows in similar fashion that the smallest subspectrum of a
spectrum $X$ containing any countable collection of simplices of $X$ is
countable.

We need a similar lemma for inclusions which are level equivalences.  To
prove such a lemma, we need to recall from the comments before
Lemma~\ref{lem-comparison-arbitrary} that homotopy of simplicial sets
commutes with transfinite compositions of monomorphisms.  The same methods
imply that relative homotopy commutes with transfinite compositions of
monomorphisms.

\begin{lemma}\label{lem-countable-equivalence}
Let $f\mathcolon X\to Y$ be a level trivial cofibration of symmetric
spectra. For every countable subspectrum $C$ of $Y$ there is a
countable subspectrum $D$ of $Y$ such that $C\subseteq D$ and $D\cap
X\to D$ is a level trivial cofibration.
\end{lemma}

\begin{proof}
Let $K\subseteq L$ be a pair of pointed simplicial sets and $v$ be a
$0$-simplex of $K$. For $n\ge 1$, let $\pi_n(L,K;v)$ denote the
relative homotopy set of the pair with the null element as the
basepoint (ignore the group structure when $n\ge 2$).  To ease
notation let $\pi_0(L,K;v)$ be the pointed set $\pi_0L/\pi_0K$.  The
inclusion $K\to L$ is a weak equivalence if and only if
$\pi_n(L,K;v)=*$ for every $v\in K_0$ and $n\ge0$.

Now, construct a countable spectrum $FC$ such that the map
$\pi_*(C_n,C_n\cap X_n;v)\to \pi_*(FC_n,FC_n\cap X_n;v)$ factors
through the basepoint $*$ for every $0$-simplex $v$ of $C_n\cap X_n$
and integer $n\ge0$.  Since $\pi_*(Y,X;v)=*$ and $\pi_*$ commutes with
filtered colimits, for each homotopy class $\alpha\in\pi_*(C_n,
C_n\cap X_n;v)$ there is a finite simplicial subset $K_\alpha\subseteq
Y_n$ such that $\pi_*(C_n,C_n\cap X_n;v)\to \pi_*(K_\alpha\cup
C_n,(K_\alpha\cup C_n)\cap X_n;v)$ sends $\alpha$ to the basepoint.
Since $C_n$ is countable, the set $\pi_*(C_n, C_n\cap X_n;v)$ is
countable. Let $B_n$ be the union of $C_n$ with all the finite
simplicial sets $K_\alpha$.  The $B_n$ are countable simplicial sets 
and generate a countable subspectrum $FC$ of $Y$ having the desired
property.

Repeat the construction to get a sequence of countable subspectra of
$Y$:
\[
C \to FC \to F^2C \to \dots \to F^nC\to \dots
\]
Let $D=\colim_nF^nC$.  The spectrum $D$ is countable. Since relative
homotopy commutes with transfinite compositions of monomorphisms, the set
$\pi_*(D_n,D_n\cap X_n;v)$ has only one element.  Therefore the
inclusion $D_n\cap X_n\to D_n$ is a weak equivalence, and so $D\cap X\to
D$ is a level equivalence.
\end{proof}

\subsection{Stable cofibrations}\label{subsec-stable-cofib}

The object of this section is to give a characterization of stable
cofibrations in $\spec $.  To this end, we introduce the latching
space.

\begin{definition}\label{def-latching}
Define $\ov{S}$ to be the symmetric spectrum such that
$\ov{S}_{n}=S^{n}$ for $n>0$ and $\ov{S}_{0}=*$.  The
structure maps are the evident ones.  Given a symmetric spectrum $X$,
define the $n$th \emph{latching space}, $L_{n}X$, to be
$\Ev _{n}(X\sm \ov{S})$.   
\end{definition}

There is a map of symmetric spectra $i\mathcolon \ov{S}\to S$
which is the identity on positive levels.  This induces a natural
transformation $L_{n}X\to X_{n}$ of pointed $\Sigma _{n}$ simplicial
sets.

The following proposition uses a model structure on the category of
pointed $\Sigma _{n}$ simplicial sets.  A map $f\mathcolon X\to Y$ of
pointed $\Sigma _{n}$ simplicial sets is a \emph{$\Sigma
_{n}$-fibration} if it is a Kan fibration of the underlying simplicial
sets.  Similarly, $f$ is a weak equivalence if it is a weak equivalence
of the underlying simplicial sets.  The map $f$ is a \emph{$\Sigma
_{n}$-cofibration} if it is a monomorphism such that $\Sigma _{n}$ acts
freely on the simplices of $Y$ not in the image of $f$.  It is
well-known, and easy to check, that the $\Sigma _{n}$-cofibrations, the
$\Sigma _{n}$-fibrations, and the weak equivalences define a model
structure on the category of pointed $\Sigma _{n}$-simplicial sets.  

\begin{proposition}\label{prop-stable-cofib}
A map $f\mathcolon X\to Y$ in $\spec $ is a stable cofibration if and
only if for all $n\geq 0$ the induced map $\Ev _{n}(f\boxprod
i)\mathcolon X_{n}\amalg _{L_{n}X}L_{n}Y\to Y_{n}$ is a $\Sigma
_{n}$-cofibration.  
\end{proposition}

\begin{proof}
Suppose first that $\Ev _{n}(f\boxprod i)$ is a $\Sigma
_{n}$-cofibration for all $n$.  Suppose $g\mathcolon E\to B$ is a level
trivial fibration.  We show that $f$ has the \llp $g$ by constructing a
lift using induction on $n$.  A partial lift defines a commutative square
\[
\begin{CD}
X_{n}\amalg_{L_{n}X}L_{n}Y @>>> E_{n} \\
@VVV @VVV \\
Y_{n} @>>> B_{n}
\end{CD}
\]
Since the left vertical map is a $\Sigma _{n}$-cofibration and the right
vertical map is a trivial $\Sigma _{n}$-fibration, there is a lift in
this diagram and so we can extend our partial lift.  Hence $f$ has the
\llp $g$, and so $f$ is a stable cofibration.  

To prove the converse, note that $\Ev _{n}$ is a left adjoint as a
functor to pointed $\Sigma _{n}$ simplicial sets.  Since the class of
stable cofibrations is the class $FI_{\partial }\cof$, it suffices to
check that $\Ev _{n}(f\boxprod i)$ is a $\Sigma _{n}$-cofibration for
$f\in FI_{\partial }$.  More generally, suppose $g\mathcolon A\to B$ is
a monomorphism of pointed simplicial sets.  Since $F_{m}g=S\otimes
G_{m}g$, we have $F_{m}g\boxprod i=G_{m}g\boxprod i$, where the second
$\boxprod $ is taken in $\symseq $.  One can easily check that $\Ev
_{n}(G_{m}g\boxprod i)$ is an isomorphism when $n\neq m$ and is the map
$(\Sigma _{n})_{+}\sm g$ when $n=m$.  In both cases, $\Ev
_{n}(G_{m}g\boxprod i)$ is a $\Sigma _{n}$-cofibration, as required. 
\end{proof}

\subsection{Pushout smash product}\label{subsec-pushout-smash}

In this section we consider the pushout smash product in an arbitrary
symmetric monoidal category and apply our results to $\spec $.  We show
that the projective level structure and the stable model structure on
$\spec $ are both compatible with the symmetric monoidal structure.  A
monoid $E$ in $\spec $ is called a \emph{symmetric ring spectrum}, and
is similar to an $A_{\infty }$-ring spectrum.  Thus, there should be a
stable model structure on the category of $E$-modules.  Similarly, there
should be a model structure on the category of symmetric ring spectra
and the category of commutative symmetric ring spectra.  These issues
are dealt with more fully in~\cite{schwede-shipley-monoids} and in work
in progress of the third author.  Their work depends heavily on the
results in this section and in Section~\ref{subsec-monoids}.  The
results of this section alone suffice to construct a stable model
structure on the category of modules over a symmetric ring spectrum
which is stably cofibrant.  This section also contains brief
descriptions of two other stable model structures on $\spec $.  

\begin{definition}
Let $f\mathcolon U\to V$ and $g\mathcolon X\to Y$ be maps in a symmetric
monoidal category $\cc$. The \emph{pushout smash product}
\[
f\boxprod g\mathcolon V\sm X\amalg_{U\sm X}U\sm Y\to V\sm Y.
\]
is the natural map on the pushout defined by the commutative square
\[
\begin{CD}
U\sm X @>f\sm X>> V\sm X \\ @VU\sm gVV @VV{V\sm g}V \\ U\sm Y @>>f\sm
Y> V\sm Y.
\end{CD}
\]
If $\cc$ is a closed symmetric monoidal category, 
\[
\ihom_{\square}(f,g)\mathcolon\ihom(V,X)\to
\ihom(U,X)\times_{\ihom(U,Y)}\ihom(V,Y). 
\]
is the natural map to the fiber product defined by the commutative
square
\[
\begin{CD}
\ihom(V,X) @>f^*>> \ihom(U,X)\\ @Vg_*VV @VVg_*V\\ \ihom(V,Y) @>>f^*> 
\ihom(U,Y).
\end{CD}
\]
\end{definition}

\begin{definition}
A model structure on a symmetric monoidal category is called
\emph{monoidal} if the pushout smash product $f\boxprod g$ of two
cofibrations $f$ and $g$ is a cofibration which is trivial if either $f$
or $g$ is.  
\end{definition}

In our situation, this is the correct condition to require
so that the model structure is compatible with the symmetric monoidal
structure.  Since the unit, $S$, is cofibrant in symmetric spectra this
condition also ensures that the symmetric monoidal structure induces a
symmetric monoidal structure on the homotopy category.  For a more general 
discussion of monoidal model structures, see~\cite{hovey-model}.  

Recall, from Definition~\ref{def-boxHom-category}, the map of sets $\cc
_{\square}(f,g)$.  

\begin{proposition}\label{prop-adjunctions-symmetric}
Let $f,g$ and $h$ be maps in a closed symmetric monoidal category
$\cc$. There is a natural isomorphism
\[
\cc_\square(f\boxprod g,h)\natiso\cc_\square(f,\ihom_\square(g,h))
\]
\end{proposition}

\begin{proof}
Use the argument in the proof of
Proposition~\ref{prop-pushout-smash-adjunctions}.
\end{proof}

\begin{proposition}\label{prop-check-generators}
Let $I$ and $J$ be classes of maps in a closed symmetric monoidal
category $\cc$.  Then
\[
I\cof\boxprod J\cof\subseteq (I\boxprod J)\cof
\]
\end{proposition}

\begin{proof}
Let $K=(I\boxprod J)\inj$. By hypothesis, $(I\boxprod J,K)$ has the
lifting property. By adjunction, $(I,\ihom_\square(J,K))$ has the
lifting property. By Proposition~\ref{prop-lifting-properties},
$(I\cof,\ihom_\square(J,K))$ has the lifting property. Then
$(J,\ihom_\square(I\cof,K))$ has the lifting property, by using
adjunction twice.  Thus $(J\cof,\ihom_\square(I\cof,K))$ has the lifting
property, by Proposition~\ref{prop-lifting-properties}. By adjunction,
$(I\cof\boxprod J\cof,K)$ has the lifting property. So $I\cof\boxprod
J\cof\subseteq (I\boxprod J)\cof$ and the proposition is proved.
\end{proof}

\begin{corollary}\label{cor-check-generators}
For classes $I$, $J$ and $K$ in $\spec$, if $I\boxprod J\subseteq K\cof$
then $I\cof\boxprod J\cof\subseteq K\cof$.
\end{corollary}

We now examine to what extent the pushout smash product preserves stable
cofibrations and stable equivalences.  To do so, we introduce a new
class of maps in $\spec $. 

\begin{definition}\label{defn-S}
Let $M$ be the class of monomorphisms in the category of symmetric
sequences $\symseq$.  A map of symmetric spectra is an
\emph{$S$-cofibration} if it is an $S\otimes M$-cofibration. A symmetric
spectrum $X$ is \emph{$S$-cofibrant} if the map $*\to X$ is an
$S$-cofibration.  A map is an \emph{$S$-fibration} if it has the \rlp
every map which is both an $S$-cofibration and a stable equivalence. 
\end{definition} 

Note that every stable cofibration is an $S$-cofibration, since
$FI_{\partial }= S\otimes \bigcup _{n}G_{n}I_{\partial }$.  On the other
hand, by Proposition~\ref{prop-otimes-preserves}, every element of
$S\otimes M$ is a monomorphism, and so every $S$-cofibration is a level
cofibration.  There is a model structure on $\spec $, called the
\emph{$S$ model structure}, where the cofibrations are the
$S$-cofibrations and the weak equivalences are the stable equivalences.
The fibrations, called \emph{$S$-fibrations} are those maps with the
\rlp $S$-cofibrations which are also stable equivalences.  Every
$S$-fibration is a stable fibration.  This model structure will be used
in a forthcoming paper by the third author to put a model structure on
certain commutative $S$-algebras.

We mention as well that there is a third model structure on $\spec $
where the weak equivalences are the stable equivalences, called the
\emph{injective \textup{(}stable\textup{)} model structure}.  The
injective cofibrations are the level cofibrations and the injective
stable fibrations are all maps which are both injective fibrations and
stable fibrations.  In particular, the fibrant objects are the injective
$\Omega $-spectra.  The interested reader can prove this is a model
structure using the methods of Section~\ref{subsec-stable-model-cat},
replacing the set $I$ with the union of $I$ and the countable level
cofibrations.

\begin{theorem}\label{thm-boxprod}
Let $f$ and $g$ be maps of symmetric spectra.
\begin{enumerate}
\item If $f$ and $g$ are stable cofibrations then $f\boxprod g$ is a
stable cofibration.
\item If $f$ and $g$ are $S$-cofibrations then $f\boxprod g$ is an
$S$-cofibration.
\item If $f$ is an $S$-cofibration and $g$ is a level cofibration, then
$f\boxprod g$ is a level cofibration.
\item If $f$ is an $S$-cofibration, $g$ is a level cofibration, and
either $f$ or $g$ is a level equivalence, then $f\boxprod g$ is a level
equivalence.
\item If $f$ is an $S$-cofibration, $g$ is a level cofibration, and
either $f$ or $g$ is a stable equivalence, then $f\boxprod g$ is a stable
equivalence.
\end{enumerate}
\end{theorem}

\begin{proof}
Parts 1 through 4 of the Proposition are proved using
Corollary~\ref{cor-check-generators}.

Part 1: Let $I=J=K=FI_\partial$. Then $K\cof$ is the class of stable
cofibrations.  We have a natural isomorphism
\[
F_pf\boxprod F_qg=F_{p+q}(f\boxprod g)
\]
for $f,g\in\sset$.  By Proposition~\ref{prop-SM7-sset}, $f\boxprod g$ is
a monomorphism when $f$ and $g$ are.  Part~3 of
Proposition~\ref{prop-stable-cof} the shows that $I\boxprod J\subseteq
K\cof$.  Now use the corollary.

Part 2: Let $I=J=K=S\otimes M$ (recall that $M$ is the class of
monomorphisms in $\symseq$). By definition, $K\cof$ is the class of
$S$-cofibrations. For $f$ and $g$ in $\symseq $, we have a natural
isomorphism
\[
S\otimes f\boxprod S\otimes g=S\otimes(f\boxprod g)
\]
where the first $\boxprod $ is taken in $\spec $ and the second
$\boxprod $ is taken in $\symseq $.  In degree $n$, 
\[
(f\boxprod g)_{n}=\bigvee _{p+q=n} (\Sigma _{p+q})_{+}\sm _{\Sigma
_{p}\times \Sigma _{q}}(f_{p}\boxprod g_{q})
\]
For $f, g \in M$, each map $f_{p}\boxprod g_{q}$ is a monomorphism, so
it follows that $f\boxprod g$ is a monomorphism.  Thus $I\boxprod
J\subseteq K\cof$. Now use the corollary.

Part 3: Let $I=S\otimes M$ and let $J=K$ be the class of level
cofibrations. By Part~\ref{RtK} of Lemma~\ref{lem-level-properties},
$K\cof=K$.  For $f\in \symseq $ and $g\in \spec $, we have a natural
isomorphism of maps of symmetric sequences
\[
(S\otimes f)\boxprod g=f\boxprod g
\]
where the first $\boxprod $ is taken in $\spec $ and the second
$\boxprod$ is taken in $\symseq$.  We have seen in the proof of part~2
that $f\boxprod g$ is a monomorphism of symmetric sequences if $f$ and
$g$ are monomorphisms.  Hence $I\boxprod J\subseteq K\cof$. Now use the
corollary.

Part 4: First assume $g$ is a level trivial cofibration. Let
$I=S\otimes M$ and let $J=K$ be the class of level trivial
cofibrations. By Part~\ref{RK} of Lemma~\ref{lem-level-properties},
$K=K\cof$.  Proposition~\ref{prop-SM7-sset} and the method used in the
proof of part~2 imply that, if $f$ and $g$ are monomorphisms of
symmetric sequences and $g$ is a level equivalence, then $f\boxprod g$
is a level equivalence.  Recall that, for $h\in \symseq $ and $g\in
\spec $, we have a natural isomorphism of maps of symmetric sequences
\[
(S\otimes h)\boxprod g=h\boxprod g
\]
where the first $\boxprod $ is taken in $\spec $ and the second
$\boxprod$ is in $\symseq$.  Hence $I\boxprod J\subseteq K\cof$. Now
use the corollary to prove part~4 when $g$ is a level equivalence.

It follows that, for any injective spectrum
$E$ and an arbitrary $S$-cofibration $h$, the map $\ihom_S(h,E)$ is an
injective fibration.  Indeed, if $g$ is a level cofibration and a level
equivalence, $\spec (g,\ihom _{S}(h,E))\natiso \spec (g\boxprod h,E)$,
and we have just seen that $g\boxprod h$ is a level cofibration and a
level equivalence, so $(g\boxprod h,E)$ has the lifting property.  

Now suppose $f$ is an $S$-cofibration and a level equivalence.  Then the
map $\ihom_{\spec}(f,E)$ is an injective fibration and a level
equivalence.  Indeed, we have $\Ev _{k}\ihom _{\spec }(f,E)=\Map _{\spec
}(f\sm (S\otimes \Sigma [k]_{+}),E)$, by
Remark~\ref{rem-describe-internal-hom}.  Since $S\otimes \Sigma [k]_{+}$
is $S$-cofibrant, and $f$ is both a level equivalence and a level
cofibration, we have just proved that $f\sm (S\otimes \Sigma
[k]_{+})=f\boxprod (*\to S\otimes \Sigma [k]_{+})$ is a level
equivalence and a level cofibration.  This shows that $\pi _{0}\Ev
_{k}\ihom _{\spec }(f,E)$ is an isomorphism; smashing with $F_{0}S^{n}$
and using a similar argument shows that $\pi _{n}\Ev _{k}\ihom _{\spec
}(f,E)$ is an isomorphism.  

Thus every level cofibration $g$ has the \llp the map
$\ihom_{\spec}(f,E)$.  By adjunction, $f\boxprod g$ and $f\boxprod
(g\boxprod j)$ where $j\mathcolon \partial\Delta[1]_+\to\Delta[1]_+$ is the
inclusion, have the extension property with respect to every injective
spectrum $E$.  It follows that $E^0(f\boxprod g)$ is an isomorphism for
every injective spectrum $E$ and hence that $f\boxprod g$ is a level
equivalence.

Part 5: Because we are working in a stable situation, a level
cofibration $i\mathcolon X\to Y$ is a stable equivalence if and only if
its cofiber $C_i=Y/X$ is stably trivial. The map $f\boxprod g$ is a
level cofibration by part~(3).  By commuting colimits, the cofiber of
$f\boxprod g$ is the smash product $Cf\sm Cg$ of the cofiber $Cf$
of $f$, which is $S$-cofibrant, and the cofiber $Cg$ of $g$.  Let $E$
be an injective $\Omega$-spectrum. We will show that $\ihom
_{S}(Cf\sm Cg,E)$ is a level trivial spectrum, and thus that
$Cf\sm Cg$ is stably trivial.  

First suppose that $f$ is a stable equivalence.  Then $\ihom_S(Cf,E)$
is a level trivial spectrum which is also injective, by part~4 and the
fact that $Cf$ is $S$-cofibrant.  Therefore $\ihom _{S}(Cf\sm
Cg,E)\natiso \ihom_S(Cg,\ihom_S(Cf,E))$ is a level trivial
spectrum, so $Cf\sm Cg$ is stably trivial and thus $f\boxprod g$ is a
stable equivalence.

Now suppose that $g$ is a stable equivalence, so that $Cg$ is stably
trivial.  By adjunction $\ihom_S(Cf\sm Cg,
E)=\ihom_S(Cg,\ihom_S(Cf,E))$.  We claim that $D=\ihom _{S}(Cf,E)$
is an injective $\Omega $-spectrum.  Indeed, we have already seen that
$D$ is injective.  From Remark~\ref{rem-describe-internal-hom}, we have
\[
\Ev _{n}D\natiso \shom _{\spec }(Cf\sm F_{n}S^{0},E) \natiso \shom
_{\spec } (Cf,\ihom _{S}(F_{n}S^{0},E))
\]
Similarly, we have 
\[
(\Ev _{n+1}D)^{S^{1}}\natiso \shom _{\spec }(Cf,\ihom
_{S}(F_{n+1}S^{1},E)).
\]
Since $E$ is an $\Omega $-spectrum, $(F_nS^0 \sm \badmap)^*\mathcolon
\ihom _{S}(F_{n}S^{0},E)\xrightarrow{} \ihom _{S}(F_{n+1}S^{1},E)$
is a level equivalence.  Since $E$ is injective, both the source and 
target are injective, and so this map is a simplicial homotopy
equivalence by Lemma~\ref{lem-injective-spectra}.  Hence $\Ev
_{n}D\xrightarrow{}(\Ev_{n+1}D)^{S^{1}}$ is still a level equivalence,
so $D=\ihom _{S}(Cf,E)$ is an injective $\Omega $-spectrum.  Since
$Cg$ is stably trivial, $\ihom _{S}(Cf\sm Cg,E)\natiso
\ihom_S(Cg,\ihom_S(Cf,E))$ is a level trivial spectrum, so $Cf\sm
Cg$ is stably trivial and thus $f\boxprod g$ is a stable
equivalence.  
\end{proof}

\begin{corollary}
The projective model structure and the stable model structure on $\spec
$ are monoidal.  
\end{corollary}

It also follows that the $S$ model structure on $\spec $ is monoidal,
once it is proven to be a model structure.  Neither the injective level
structure nor the injective stable structure is monoidal.  

Adjunction then gives the following corollary.  

\begin{corollary}\label{cor-Hom-SM7}
Let $f$ and $g$ be maps of symmetric spectra. 
\begin{enumerate}
\item If $f$ is a stable cofibration and $g$ is a stable fibration, then
$\Hom _{\square}(f,g)$ is a stable fibration, which is a level
equivalence if either $f$ or $g$ is a stable equivalence.  
\item If $f$ is a stable cofibration and $g$ is a level fibration, then
$\Hom _{\square}(f,g)$ is a level fibration, which is a level
equivalence if either $f$ or $g$ is a level equivalence.  
\item If $f$ is an $S$-cofibration and $g$ is an $S$-fibration, then
$\Hom _{\square}(f,g)$ is an $S$-fibration, which is a level
equivalence if either $f$ or $g$ is a stable equivalence.  
\item If $f$ is an $S$-cofibration and $g$ is an injective fibration,
then $\Hom _{\square}(f,g)$ is an injective fibration, which is a level
equivalence if either $f$ or $g$ is a level equivalence.  
\end{enumerate}
\end{corollary}

\begin{corollary}\label{cor-X-smash}
If $X$ is an $S$-cofibrant symmetric spectrum, the functor $X\sm~-$
preserves level equivalences and it preserves stable equivalences.
\end{corollary}

\begin{proof}
Part four of Theorem~\ref{thm-boxprod} implies that $X\sm~-$ preserves
level trivial cofibrations.  Lemma~\ref{lem-ken-brown} then implies that
it preserves level equivalences, since every symmetric spectrum is level
cofibrant.  An arbitrary stable equivalence can be factored as a stable
trivial cofibration followed by a level equivalence.  Part five of
Theorem~\ref{thm-boxprod} implies that $X\sm~-$ takes stable trivial
cofibrations to stable equivalences.  
\end{proof}

\subsection{Proper model categories}\label{subsec-proper}

In this section we recall the definition of a proper model category and
show that the stable model category of symmetric spectra is proper.

\begin{definition}
\begin{enumerate}
\item
A model category is \emph{left proper} if for every pushout square
\[
\begin{CD}
A @>g>> B\\
@VfVV   @VVhV\\
X @>>>  Y
\end{CD}
\]
with $g$ a cofibration and $f$ a weak equivalence, $h$ is a weak
equivalence.
\item
A model category is \emph{right proper} if for every pullback square
\[
\begin{CD}
A @>>> B\\
@VhVV   @VVfV\\
X @>>g>  Y
\end{CD}
\]
with $g$ a fibration and $f$ a weak equivalence, $h$ is a weak
equivalence.
\item
A model category is \emph{proper} if it is both left proper and right
proper.  
\end{enumerate}
\end{definition}

The category of simplicial sets is a proper model
category~\cite{bousfield-friedlander} (see~\cite{hirschhorn} for more
details).  Hence the category of pointed simplicial sets and both level
model structures on $\spec $ are proper.

\begin{theorem}\label{thm-proper}
The stable model category of symmetric spectra is proper.
\end{theorem}

\begin{proof}
Since every stable cofibration is a level cofibration, the stable
model category of symmetric spectra is left proper by part one of
Lemma~\ref{lem-proper}. Since every stable fibration is a level
fibration, the stable model category of symmetric spectra is right
proper by part two of Lemma~\ref{lem-proper}.
\end{proof}

\begin{lemma}\label{lem-proper}
\begin{enumerate}
\item[]
\item Let 
\[
\begin{CD}
A @>g>> B\\
@VfVV   @VVhV\\
X @>>>  Y
\end{CD}
\]
be a pushout square with $g$ a level cofibration and $f$  a stable
equivalence.  Then $h$ is a stable equivalence.
\item 
Let 
\[
\begin{CD}
A @>k>> B\\
@VfVV   @VVhV\\
X @>>g>  Y
\end{CD}
\]
be a pullback square with $g$ a level fibration and $h$ a stable
equivalence.  Then $f$ is a stable equivalence.
\end{enumerate}
\end{lemma}

\begin{proof}
Part one: Let $E$ be an injective $\Omega$-spectrum.  Apply the functor
$\shom_{\spec }(-,E)$ to the pushout square. The resulting commutative
square
\[
\begin{CD}
\shom(A,E) @<\shom(g,E)<< \shom(B,E)\\
@A\shom(f,E)AA   @AA\shom(h,E)A\\
\shom(X,E) @<<<  \shom(Y,E)
\end{CD}
\]
is a pullback square of pointed simplicial sets with $\shom(f,E)$ a weak
equivalence, by Proposition~\ref{prop-def-stable-equivalence}.  We claim
that $\shom(g,E)$ is a Kan fibration.  Indeed, let $k\mathcolon E\to *$
denote the obvious map, and let $c$ denote a trivial cofibration of
pointed simplicial sets.  Then $\shom (g,E)=\shom _{\square}(g,k)$.  We
must show that $(c,\shom_{\square}(g,k))$ has the lifting property.  By
Corollary~\ref{cor-boxhom-lifting}, this is equivalent to showing that
$(c\boxprod g,k)$ has the lifting property.  But, by
Proposition~\ref{prop-SM7-sset}, $c\boxprod g$ is a level equivalence
and a level cofibration.  Since $E$ is injective, it follows that
$(c\boxprod g,k)$ has the lifting property, and so $\shom(g,E)$ is a Kan
fibration.  By properness for simplicial sets, $\shom(h,E)$ is a weak
equivalence. It follows that $h$ is a stable equivalence.

Part two: Let $F$ be the fiber over the basepoint of the map
$g\mathcolon X\to Y$. Since $k$ is a pullback of $g$, $F$ is isomorphic to
the fiber over the basepoint of the map $k\mathcolon A\to B$. The maps
$X/F\to Y$, $A/F\to B$ and $B\to Y$ are stable equivalences; so
$A/F\to X/F$ is a stable equivalence. Consider the Barratt-Puppe
sequence (considered in the proof of
Lemma~\ref{stable-fib-equiv-level}) 
\[
F \to A \to A/F \to F\sm S^{1} \to A\sm S^{1} \to A/F \sm S^{1} \to F
\sm S^{2} 
\]
and the analogous sequence for the pair $(F,X)$.  Given an injective
$\Omega $-spectrum $E$, apply the functor $E^{0}(-)$, and note that
$E^{0}(Z\sm S^{1})\natiso \pi _{1}\shom (Z,E)$ is naturally a group.
The five-lemma then implies that $f\sm S^{1}\mathcolon A\sm
S^{1}\xrightarrow{}X\sm S^{1}$ is a stable equivalence.  Part two of
Theorem~\ref{thm-traditional-results} shows that $f$ is a stable
equivalence.
\end{proof}

\subsection{The monoid axiom}\label{subsec-monoids}

In~\cite{schwede-shipley-monoids}, techniques are developed to form
model category structures for categories of monoids, algebras, and
modules over a monoidal model category.  One more axiom is required
which is referred to as the \emph{monoid axiom}.  In this section we
verify the monoid axiom for symmetric spectra.  The results
of~\cite{schwede-shipley-monoids} then immediately give a model
structure on symmetric ring spectra.  See also~\cite{shipley-thh} for
more about symmetric ring spectra.  After proving the monoid axiom, we
discuss the homotopy invariance of the resulting model categories of
modules and algebras.  

Let $K$ denote the class in $\spec $ consisting of all maps $f\sm X$,
where $f$ is a stable trivial cofibration and $X$ is some symmetric
spectrum.  The following theorem is the monoid axiom for symmetric
spectra.  

\begin{theorem}\label{thm-monoid-axiom}
Transfinite compositions of pushouts of maps of $K$ are stable equivalences.
That is, suppose $\alpha $ is an ordinal and $A\mathcolon \alpha \to
\spec $ is a functor which preserves colimits and such that each map
$A_{i}\xrightarrow{}A_{i+1}$ is a pushout of a map of $K$.  Then the map
$A_{0}\xrightarrow{}\colim _{i<\alpha }A_{i}$ is a stable equivalence
\textup{(} and a level cofibration \textup{)}.
\end{theorem}

We then have the following two corollaries, which follow
from~\cite[Theorem 3.1]{schwede-shipley-monoids}.

\begin{corollary}
Suppose $R$ is a monoid in the category of symmetric spectra.  Then
there is a model structure on the category of $R$-modules where a map
$f\mathcolon X\xrightarrow{}Y$ is a weak equivalence
\textrm{(}fibration\textrm{)} if and only if the underlying map of
symmetric spectra is a stable equivalence \textrm{(}stable
fibration\textrm{)}.
Moreover, if $R$ is a commutative monoid then this is a monoidal model
category satisfying the monoid axiom.
\end{corollary}

\begin{corollary}
Suppose $R$ is a commutative monoid in the category of symmetric
spectra.  Then there is a model structure on the category of
$R$-algebras where a homomorphism $f\mathcolon X\xrightarrow{}Y$ is a
weak equivalence \textrm{(}fibration\textrm{)} if and only if the
underlying map of symmetric spectra is a stable equivalence
\textrm{(}stable fibration\textrm{)}.  Any cofibration of $R$-algebras
whose source is cofibrant as an $R$-module is a cofibration of
$R$-modules.  
\end{corollary}

Taking $R=S$ gives a model structure on the category of monoids of
symmetric spectra, the symmetric ring spectra.  

We will prove Theorem~\ref{thm-monoid-axiom} in a series of lemmas.  

\begin{lemma}\label{lem-K-stable-equiv}
Every map of $K$ is a level cofibration and a stable equivalence.  
\end{lemma}

\begin{proof}
This follows from parts three and five of Theorem~\ref{thm-boxprod}
applied to a $S$-cofibration and stable equivalence $f$ and the level
cofibration $*\to X$.
\end{proof}

This lemma means that every map of $K$ is a trivial cofibration in the
injective stable model structure.  Thus every colimit of pushouts of
maps of $K$ will also be a trivial cofibration in the injective stable
model structure.  This is the real reason that the monoid axiom is true
for $\spec $, but since we have not proved that the injective stable
structure really is a model structure, we give a direct proof of
Theorem~\ref{thm-monoid-axiom}.  

\begin{lemma}\label{lem-pushouts-of-level-cofibs}
Let 
\[
\begin{CD}
A @>f>> B \\
@VhVV @VVV \\
C @>>g> D
\end{CD}
\]
be a pushout square in $\spec $ where $f$ is a level cofibration and a
stable equivalence.  Then $g$ is a level cofibration and a stable
equivalence.  
\end{lemma}

\begin{proof}
The map $g$ is certainly a level cofibration.  There are factorizations
$f=pi$ and $h=qj$ where $i$ and $j$ are stable cofibrations and $p$ and
$q$ are level trivial fibrations.  Consider the resulting diagram of
pushout squares:
\[
\begin{CD}
A @>i>> A' @>p>> B \\
@VjVV @Vj'VV @VVj''V \\
C' @>i'>> Q' @>p'>> Q \\
@VqVV @Vq'VV @VVq''V \\
C @>>i''> D' @>>p''> D
\end{CD}
\]
Then $g=p''i''$.  We will show that $q$, $q''$, and $p'i'$ are stable
equivalences.  The two out of three property will then guarantee that
$g=p''i''$ is a stable equivalence.  The map $q$ is a level equivalence
by hypothesis.  The map $p'i'$ is a pushout of a level cofibration, so
is a level cofibration.  The properness of the injective model structure
then implies that $q''$ is a level equivalence.  Since $p'i'$ is a
pushout of the stable equivalence $f$ though the stable cofibration $j$,
it is a stable equivalence by the properness of the stable model
structure.  Thus $g=p''i''$ is a stable equivalence, as required. 
\end{proof}

The following lemma completes the proof of
Theorem~\ref{thm-monoid-axiom}. 

\begin{lemma}\label{lem-colimits-of-stable-equivalences}
Suppose $\alpha $ is an ordinal and $A\mathcolon \alpha
\xrightarrow{}\spec $ is a functor that preserves colimits such that the
map $A_{i}\xrightarrow{}A_{i+1}$ is a level cofibration and a stable
equivalence for all $i$.  Then the composite $A_{0}\xrightarrow{}\colim
A$ is a level cofibration and a stable equivalence.
\end{lemma}

\begin{proof}
Inductively define a new sequence $B_{i}$ by letting $B_{0}=A_{0}$ and
factoring the composite $B_{i}\xrightarrow{}A_{i}\xrightarrow{}A_{i+1}$
into a stable cofibration $B_{i}\xrightarrow{}B_{i+1}$ followed by a
level equivalence $B_{i+1}\xrightarrow{}A_{i+1}$.  For limit ordinals
$i$, define $B_{i}$ to be the colimit of the $B_{j}$ for $j<i$.  Then
each map $B_{i}\xrightarrow{}B_{i+1}$ is a stable trivial cofibration,
so the composite $A_{0}=B_0\xrightarrow{}\colim B$ is a stable trivial
cofibration.  On the other hand, the map $\colim B\xrightarrow{}\colim
A$ is a level equivalence, since homotopy of simplicial sets commutes
with transfinite compositions of monomorphisms.
\end{proof}

We now show that the model categories of modules and algebras are
homotopy invariant.  

\begin{lemma}\label{lem-level-quillen-equiv}
Suppose $R$ is a monoid in $\spec $, $f\mathcolon M\xrightarrow{}N$ is a
level equivalence of right $R$-modules, and $P$ is cofibrant in the
model category of left $R$-modules.  Then the map $f\sm_{R} P\mathcolon
M\sm _{R}P\xrightarrow{}N\sm _{R}P$ is a level equivalence.
\end{lemma}

\begin{proof}
The factorizations in the model category of $R$-modules are obtained by
applying the argument of Lemma~\ref{lem-factorization} to the sets $R\sm
FI_{\partial }$ and $R\sm FI_{\Lambda }$.  In particular, $P$ is a
retract, as an $R$-module, of a transfinite composition of pushouts of
coproducts of $R\sm FI_{\partial }$.  Thus we can assume that we have an
ordinal $\alpha $ and an $\alpha $-indexed diagram $X$ of $R$-modules
with colimit $P$, such that $X_{0}=0$ and each map $i_{\beta }\mathcolon
X_{\beta }\xrightarrow{}X_{\beta +1}$ is a pushout of a coproduct of
maps of $R\sm FI_{\partial }$.  We will show that the map $f\sm
_{R}X_{\beta }$ is a level equivalence for all $\beta \leq \alpha $ by
transfinite induction.

Getting the induction started is trivial.  For the successor ordinal
case of the induction, suppose $f\sm _{R}X_{\beta }$ is a level
equivalence.  We have a pushout square 
\[
\begin{CD}
M\sm _{R} A @>M\sm _{R}g>> M\sm _{R} B \\
@VVV @VVV \\
M\sm _{R}X_{\beta } @>>M\sm _{R}i_{\beta }> M\sm _{R}X_{\beta +1}
\end{CD}
\]
and an analogous pushout square with $N$ replacing $M$.  There is a map
from the pushout square with $M$ to the pushout square with $N$.  By
assumption, this map is a level equivalence on the lower right corner.
Since $M\sm _{R}(R\sm F_{n}K)=M\sm F_{n}K$ and similarly for $N$, this
map is also a level equivalence on the upper right and upper left
corners by Corollary~\ref{cor-X-smash} .  Furthermore the map $M\sm
_{R}g$ is a level cofibration by part three of
Theorem~\ref{thm-boxprod}.  Dan Kan's cubes lemma (see~\cite[Section
5.2]{hovey-model} or~\cite{kan-model}), applied one level at a time,
then implies that $f\sm _{R}X_{\beta +1}$ is a level equivalence.

Now suppose $\beta $ is a limit ordinal and $f\sm _{R}X_{\gamma }$ is a
level equivalence for all $\gamma <\beta $.  Since each map $M\sm
_{R}i_{\gamma }$ is a level cofibration, since it is a pushout of $M\sm
_{R}g$, the homotopy of $M\sm _{R}X_{\beta }$ is the colimit of the
homotopy of the $M\sm _{R}X_{\gamma }$.   A similar comment holds for
$N$ replacing $M$.  By comparing homotopy, we find that $f\sm
_{R}X_{\beta }$ is a level equivalence, as required. 
\end{proof}

\begin{proposition}\label{lem-stable-quillen-equiv}
Suppose $R$ is a monoid in $\spec $, $f\mathcolon M\xrightarrow{}N$ is a
stable equivalence of right $R$-modules, and $P$ is cofibrant in the
model category of left $R$-modules.  Then the map $f\sm_{R} P\mathcolon
M\sm _{R}P\xrightarrow{}N\sm _{R}P$ is a stable equivalence.
\end{proposition}

\begin{proof}
Factor $f$ in the model category of right $R$-modules into a stable
trivial cofibration followed by a level trivial fibration.  In view of
Lemma~\ref{lem-level-quillen-equiv}, we can then assume $f$ is a stable
trivial cofibration.  Any stable trivial cofibration is in $(J\sm R)\cof
$, so is a retract of a transfinite composition of pushouts of
coproducts of $J\sm R$.  Since retracts of stable equivalences are
stable equivalences, we may as well assume that $f$ is the transfinite
composition of an $\alpha $-indexed sequence $X_{\beta }$, where
$X_{0}=M$, $\colim _{\beta <\alpha }X_{\beta }=N$, and each map
$X_{\beta }\xrightarrow{}X_{\beta +1}$ is a pushout of a coproduct of
maps of $J\sm R$.  Note that $(J\sm R)\sm _{R}P=J\sm P$ consists of maps
which are level cofibrations (by part~3 of Theorem~\ref{thm-boxprod}
applied to the level cofibration $0\xrightarrow{}P$ and the stable
cofibrations of $J$) and stable equivalences (by part~5 of
Theorem~\ref{thm-boxprod}).  It follows from
Lemma~\ref{lem-pushouts-of-level-cofibs} that the map $X_{\beta }\sm
_{R}P\xrightarrow{}X_{\beta +1}\sm _{R}P$ is a level cofibration and a
stable equivalence.  Then
Lemma~\ref{lem-colimits-of-stable-equivalences} shows that $f\sm _{R}P$
is a level cofibration and stable equivalence, completing the proof. 
\end{proof}

Since $S$, the sphere spectrum, is cofibrant in $\spec$, Theorems~3.3
and~3.4 of~\cite{schwede-shipley-monoids}, apply to give the following
theorem.

\begin{theorem}\label{thm-monoids-and-homotopy}
Suppose $f\mathcolon R\xrightarrow{}R'$ is a stable equivalence of
monoids of symmetric spectra.  Then induction and restriction induce a
Quillen equivalence from the model category of $R$-modules to the model
category of $R'$-modules.  If, in addition, $R$ and $R'$ are commutative,
then induction and restriction induce a Quillen equivalence from the
model category of $R$-algebras to the model category of $R'$-algebras. 
\end{theorem}

\subsection{Semistable spectra}\label{sec-semistable} 

In this section we consider a subcategory of symmetric spectra called
the semistable spectra.  This subcategory sheds light on the difference
between stable equivalences and stable homotopy equivalences of
symmetric spectra.  As in Section~\ref{subsec-stable-category}, the
stable homotopy category is equivalent to the homotopy category of
semistable spectra obtained by inverting the stable homotopy
equivalences.  Semistable spectra also play a role
in~\cite{shipley-thh}.

Because stable equivalences are not always stable homotopy equivalences,
the stable homotopy groups are not algebraic invariants of stable
homotopy types.  So the stable homotopy groups of a spectrum may not be
``correct."  For any symmetric spectrum $X$, though, if there is a map
from $X$ to an $\Omega$-spectrum which induces an isomorphism in stable
homotopy, then the stable homotopy groups of $X$ are the ``correct"
stable homotopy groups.  In other words, these groups are isomorphic to
the stable homotopy groups of the stable fibrant replacement of $X$.
Spectra with this property are called semistable.

Let $L$ denote a fibrant replacement functor in $\spec$, obtained by
factoring $X \to *$ into a stable trivial cofibration followed by
a stable fibration, as in Section~\ref{subsec-equivalence}.

\begin{definition}\label{semistable}
A  \emph{semistable} symmetric spectrum is one for which the stable
fibrant replacement map, $X \xrightarrow{} LX$, is a stable homotopy
equivalence.
\end{definition}

Of course $X \xrightarrow{} LX$ is always a stable equivalence, but not
all spectra are semistable.  For instance, $F_1S^1$ is not semistable.
Certainly any stably fibrant spectrum, \ie an $\Omega $-spectrum, is
semistable.  In Section~\ref{subsec-stable-equivalence} we defined the
functor $R^{\infty}$ and noted that, although it is similar to the
standard $\Omega$-spectrum construction for (non-symmetric) spectra, it
is not always an $\Omega $-spectrum and $X \xrightarrow{} R^{\infty}X$ is
not always a stable homotopy equivalence, even if $X$ is level fibrant.
Let $K$ be a level fibrant replacement functor, the prolongation of a
fibrant replacement functor for simplicial sets.  The following
proposition shows that on semistable spectra $\RK$ does have these
expected properties.

\begin{proposition}\label{sf}
The following are equivalent.
\begin{enumerate}
\item The symmetric spectrum $X$ is semistable.
\item The map $X \xrightarrow{} RK X=\Hom(F_1S^1, KX)$ is a stable
homotopy equivalence.
\item $X \xrightarrow{} \RK X$ is a stable homotopy equivalence.
\item $\RK X$ is an $\Omega $-spectrum.
\end{enumerate}
\end{proposition}

Before proving this proposition we need the following lemma.

\begin{lemma}\label{lem-spec-Q}
Let $X\in \spec $.  Then $\pi _{k}(\RK X)_{n}$ and $\pi _{k+1}(\RK
X)_{n+1}$ are isomorphic groups, and $i_{\RK X}\colon (\RK X)_n \to
\Hom(F_1S^1, \RK X)_{n+1}$ induces a monomorphism $\pi _{k}(\RK
X)_{n}\xrightarrow{}\pi _{k+1}(\RK X)_{n+1} $.
\end{lemma}

\begin{proof}
As noted in the proof of Theorem~\ref{thm-U-detects-stable},
$\pi_{k}(\RK X)_{n}$ and $\pi_{k+1}(\RK X)_{n+1}$ are isomorphic to the
$(k-n)$th classical stable homotopy group of $X$.  However, the map
$\pi_{k}i_{\RK X}$ need not be an isomorphism.  Indeed, $\pi_{k}i_{\RK
X}$ is the map induced on the colimit by the vertical maps in the
diagram \[
\begin{CD}
\pi _{k}X_{n} @>>> \pi _{k+1}X_{n+1} @>>> \pi _{k+2}X_{n+2} @>>>
\cdots \\ @VVV @VVV @VVV \\ \pi _{k+1}X_{n+1} @>>> \pi
_{k+2}X_{n+2} @>>> \pi _{k+3} X_{n+3} @>>> \cdots
\end{CD}
\]
where the vertical maps are not the same as the horizontal maps, but
differ from them by isomorphisms.  The induced map on the colimit is
injective in such a situation, though not necessarily surjective.  For
an example of this phenomenon, note that we could have $\pi
_{k}X_{n}\cong \mathbb{Z}^{n-k}$, with the horizontal maps being the
usual inclusions, but the vertical maps begin the inclusions that take
$(x_{1},\dots ,x_{n-k})$ to $(0,x_{1},x_{2},\dots ,x_{n-k})$.  Then 
the element $(1,0,\dots ,0,\dots )$ of the colimit is not in the image
of the colimit of the vertical maps.  
\end{proof}

\begin{proof}[Proof of Proposition~\ref{sf}]
First we show that 1 implies 2 by using the following diagram.
$$
\begin{CD}
X @>>>  \Hom(F_1S^1, KX)\\
@VVV       @VVV\\
LX @>>>  \Hom(F_1S^1, KLX)
\end{CD}
$$
Since $\Hom(F_1S^1,K( -))$ preserves stable homotopy equivalences, both
vertical arrows are stable homotopy equivalences.  The bottom map is a
level equivalence since $LX$ is an $\Omega $-spectrum.  Hence the top map is
also a stable homotopy equivalence.

Also, 2 easily implies 3. Since $\Hom(F_1S^1,K(-))$ preserves stable
homotopy equivalences, $X \to \RK X$ is a colimit of stable homotopy
equivalences provided $X \to \Hom(F_1S^1,KX)$ is a stable homotopy
equivalence.

Next we show that 3 and 4 are equivalent.  The map $\pi_*X \to \pi_*\RK
X$ factors as $\pi_*X \to \pi_*(\RK X)_0 \to \pi_*\RK X$ where the first
map here is an isomorphism by definition.  Then by
Lemma~\ref{lem-spec-Q} we see that $\pi_*(\RK X)_0 \to \pi_* \RK X$ is
an isomorphism if and only if $\pi_*(\RK X)_n \to \pi_{*+1}(\RK
X)_{n+1}$ is an isomorphism for each $n$.

To see that 3 implies 1, consider the following diagram.
$$
\begin{CD}
X @>>> LX\\
@VVV       @VVV\\
\RK X @>>> \RK LX
\end{CD}
$$
By 3 and 4 the left arrow is a stable homotopy equivalence to an $\Omega
$ spectrum, $\RK X$.  Since $LX$ is an $\Omega $-spectrum the right arrow is a
level equivalence.  Since $X \to LX$ is a stable equivalence, the bottom
map must also be a stable equivalence.  But a stable equivalence between
$\Omega $-spectra is a level equivalence, so the bottom map is a level
equivalence.  Hence the top map is a stable homotopy equivalence.
\end{proof}

Two classes of semistable spectra are described in the following
proposition.  The second class includes the connective and convergent
spectra.

\begin{proposition}\label{sf-ex}
\begin{enumerate}
\item If the classical stable homotopy groups of $X$ are all finite then
$X$ is semistable.
\item Suppose that $X$ is a level fibrant symmetric spectrum and there
exists some $\alpha > 1$ such that $X_{n}\xrightarrow{}\Omega X_{n+1}$
induces an isomorphism $\pi_{k}X_{n}\xrightarrow{}\pi _{k+1}X_{n+1}$ for
all $k\leq \alpha n$ for sufficiently large $n$.  Then $X$ is
semistable.
\end{enumerate}
\end{proposition}

\begin{proof}
By Lemma~\ref{lem-spec-Q}, $\pi_k(\RK X)_n \to \pi_{k+1}(\RK X)_{n+1}$
is a monomorphism between two groups which are isomorphic.  In the first
case these groups are finite, so this map must be an isomorphism.  Hence
$\RK X$ is an $\Omega $-spectrum, so $X$ is semistable.

For the second part we also show that $\RK X$ is an $\Omega $-spectrum.
Since for fixed $k$ the maps $\pi_{k+i}X_{n+i} \to \pi_{k+1+i}X_{n+1+i}$
are isomorphisms for large $i$, $\pi_k(\RK X)_n \to \pi_{k+1}(\RK
X)_{n+1}$ is an isomorphism for each $k$ and $n$.
\end{proof}

The next proposition shows that stable equivalences between semistable
spectra are particularly easy to understand.

\begin{proposition}\label{sf-equivs}
Let $f\mathcolon X \xrightarrow{} Y$ be a map between two semistable
symmetric spectra.  Then $f$ is a stable equivalence if and only if it
is a stable homotopy equivalence.
\end{proposition}

\begin{proof}
Every stable homotopy equivalence is a stable equivalence, by
Theorem~\ref{thm-U-detects-stable}.  Conversely, if $f$ is a stable
equivalence, so is $Lf$.  Since stable equivalences between stably
fibrant objects are level equivalences by
Lemma~\ref{lem-stably-fibrant-equiv}, $Lf$ is in particular a stable
homotopy equivalence.  Since $X$ and $Y$ are semistable, both maps $X
\to LX$ and $Y \to LY$ are stable homotopy equivalences.  Hence $f$ is a
stable homotopy equivalence.
\end{proof}

\section{Topological spectra}\label{sec-topological-spectra}

In this final section, we describe symmetric spectra based on
topological spaces.  The advantage of symmetric spectra over ordinary
spectra is that symmetric spectra form a symmetric monoidal category,
whereas spectra do not.  Thus, we are interested in a setting where
symmetric spectra of topological spaces also form a symmetric monoidal
category.  This requires that we begin with a symmetric monoidal
category of pointed topological spaces.  The category $\top $ of all
pointed topological spaces is not symmetric monoidal under the smash
product; the associativity isomorphism is not always continuous.  So we
begin in Section~\ref{subsec-compactly-generated} by describing
compactly generated spaces, a convenient category of topological spaces
which is closed symmetric monoidal.  Then in
Section~\ref{subsec-top-spectra} we define topological spectra, and in
Section~\ref{subsec-top-stable} we discuss the stable model structure on
topological symmetric spectra.  In Section~\ref{subsec-top-prop} we show
that the stable model structure is monoidal and proper.  We also
consider monoids and modules over them, though here our results are much
less strong than in the simplicial case.  Indeed, in this regard, we do
not know more than we do in a general monoidal model category, where we
can appeal to the results of~\cite{hovey-monoids}.

\subsection{Compactly generated spaces}\label{subsec-compactly-generated}

In this section, we describe a closed symmetric monoidal category of
topological spaces.  There are many different choices for such a
category~\cite{wyler}, and any one of them suffices to construct
topological symmetric spectra.  However, in order to get model
categories of modules over a monoid, one needs especially good
properties of the underlying category of topological spaces.  For
this reason, we work in the category of compactly generated spaces, for
which the basic reference is the appendix of~\cite{lewis-thesis}.  This
is the same category used in~\cite{lewis-may-steinberger}
and~\cite{elmendorf-kriz-mandell-may}.  We also discuss the model
structure on the category of compactly generated spaces.

\begin{definition}\label{defn-compactly-generated}
Suppose $X$ is a topological space.  A subset $A$ of $X$ is called
\emph{compactly open} if, for every continuous map $f\mathcolon
K\xrightarrow{}X$ where $K$ is compact Hausdorff, $f^{-1}(A)$ is open.
We denote by $kX$ the set $X$ with the topology consisting of the
compactly open sets, and we say that $X$ is a \emph{$k$-space} if
$kX=X$.  We denote the full subcategory of $\top $ consisting of pointed
$k$-spaces by $\cg $.  The space $X$ is \emph{weak Hausdorff} if for
every continuous map $f\mathcolon K\xrightarrow{}X$ where $K$ is compact
Hausdorff, the image $f(K)$ is closed in $X$.  We denote the full
subcategory of $\cg $ consisting of weak Hausdorff spaces by $\spaces $, and
refer to its objects as \emph{compactly generated} pointed spaces.
\end{definition}

Note that every locally compact Hausdorff space is compactly generated.
The functor $k\mathcolon \top \xrightarrow{}\cg $ is a right adjoint to
the inclusion functor.  It follows that $\cg $ is bicomplete.  Colimits
are taken in $\top $, so in particular a colimit of $k$-spaces is again
a $k$-space.  One must apply $k$ to construct limits.  Thus, the product
in $\cg $ is the topological space $k(X\times Y)$.  If $X$ is locally
compact Hausdorff, then $k(X\times Y)=X\times Y$.

One can easily see that limits of weak Hausdorff spaces are again weak
Hausdorff.  Thus, $\spaces $ is complete and the inclusion functor
$\spaces \xrightarrow{}\cg $ preserves all limits.  Freyd's adjoint
functor theorem then guarantees the existence of a left adjoint
$w\mathcolon \cg \xrightarrow{}\spaces $ to the inclusion functor.  The
space $wX$ is the maximal weak Hausdorff quotient of $X$, but it is very
hard to describe $wX$ in general.  The functor $w$ creates colimits in
$\spaces $, so that $\spaces $ is cocomplete.  However, very often the
colimit of a diagram in $\spaces $ is the same as the colimit of the
diagram in $\cg $, so that we never need to apply $w$ in practice.  
This is the case for directed systems of injective maps~\cite[Proposition
A.9.3]{lewis-thesis}, and for pushouts of closed
inclusions~\cite[Proposition A.7.5]{lewis-thesis}.  In particular,
cofibrations are both injective and closed inclusions~\cite[Proposition
A.8.2]{lewis-thesis}, so these remarks will apply to transfinite
compositions and pushouts of cofibrations.  

Both the category $\mathcal{K}$ of unpointed $k$-spaces and the category
$\mathcal{T}$ of unpointed compactly generated spaces are closed
symmetric monoidal under the $k$-space version of the
product~\cite[Theorem A.5.5]{lewis-thesis}.  The right adjoint of
$X\times -$ is the functor $kC(X,-)$.  Here $C(X,Y)$ is the set of
continuous maps from $X$ to $Y$.  A subbasis for the topology on
$C(X,Y)$ is given by the sets $S(f,U)$, where $f\mathcolon
K\xrightarrow{}X$ is a continuous map, $K$ is compact Hausdorff, and $U$
is open in $Y$.  A map $g\mathcolon X\xrightarrow{}Y$ is in $S(f,U)$ if
and only if $(g\circ f)(K)\subseteq U$.  If $X$ is locally compact
Hausdorff, then this topology on $C(X,Y)$ is the same as the
compact-open topology.  However, there is no guarantee that $C(X,Y)$ is
a $k$-space, even if $X$ is locally compact Hausdorff, and so we must
apply $k$.  If $X$ and $Y$ are both second countable Hausdorff, and $X$
is locally compact, then $kC(X,Y)=C(X,Y)$~\cite[p. 161]{lewis-thesis}.
Fortunately, $C(X,Y)$ is weak Hausdorff as long as $Y$ is~\cite[Lemma
A.5.2]{lewis-thesis}, so $kC(X,Y)$ is already weak Hausdorff and there
is no need to apply $w$ to get the closed structure in $\mathcal{T}$.

The product induces a smash product $X\sm Y=k(X\times Y)/(X\vee Y)$ on
$\cg $.  Then $\cg $ is a closed symmetric monoidal category; the right
adjoint of $X\sm -$ is $kC_{*}(X,-)$, $k$ applied to the subspace of
$kC(X,Y)$ consisting of pointed continuous maps.  The subcategory
$\spaces $ of $\cg $ is closed under the smash product~\cite[Lemma
A.6.2]{lewis-thesis}, so $\spaces $ is also a closed symmetric monoidal
category.  In fact, the functor $w\mathcolon \cg \xrightarrow{}\spaces $
is symmetric monoidal.  This fact does not appear
in~\cite{lewis-thesis}; its proof uses adjointness and the fact that
$kC_{*}(X,Y)$ is already weak Hausdorff if $Y$ is so.

We now discuss the geometric realization.  One source for the
geometric realization is~\cite[III]{gabriel-zisman}.  Note that the
geometric realization of a finite simplicial set is compact Hausdorff,
and so is a $k$-space.  Since the geometric realization is a left
adjoint, it commutes with colimits, and so the geometric realization of
any simplicial set is a $k$-space.  Also, the geometric realization of
any simplicial set is a cell complex, so is Hausdorff.  We can therefore
consider the geometric realization as a functor to $\spaces $.  We
denote the geometric realization by $\georeal{-}\mathcolon \sset
\xrightarrow{}\spaces $, and its right adjoint, the singular functor, by
$\Sing \mathcolon \spaces \xrightarrow{}\sset $.

The geometric realization has a number of properties not expected in a
general left adjoint.  For example, it preserves finite
limits~\cite[III.3]{gabriel-zisman}.  Gabriel and Zisman prove this
using a different category of topological spaces, but their proof relies
on first working with certain finite simplicial sets and their geometric
realizations, and then using the closed symmetric monoidal structure to
extend to general simplicial sets.  Thus their argument works in any
closed symmetric monoidal category of topological spaces that contains
finite CW complexes.  It follows easily from this that the geometric
realization preserves the smash product, and so defines a symmetric
monoidal functor $\georeal{-}\mathcolon \sset \xrightarrow{}\spaces $.  

The categories $\spaces $, $\cg $ and $\top $ are all model categories,
and both $w\mathcolon \spaces \xrightarrow{}\cg $ and the inclusion $\cg
\xrightarrow{}\top $ are Quillen equivalences.  The easiest way to
describe the model structure is to define $f$ to be a weak equivalence
(fibration) if and only $\Sing f$ is a weak equivalence (fibration) of
simplicial sets.  Then $f$ is a cofibration if and only if it has the
\llp all trivial fibrations.  Of course, the weak equivalences are the
maps which induce isomorphisms on homotopy groups at all possible
basepoints.  All inclusions of relative CW-complexes are cofibrations,
and any cofibration is a retract of a transfinite composition of
inclusions of relative CW-complexes.  In fact, the cofibrations are the
class $\georeal{I_{\partial }}\cof $ and the trivial cofibrations are
the class $\georeal{I_{\Lambda }}\cof $.  Every topological space is
fibrant, and the geometric realization is a Quillen equivalence.  Both
the geometric realization and the singular functor preserve and reflect
weak equivalences, and the maps $\georeal{\Sing X}\xrightarrow{}X$ and
$K\xrightarrow{}\Sing \georeal{K}$ are weak equivalences.

In proving the claims of the last paragraph, one needs the small object
argument in each of the categories involved.  The most obvious smallness
statement is the following which states that all spaces are \emph{small}
relative to inclusions.  

\begin{lemma}\label{lem-spaces-are-small}
Suppose $K$ is a topological space with cardinality $\gamma $, and
suppose $X\mathcolon \alpha \xrightarrow{}\top $ is a colimit-preserving
functor, where $\alpha $ is a $\gamma $-filtered ordinal, and suppose
that each map $X_{\beta }\xrightarrow{}X_{\beta +1}$ is an inclusion.
Then the map 
\[
\colim_{\alpha}\top (K,X)\xrightarrow{} \top
(K,\colim_{\alpha } X)
\]
is an isomorphism.
\end{lemma}

\begin{proof}
One can verify using transfinite induction that each map $X_{\beta
}\xrightarrow{}\colim _{\alpha}X$ is an inclusion.  Given
a map $f\mathcolon K\xrightarrow{}\colim _{\alpha}X$, there is
a $\beta <\alpha $ such that $f$ factors as a map of sets through
$X_{\beta }$, since $\alpha $ is $\gamma $-filtered.  But then $f$ is
automatically continuous as a map to $X_{\beta }$ since the map
$X_{\beta }\xrightarrow{}\colim _{\alpha}X$ is an inclusion.
\end{proof}

For Lemma~\ref{lem-top-factorization} we need the class of inclusions to
be preserved under various constructions.  Inclusions are not
well behaved in general.  However, in $\spaces $, \emph{closed}
inclusions are well behaved: they are closed under
coproducts~\cite[Proposition A.7.1(c)]{lewis-thesis},
pushouts~\cite[Proposition A.7.5]{lewis-thesis}, and transfinite
compositions~\cite[Proposition A.9.3]{lewis-thesis}.  Furthermore, the
smash product preserves closed inclusions~\cite[Proposition
A.7.3]{lewis-thesis}.  Given this, the proof of the Factorization
Lemma~\ref{lem-factorization} applies to prove the following analogue.

\begin{lemma}\label{lem-top-factorization}
Let $I$ be a set of closed inclusions in the category $\spaces $.  There
is a functorial factorization of every map in $\spaces $ as an
$I$-cofibration followed by an $I$-injective map.  
\end{lemma}

The model structure on $\spaces $ is monoidal and proper.  The fact that
it is monoidal follows immediately from the corresponding fact for
simplicial sets, since the geometric realization is monoidal and the
generating cofibrations and trivial cofibrations are in the image of the
geometric realization.  Properness for $\spaces $ is a little harder:
right properness is formal, since all objects are fibrant, but left
properness requires some care with local coefficients~\cite[Proposition
11.1.11]{hirschhorn}.  The monoid axiom holds in both $\spaces$ and $\cg
$.  We get a model category of monoids in $\spaces $, but not in $\cg $.
This is the reason we choose to use $\spaces$ here and is explained 
in~\cite{hovey-monoids}.  

\subsection{Topological spectra}\label{subsec-top-spectra}

In this section we define the category of topological symmetric
spectra.  The basic definitions are completely analogous to the
definitions for simplicial symmetric spectra.  

The basic definitions for topological spectra and topological symmetric
spectra are straightforward.  Let $S^{1}\in \spaces $ be the unit circle
in $\mathbf{R}^{2}$ with basepoint $(1,0)$.

\begin{definition}\label{def-top-spectrum}
A \emph{topological spectrum} is 
\begin{enumerate}
\item a sequence $X_0,X_1,\dots,X_n,\dots$ in $\spaces $; and  
\item a pointed continuous map $\sigma\mathcolon S^1\sm X_n\to X_{n+1}$
for each $n\ge0$.
\end{enumerate}
The maps $\sigma$ are the \emph{structure maps} of the spectrum.  A
\emph{map of topological spectra} $f\mathcolon X\to Y$ is a sequence of
pointed continuous maps $f_n\mathcolon X_n\to Y_n$ such that the diagram
\[
\begin{CD}
S^1\sm X_n @>\sigma>> X_{n+1}\\ 
@VS^1\sm f_nVV @Vf_{n+1}VV \\ 
S^1\sm Y_n @>\sigma>> Y_{n+1}
\end{CD}
\]
is commutative for each $n\ge0$. Let $\topBF$ denote the category of
topological spectra.
\end{definition}

Let $S^{p}=(S^{1})^{\sm p}$ with the (continuous) left permutation
action of $\Sigma _{p}$.

\begin{definition}\label{defn-top-symm-spectrum}
A \emph{topological symmetric spectrum} is
\begin{enumerate}
\item a sequence $X_{0},X_{1},\dots ,X_{n},\dots $ in $\spaces $;
\item a pointed continuous map $\sigma\mathcolon S^1\sm X_n\to X_{n+1}$
for each $n\ge0$; and
\item a basepoint preserving continuous left action of $\Sigma_n$ on
$X_n$ such that the composition
\[
\sigma^p=\sigma\circ (S^1\sm\sigma)\circ\dots\circ (S^{p-1}\sm
\sigma)\mathcolon S^p\sm X_n\to X_{n+p},
\] 
of the maps $S^i\sm S^1\sm X_{n+p-i-1}\xrightarrow{S^i\sm\sigma} S^{
i}\sm X_{n+p-i}$ is $\Sigma_p\times\Sigma_n$-equivariant for $p\ge1$
and $n\ge0$.
\end{enumerate}
A \emph{map of topological symmetric spectra} $f\mathcolon X\to Y$ is a
sequence of pointed continuous maps $f_n\mathcolon X_n\to Y_n$ such that
$f_n$ is $\Sigma_n$-equivariant and the diagram
\[
\begin{CD}
S^1\sm X_n @>\sigma>> X_{n+1}\\ 
@VS^1\sm f_nVV @Vf_{n+1}VV \\ 
S^1\sm Y_n @>\sigma>> Y_{n+1}
\end{CD}
\] 
is commutative for each $n\ge0$. Let $\topspec$ denote the category of
symmetric spectra.
\end{definition}


\begin{example}\label{bordism}
The {\em real bordism spectrum} $MO$ is the sequence of spaces $MO_n=
S^n \sm_{O_n} (EO_n)_+$.  Here we require that $EG$ be a functorial
construction of a free contractible $G$-space, for example Segal's
construction,~\cite{segal}.  Since there are commuting left and right
actions of $O_{n}$ on $S^{n}$, the symmetric group acts on the left, 
via inclusion into $O_n$,  on
$S^{n}$ and on $MO_{n}$.  One can check that the usual structure maps
$S^{1}\sm MO_n\xrightarrow{}MO_{n+1}$ make $MO$ into a topological
symmetric spectrum.  Similarly, the other bordism spectra such as $MU$
are naturally topological symmetric spectra.  One can check that these
examples are even commutative symmetric ring spectra.
\end{example}

With these definitions, most of the basic results of
Sections~\ref{sec-spectra} and~\ref{sec-smash-product} go through
without difficulty, replacing simplicial sets by topological spaces.
That is, the category of topological symmetric spectra is the category
of $S$-modules, where $S$ is the topological symmetric sequence
$(S^{0},S^{1},S^{2},\dots )$, and $S^{n}=(S^{1})^{\sm n}$ with the
permutation action of $\Sigma _{n}$.  The category of topological
symmetric spectra is a closed symmetric monoidal category.  

For topological symmetric spectra $X$ and $Y$, we can form the mapping
space $C_{\topspec}(X,Y)$ whose points are maps of symmetric spectra
from $X$ to $Y$, topologized as $k$ applied to the subspace of $k\prod
kC_{*}(X_{n},Y_{n})$.  We get the usual adjunction properties.
Similarly, we can form function objects $\ihom _{\Sigma }(X,Y)$ and
$\ihom _{S}(X,Y)$.  Indeed, given topological symmetric sequences $X$
and $Y$, the $n$th space of the symmetric sequence $\ihom _{\Sigma
}(X,Y)$ is $C_{\topsymseq}(X\otimes G_{n}S^{0},Y)$.  A point is a map of
symmetric sequences that raises degree by $n$; the leftover $\Sigma
_{n}$ action defines a symmetric sequence.  Similarly, given topological
symmetric spectra $X$ and $Y$, the $n$th space of the symmetric spectrum
$\ihom _{S}(X,Y)$ is $C_{\topspec }(X\sm F_{n}S^{0},Y)$.  A point is a
map of symmetric spectra that raises degree by $n$, and the topology is
$k$ applied to the subspace topology inherited from $C_{\topsymseq }(X\otimes
G_{n}S^{0},Y)$.  The structure map $S^{1}\sm C_{\topspec}(X\sm
F_{n}S^{0},Y)\xrightarrow{}C_{\topspec }(X\sm F_{n+1}S^{0},Y)$ takes
$(t,f)$ to $\sigma _{t}(f)$, where $t\in S^{1}$ and $f$ is a map of
degree $n$, and $\sigma _{t}(f)(x)=f(\sigma _{X}(t,x))=\sigma
_{Y}(t,f(x))$.

\subsection{Stable model structure}\label{subsec-top-stable}

In this section we discuss the stable homotopy theory of topological
symmetric spectra.  Rather than proceeding analogously to
Section~\ref{sec-homotopy}, we use the geometric realization and
singular functor to lift the stable model structure on simplicial
symmetric spectra to topological symmetric spectra.  

Recall that both the geometric realization and the singular functor
reflect and preserve weak equivalences, and that maps 
$A\xrightarrow{}\Sing \georeal{A}$ and the counit $\georeal{\Sing
X}\xrightarrow{}X$ of the adjunction are both natural weak equivalences.

The following proposition is then immediate.  

\begin{proposition}\label{top-geo-real}
The geometric realization induces a functor $\georeal{-}\mathcolon \spec
\xrightarrow{}\topspec $ such that $\georeal{X}_{n}=\georeal{X_n}$.  The
geometric realization is symmetric monoidal.  The singular functor
induces a functor $\Sing \mathcolon \topspec \xrightarrow{}\spec $,
right adjoint to $\georeal{-}$, such that $(\Sing X)_{n}=\Sing X_{n}$.
Both the geometric realization and the singular complex detect and
preserve level equivalences, and the unit and counit maps of the
adjunction are level equivalences.
\end{proposition}

We then make the following definition.  

\begin{definition}\label{defn-top-model}
Suppose $f\mathcolon X\xrightarrow{}Y$ is a map of topological symmetric
spectra.  
\begin{enumerate}
\item Define $f$ to be a \emph{stable equivalence} if $\Sing f$ is a
stable equivalence in $\spec $. 
\item Define $f$ to be a \emph{stable fibration} if $\Sing f$ is a
stable fibration in $\spec $. 
\item Define $f$ to be a \emph{stable cofibration} if $f$ has the \llp
all stable trivial fibrations (\ie maps which are both stable
equivalences and stable fibrations).  
\end{enumerate}
\end{definition}

We will prove that these structures form a model structure on $\topspec
$, called the \emph{stable model structure}.  Here are some of their
basic properties.  

\begin{lemma}\label{lem-top-basic-stable}
Suppose $f\mathcolon X\xrightarrow{}Y$ is a map of topological symmetric
spectra.  
\begin{enumerate}
\item If $f$ is a level equivalence, then $f$ is a stable equivalence.  
\item The map $f$ is a stable trivial fibration if and only if $f$ is a level
trivial fibration.  
\item The map $f$ is a stable fibration if and only if $f$ is a level fibration
and the diagram 
\[
\begin{CD}
X_{n} @>>> X_{n+1}^{S^{1}} \\
@VVV @VVV \\
Y_{n} @>>> Y_{n+1}^{S^{1}}
\end{CD}
\]
is a homotopy pullback square. In particular, a topological symmetric
spectrum $X$ is stably fibrant if and only if the map
$X_{n}\xrightarrow{}X_{n+1}^{S^{1}}$ is a weak equivalence for all
$n\geq 0$.  Such $X$ are called \emph{$\Omega $-spectra}.  
\item The map $f$ is a stable equivalence between $\Omega $-spectra if and only
if $f$ is a level equivalence between $\Omega $-spectra.  
\item The map $f$ is a stable fibration between $\Omega $-spectra if and
only if $f$ is a level fibration between $\Omega $-spectra.
\item The map $f$ is a stable cofibration if and only if $f\in
\georeal{FI_{\partial }}\cof $.  
\item If $f$ is a stable cofibration, then $f$ is a level cofibration.  
\end{enumerate}
\end{lemma}

\begin{proof}
Part 1 follows from the corresponding statement for simplicial symmetric
spectra.  Part~2 follows in the same way.  For part~3, note that
$\georeal{-}$ commutes with smashing with $S^{1}$, so, by adjointness,
$\Sing (X_{n+1}^{S^{1}})$ is isomorphic to $(\Sing X_{n+1})^{S^{1}}$.
Since $\Sing $ commutes with pullbacks and preserves and reflects weak
equivalences, the square in question is a homotopy pullback square if
and only if $\Sing f$ is a stable fibration.  Parts~4 and~5 follow
immediately from the corresponding statements for simplicial symmetric
spectra.  For part~6, note that $g$
is a stable trivial fibration if and only if $\Sing g$ is in
$FI_{\partial }\inj $, which, by adjointness, holds if and only if $g$
is in $\georeal{FI_{\partial }}\inj $.  The result follows.  For part~7,
note that level cofibrations are the class of $R(tK)$-projective maps,
just as in part~\ref{RtK} of Lemma~\ref{lem-level-properties}, where
$tK$ is the class of trivial fibrations in $\spaces $, and
$R(tK)=\bigcup_{n} R_{n}(tK)$, and $R_{n}\mathcolon \spaces
\xrightarrow{}\topspec $ is the right adjoint to the evaluation functor
$\Ev _{n}$.  Thus, it suffices to show that each map in
$\georeal{FI_{\partial }}$ is a level cofibration, but this is clear.  
\end{proof}
 
To characterize the stable cofibrations, we can introduce the latching
space $L_{n}X$ and the natural transformation $i\mathcolon
L_{n}X\xrightarrow{}X_{n}$ just as we did in
Section~\ref{subsec-stable-cofib}.  In the category of $\Sigma
_{n}$-spaces, a \emph{$\Sigma _{n}$-cofibration} is an equivariant map
$f\mathcolon X\xrightarrow{}Y$ which has the \llp all equivariant maps
that are both weak equivalences and fibrations.  In particular, if
$g\mathcolon X\xrightarrow{}Y$ is an equivariant relative CW complex on
which $\Sigma _{n}$ acts cellularly and also acts freely on the cells
not in $X$, then $f$ is a $\Sigma _{n}$-cofibration.

\begin{lemma}\label{lem-top-stable-cofib}
A map $f\mathcolon X\xrightarrow{}Y$ of topological symmetric spectra
is a stable cofibration if and only if for all $n\geq 0$ the induced map
$\Ev_{n}(f\boxprod i)\mathcolon X_{n}\amalg _{L_{n}X}
L_{n}Y\xrightarrow{}Y_{n}$ is a $\Sigma _{n}$-cofibration.  
\end{lemma}

The proof of this lemma is exactly the same as the proof of
Proposition~\ref{prop-stable-cofib}. 

In order to prove that these classes define a model structure on
$\topspec $, the first thing we need is a version of the Factorization
Lemma~\ref{lem-factorization}.  The following proposition is the
analogue of Proposition~\ref{prop-spec-smallness}, and its proof is
similar, given Lemma~\ref{lem-spaces-are-small}.  Define the
\emph{cardinality} of a topological symmetric spectrum $X$ to be the
cardinality of the disjoint union of its spaces $X_{n}$.  

\begin{proposition}\label{prop-topspec-smallness}
Let $X$ be a topological symmetric spectrum of cardinality $\gamma $.
Let $\alpha $ be a $\gamma $-filtered ordinal and let $D\mathcolon
\alpha \xrightarrow{} \topspec $ be colimit-preserving $\alpha $-indexed
diagram of level inclusions.  Then the natural map
\[
\colim _{\alpha }\topspec (X,D)\xrightarrow{}\topspec (X,\colim _{\alpha
}D) 
\]
is an isomorphism.  
\end{proposition}

This proposition immediately leads to a factorization lemma for
topological symmetric spectra, analogous to the Factorization
Lemmas~\ref{lem-factorization} and~\ref{lem-top-factorization} for
simplicial symmetric spectra and for compactly generated spaces.  

\begin{lemma}\label{lem-topspec-factorization}
Let $I$ be a set of level closed inclusions in $\topspec $.  There is a
functorial factorization of every map of topological symmetric spectra
as an $I$-cofibration followed by an $I$-injective.
\end{lemma}

We can now establish the projective and stable model structures on
$\topspec $.  We begin with the projective model structure.  

\begin{theorem}\label{thm-top-projective-model}
The stable cofibrations, level fibrations, and level equivalences define
a model structure on $\topspec $, called the \emph{projective level
structure}.  The geometric realization $\georeal{-}\mathcolon \spec
\xrightarrow{}\topspec $ defines a Quillen equivalence between the
projective level structure on $\spec $ and the projective level
structure on $\topspec $.  
\end{theorem}

\begin{proof}
Certainly $\topspec $ is bicomplete, and the retract and two out of
three axioms are immediate consequences of the definitions.  The lifting
axiom for stable cofibrations and level trivial fibrations follows from
the definitions as well.  Adjointness implies that the level trivial
fibrations are the class $\georeal{FI_{\partial }}\inj $.  Part two of
Lemma~\ref{lem-top-basic-stable} then implies that the stable
cofibrations are the class $\georeal{FI_{\partial }}\cof $.
Lemma~\ref{lem-topspec-factorization} produces the required
factorization into a stable cofibration followed by a level trivial
fibration.  

Adjointness also implies that the level fibrations are the class
$\georeal{FI_{\Lambda }}\inj $.  Each map of $\georeal{FI_{\Lambda }}$
is a stable cofibration and a level equivalence, and in particular a
level trivial cofibration.  Just as in part~\ref{RK} of
Lemma~\ref{lem-level-properties}, the level trivial cofibrations are the
$RK$-projective maps.  It follows that every map of
$\georeal{FI_{\Lambda }}\cof $  is a level trivial cofibration.  
We knew already that each map of $\georeal{FI_{\Lambda }}\cof $ is a
stable cofibration.  Thus, Lemma~\ref{lem-topspec-factorization} gives
us a factorization of any map into a stable cofibration and level
equivalence followed by a level fibration.  

We are left with the lifting axiom for stable cofibrations and level
equivalences, and level fibrations.  Let $f$ be a stable cofibration and
level equivalence.  Write $f=pg$, where $g\in \georeal{FI_{\Lambda
}}\cof $ and $p$ is a level fibration.  Then $p$ is a level equivalence,
so $f$ has the \llp $p$.  It follows from Proposition~\ref{prop-retract}
that $f$ is a retract of $g$, so $f$ has the \llp level fibrations.  

The singular functor is obviously a right Quillen functor, and reflects
all level equivalences.  Now let $K$ be a fibrant replacement functor on
$\topspec $, obtained, for example, from the Factorization
Lemma~\ref{lem-topspec-factorization} applied to $\georeal{FI_{\Lambda
}}$.  Then in order to show that the pair $(\georeal{-},\Sing )$ is a
Quillen equivalence, it suffices by Lemma~\ref{lem-Quillen-equivalences}
to show that the map $X\xrightarrow{}\Sing K\georeal{X}$ is a level
equivalence for all simplicial symmetric spectra $X$.  However, the map
$X\xrightarrow{}\Sing \georeal{X}$ is a level equivalence.  Since the
map $\georeal{X}\xrightarrow{}K\georeal{X}$ is a level equivalence, so
is the map $\Sing \georeal{X}\xrightarrow{}\Sing K\georeal{X}$.  By the
two out of three axiom, the composite $X\xrightarrow{}\Sing
K\georeal{X}$ is a stable equivalence, as required.
\end{proof}

We do not know if the injective level structure on $\topspec $ is a
model structure.  

\begin{theorem}\label{thm-top-model}
The stable model structure on $\topspec $ is a model structure.
We have the following properties\textup{:} 
\begin{enumerate}
\item A map $f\in \spec $ is a stable equivalence if and only if
$\georeal{f}$ is a stable equivalence in $\topspec $;
\item A map $g\in \topspec $ is a stable equivalence if and only if
$\Sing g$ is a stable equivalence in $\spec $; and
\item The geometric realization and singular functor define a Quillen
equivalence between the stable model structure on simplicial symmetric
spectra and the stable model structure on topological symmetric
spectra.  
\end{enumerate}
\end{theorem}

\begin{proof}
The category $\topspec $ is bicomplete, and the two out of three and
retract axioms are immediate consequences of the definitions.  The
stable trivial fibrations and the level trivial fibrations coincide, by
Lemma~\ref{lem-top-basic-stable}.  Thus the lifting and factorization
axioms for stable cofibrations follow from the corresponding axioms in
the projective level structure.  

The stable fibrations form the class $\georeal{J}\inj $,
where $J$ is the set of generating stable trivial cofibrations used in
Section~\ref{subsec-stable-model-cat}.  In
Proposition~\ref{prop-top-J}, we show that every map of $\georeal{J}\cof
$ is a stable cofibration and a stable equivalence.  Hence the
Factorization Lemma~\ref{lem-topspec-factorization} gives us the required
(functorial) factorization of an arbitrary map into a stable trivial
cofibration followed by a stable fibration.  

Only one lifting property remains.  Suppose $i$ is a stable trivial
cofibration.  We can factor $i$ as $i=qi'$, where $q$ is a stable
fibration and $i'\in \georeal{J}\cof $.  By the preceding paragraph,
$i'$ is a stable equivalence.  By the two out of three axiom, $q$ is a
stable equivalence.  Thus $i$ has the \llp $q$, and so the Retract
Argument~\ref{prop-retract} implies that $i$ is a retract of $i'$.  Thus
$i\in \georeal{J}\cof $, and so $i$ has the \llp every stable fibration
$p$.  This completes the proof that $\topspec $ is a model category.

Now, by definition, $\Sing \mathcolon \topspec \xrightarrow{}\spec $ is
a right Quillen functor, and $f\in \topspec $ is a stable equivalence if
and only if $\Sing f$ is a stable equivalence in $\spec $.  Now let $K$
be a fibrant replacement functor on $\topspec $, obtained, for example,
from the Factorization Lemma~\ref{lem-topspec-factorization} applied to
$\georeal{J}$.  Then in order to show that the pair $(\georeal{-},\Sing
)$ is a Quillen equivalence, it suffices by
Lemma~\ref{lem-Quillen-equivalences} to show that the map
$X\xrightarrow{}\Sing K\georeal{X}$ is a stable equivalence for all
simplicial symmetric spectra $X$.  However, the map
$X\xrightarrow{}\Sing \georeal{X}$ is a level equivalence.  Since the
map $\georeal{X}\xrightarrow{}K\georeal{X}$ is a stable equivalence, so
is the map $\Sing \georeal{X}\xrightarrow{}\Sing K\georeal{X}$.  By the
two out of three axiom, the composite $X\xrightarrow{}\Sing
K\georeal{X}$ is a stable equivalence, as required.  

Finally, we show that $\georeal{-}$ detects and preserves stable
equivalences.  Given a map $f\mathcolon X\xrightarrow{}Y$, consider the
commutative square below.
\[
\begin{CD}
X @>f>> Y \\
@VVV @VVV \\
\Sing \georeal{X} @>>\Sing \georeal{f}> \Sing \georeal{Y}
\end{CD}
\]
The vertical maps are level equivalences.  Now $\Sing \georeal{f}$ is a
stable equivalence if and only if $\georeal{f}$ is a stable
equivalence.  Hence $f$ is a stable equivalence if and only if
$\georeal{f}$ is a stable equivalence.  
\end{proof}

We still owe the reader a proof that the maps of $\georeal{J}\cof $ are
stable cofibrations and stable equivalences.  We begin with the
following lemma. 

\begin{lemma}\label{lem-top-J-pushout}
Suppose $f\mathcolon A\xrightarrow{}B$ is an element of $J$, and suppose
the following square is a pushout square in $\topspec $. 
\[
\begin{CD}
\georeal{A} @>\georeal{f}>> \georeal{B} \\
@VrVV @VVsV \\
C @>>g> D
\end{CD}
\]
Then $g$ is a stable equivalence and a stable cofibration.  
\end{lemma}

\begin{proof}
It is clear that $\georeal{f}$, and hence $g$, is a stable cofibration.
Note that $\georeal{f}$ is itself a stable equivalence, since $\Sing
\georeal{f}$ is level equivalent to the stable equivalence $f$, so is a
stable equivalence.  

We first suppose that $r$ is a level cofibration.  Consider the pushout
square in $\spec $ below.  
\[
\begin{CD}
\Sing \georeal{A} @>\Sing \georeal{f}>> \Sing \georeal{B} \\
@V\Sing rVV @VVs'V \\
\Sing C @>>g'> E
\end{CD}
\]
Since $\Sing $ preserves level cofibrations (since it preserves
injections), each map in this square is a level cofibration.
Furthermore, $\Sing \georeal{f}$ is a stable equivalence.  By the
properness result of Theorem~\ref{thm-proper}, we find that $g'$ is a
stable equivalence.  There is a map $h\mathcolon E\xrightarrow{}\Sing D$
such that $hg'=\Sing g$, so it will suffice to prove that $h$ is a level
equivalence.  

To see this, we apply the geometric realization again to obtain a
pushout square
\[
\begin{CD}
\georeal{\Sing \georeal{A}} @>\georeal{\Sing \georeal{f}}>>
\georeal{\Sing \georeal{B}} \\
@V\georeal{\Sing r}VV @VV\georeal{s'}V \\
\georeal{\Sing C} @>>\georeal{g'}> \georeal{E}
\end{CD}
\]
Again, each of the maps in this diagram is a level cofibration.  There
is an obvious map from this pushout square to our original pushout
square induced by the natural level equivalence $\georeal{\Sing
X}\xrightarrow{}X$.  This map is a level equivalence everywhere except
possibly the bottom right.  It follows from a general model category
argument that the induced map $\georeal{E}\xrightarrow{}D$ is a level
equivalence.  To see this, we can work one level at a time.  Let $\cd $
denote the category with three objects $i,j,k$ and two non-identity maps
$i\xrightarrow{}j$ and $i\xrightarrow{}k$.  Then the pushout defines a
functor $\spaces ^{\cd }\xrightarrow{}\spaces $ left adjoint to the diagonal
functor.  There is a model structure on $\spaces ^{\cd }$ such that a map
$f$ from $X_{j}\xleftarrow{}X_{i} \xrightarrow{}X_{k}$ to
$Y_{j}\xleftarrow{}Y_{i}\xrightarrow{}Y_{k}$ is a weak equivalence
(fibration) if and only if each map $f_{i}$, $f_{j}$, and $f_{k}$ is a
weak equivalence (fibration).  Cofibrations are defined by lifting, and
one can check that $f$ is a cofibration if and only if $f_{i}$ and the
induced maps $X_{j}\amalg _{X_{i}}Y_{i}\xrightarrow{}Y_{j}$ and
$X_{k}\amalg _{X_{i}}Y_{i}\xrightarrow{}Y_{k}$ are cofibrations.  The
diagonal is easily seen to be a right Quillen functor, so the pushout is
a left Quillen functor, and hence preserves weak equivalences between
cofibrant objects.  Hence the map $\georeal{E}_{n}\xrightarrow{}D_{n}$
is a weak equivalence as desired.

Since the level equivalence $\georeal{E}\xrightarrow{}D$ is the
composite $\georeal{E}\xrightarrow{\Sing h}\Sing
\georeal{D}\xrightarrow{}D$, it follows that $\Sing h$ is a level
equivalence.  Thus $h$ is a level equivalence as required, and the proof
is complete in the case where $r$ is a level cofibration.  

In the general case, we can factor $r=r''\circ r'$, where $r'$ is a
stable cofibration and $r''$ is a stable trivial fibration, and hence a
level equivalence.  Then we get a commutative diagram 
\[
\begin{CD}
\georeal{A} @>\georeal{f}>> \georeal{B} \\
@Vr'VV @VVs'V \\
C' @>g'>> D' \\
@Vr''VV @VVs''V \\
C @>>g> D
\end{CD}
\] 
where each square is a pushout square.  Since $r'$ is a stable
cofibration, and in particular a level cofibration, $g'$ is a stable
equivalence (and a level cofibration).  The map $r''$ is a level
equivalence.  Since the model category $\spaces $ is proper, it follows
that $s''$ is a level equivalence.  The two out of three property then
guarantees that $g$ is a stable equivalence as required.  
\end{proof}

Now we need to understand stable equivalences and colimits.

\begin{lemma}\label{lem-top-J-colimit}
Suppose $\alpha $ is an ordinal, and $X\mathcolon \alpha
\xrightarrow{}\topspec $ is a functor such that each map
$X_{i}\xrightarrow{}X_{i+1}$ is a stable cofibration and a stable
equivalence.  Then the induced map $X_{0}\xrightarrow{}\colim_{i} X_{i}$
is a stable cofibration and a stable equivalence.
\end{lemma}

\begin{proof}
It is easy to see that $X_{0}\xrightarrow{}\colim_{i} X_{i}$ is a stable
cofibration.  Consider the functor $\Sing X\mathcolon \alpha
\xrightarrow{}\spec $.  Then each map $\Sing X_{i}\xrightarrow{}\Sing
X_{i+1}$ is a level cofibration and a stable equivalence.  Let
$Y=\colim_{i} \Sing X_{i}$.  Then the map $\Sing X_{0}\xrightarrow{}Y$
is a stable equivalence, by
Lemma~\ref{lem-colimits-of-stable-equivalences}.  Thus it suffices to
show that the map $Y\xrightarrow{f}\Sing \colim_{i} X_{i}$ is a level
equivalence.

Note that $\georeal{Y}=\colim_{i} \georeal{\Sing X_{i}}$.  There is a
natural map $\georeal{Y}\xrightarrow{}\colim _{i} X_{i}$ induced by the
level equivalences $\georeal{\Sing X_{i}}\xrightarrow{}X_{i}$.  Since
homotopy groups of topological spaces commute with transfinite
compositions of cofibrations, and transfinite compositions of
cofibrations of weak Hausdorff spaces are still weak Hausdorff, it
follows that $\georeal{Y}\xrightarrow{}\colim _{i}X_{i}$ is a level
equivalence.  Since this map factors as
$\georeal{Y}\xrightarrow{\georeal{f}}\georeal{\Sing \colim
_{i}X_{i}}\xrightarrow{}\colim _{i}X_{i}$, we see that $\georeal{f}$ is
a level equivalence.  Hence $f$ is a level equivalence as required.
\end{proof}

\begin{proposition}\label{prop-top-J}
Every map in $\georeal{J}\cof $ is a stable cofibration and a stable
equivalence. 
\end{proposition}

\begin{proof}
Lemma~\ref{lem-top-J-pushout} and Lemma~\ref{lem-top-J-colimit} show
that every colimit of pushouts of maps of $\georeal{J}$ is a
stable cofibration and a stable equivalence.  Now consider an arbitrary
map $f$ in $\georeal{J}\cof $.  The proof of the Factorization
Lemma~\ref{lem-top-factorization} shows that we can factor $f=qi$, where
$i$ is a colimit of pushouts of maps of $\georeal{J}$ and $q\in
\georeal{J}\inj $.  Since $f$ has the \llp $q$, the Retract
Argument~\ref{prop-retract} shows that $f$ is a retract of $i$.  Since
$i$ is a stable cofibration and a stable equivalence, so is $f$.  
\end{proof}

We note that it is also possible to develop the homotopy theory of
topological symmetric spectra in a parallel way to our development for
simplicial symmetric spectra.  However, there are some difficulties with
this approach.  We do not know if the injective level
structure is a model structure.  Therefore, in this development, a map
$f$ would be defined to be a stable equivalence if and only if $\shom
(cf,X)$ is a weak equivalence for every $\Omega $-spectrum $X$, where
$cf$ is a cofibrant approximation to $f$ in the projective level
structure. 

\subsection{Properties of topological symmetric
spectra}\label{subsec-top-prop} 

In this section we show that the stable model structure on topological
symmetric spectra is monoidal and proper.  We also describe what we know
about monoids and modules.  Our results here are much less complete than
in the simplicial case, however.

We begin with an analogue of Theorem~\ref{thm-boxprod}.  We do not
discuss $S$-cofibrations in topological symmetric spectra.  

\begin{theorem}\label{thm-top-pushout-smash}
Let $f$ and $g$ be maps of topological symmetric spectra.  
\begin{enumerate}
\item If $f$ and $g$ are stable cofibrations, so is $f\boxprod g$. 
\item If $f$ is a stable cofibration and $g$ is a level cofibration,
then $f\boxprod g$ is a level cofibration.  
\item If $f$ is a stable cofibration and $g$ is a level cofibration, and
either $f$ or $g$ is a level equivalence, then $f\boxprod g$ is a level
equivalence.
\item If $f$ and $g$ are stable cofibrations, and either $f$ or $g$ is a
stable equivalence, then $f\boxprod g$ is a stable equivalence.  
\end{enumerate}
\end{theorem}

\begin{proof}
The class of stable cofibrations is $\georeal{FI_{\partial }}\cof $.
Thus, by Corollary~\ref{cor-check-generators}, for part~1 it suffices to
show that $\georeal{FI_{\partial }}\boxprod \georeal{FI_{\partial }}$
consists of stable cofibrations.  But $\georeal{f}\boxprod
\georeal{g}=\georeal{f\boxprod g}$, so the result follows from the
corresponding result for simplicial symmetric spectra.  

For part~2, we use Corollary~\ref{cor-check-generators} with
$I=\georeal{FI_{\partial }}$ and $J=K$ the class of level cofibrations.
Then $J\cof =J$, since $J=R(tK)\inj $, as we have already seen in the
proof of Lemma~\ref{lem-top-basic-stable}.  Any map in $I$ is of the
form $S\otimes f$, where $f$ is a level cofibration of topological
symmetric sequences.  Then, just as in the proof of part~2 of
Theorem~\ref{thm-boxprod}, we have $(S\otimes f)\boxprod g=f\boxprod g$,
where the second $\boxprod $ is taken in the category of topological
symmetric sequences.  But in $\topsymseq $, we have 
\[
(f\boxprod g)_{n}=\bigvee _{p+q=n} (\Sigma _{p+q})_{+} \sm _{\Sigma
_{p}\times \Sigma _{q}} (f_{p}\boxprod g_{q}).
\]
Each map $f_{p}\boxprod g_{q}$ is a cofibration, since $\spaces $ is a
monoidal model category, $f\boxprod g$ is a level cofibration as
required. 

For part~3, if $g$ is a level equivalence, then we can use the same
method as in part~2 to conclude that $f\boxprod g$ is a level
equivalence for $f\in \georeal{FI_{\partial }}$.  Since the class of
level trivial cofibrations is $RK\inj $ (by the analogue
of part~\ref{RK} of Lemma~\ref{lem-level-properties}), this completes
the proof in case $g$ is a level equivalence.  

To do the case when $f$ is a level equivalence, we use
Corollary~\ref{cor-check-generators} with $I=\georeal{FI_{\Lambda }}$,
$J$ the class of level cofibrations, and $K$ the class of level trivial
cofibrations.  Then $I\cof $ is the class of stable cofibrations and
level equivalences, as these are the cofibrations in the projective
level structure.  Then the proof that $I\boxprod J\subseteq K=K\cof $
proceeds just as in part~2.  

Finally, for part~4 we use Corollary~\ref{cor-check-generators} with
$I=\georeal{FI_{\partial }}$, and $J=K=\georeal{J}$, where the second
$J$ is the set of generating stable trivial cofibrations in $\spec $.
Then, as in part~1, we use the fact that the geometric realization
commutes with the box product and the corresponding fact for simplicial
symmetric spectra.  
\end{proof}

\begin{corollary}\label{cor-top-monoidal}
The stable model structure and the projective level structure on
topological symmetric spectra are monoidal.
\end{corollary}

Note that Theorem~\ref{thm-top-pushout-smash} is not as strong as
Theorem~\ref{thm-boxprod}, since Theorem~\ref{thm-boxprod} implies that
$f\boxprod g$ is a stable equivalence if $f$ is a stable cofibration,
$g$ is a level cofibration, and either $f$ or $g$ is a stable
equivalence.  We do not know if this is true for topological symmetric
spectra.  

We now show that the stable model structure is proper. 

\begin{lemma}\label{lem-top-proper}
\begin{enumerate}
\item Let 
\[
\begin{CD}
A @>g>> B \\
@VfVV @VVhV \\
X @>>> Y
\end{CD}
\]
be a pushout square of topological symmetric spectra with $g$ a level
cofibration and $f$ a stable equivalence.  Then $h$ is a stable
equivalence.  
\item Let 
\[
\begin{CD}
A @>k>> B \\
@VfVV @VVhV \\
X @>>g>Y
\end{CD}
\]
be a pullback square of topological symmetric spectra with $g$ a level
fibration and $h$ a stable equivalence.  Then $f$ is a stable
equivalence.  
\end{enumerate}
\end{lemma}

\begin{proof}
Part~2 follows immediately from the fact that $\Sing $ preserves
pullbacks, level fibrations, and stable equivalences, the corresponding
fact for simplicial symmetric spectra~\ref{lem-proper}, and the fact
that $\Sing $ detects stable equivalences.  Part~1 would also follow in
the same way if $\Sing $ preserved level cofibrations and pushouts.  The
functor $\Sing $ does preserve level cofibrations, since it preserves
monomorphisms.  But $\Sing $ only preserves pushouts up to level
equivalence.  More precisely, suppose $C$ is the pushout of $\Sing g$
and $\Sing f$ in the diagram in part~1.  Then there is a map
$C\xrightarrow{j}\Sing Y$, and we claim this map is a level equivalence.
To prove this, it suffices to prove that $\georeal{j}$ is a level
equivalence.  Since the geometric realization does preserve pushouts, it
suffices to verify that the map $\georeal{C}\xrightarrow{}Y$ induced by
the adjunction level equivalences $\georeal{\Sing D}\xrightarrow{}D$ for
$D=A,X$ and $B$ is a level equivalence.  We can work one level at a
time.  We are then reduced to showing that, if we have a diagram
\[
\begin{CD}
Z @<<< X @>f>> Y \\
@VVV @VVV @VVV \\
Z' @<<< X' @>>f'> Y'
\end{CD}
\]
of topological spaces, where $f$ and $f'$ are cofibrations and the
vertical maps are weak equivalences, the induced map $W\xrightarrow{}W'$
on the pushouts is a weak equivalence.  This statement is actually true
in any left proper model category, and in particular holds in $\top $.
A proof of it appears as~\cite[Proposition 11.3.1]{hirschhorn}. 
\end{proof}

\begin{theorem}\label{thm-top-proper}
The stable model structure on topological symmetric spectra is proper.  
\end{theorem}

\begin{proof}
This now follows immediately from Lemma~\ref{lem-top-proper}.  
\end{proof}

We now discuss monoids and modules in topological symmetric spectra.
Here we run into a serious problem: we do not know whether the monoid
axiom is true or not.  The general question of how to cope with monoids
and modules in a monoidal model category when the monoid axiom may not
hold is considered in~\cite{hovey-monoids}, where special attention is
given to the example of topological symmetric spectra.  In particular,
the following theorem is proved in~\cite{hovey-monoids}.  

\begin{theorem}\label{thm-top-modules}
Suppose $R$ is a monoid in $\topspec $.
\begin{enumerate}
\item There is a monoid $R'$ which is stably cofibrant in $\topspec $
and a level equivalence and homomorphism $f\mathcolon
R'\xrightarrow{}R$. 
\item If $R$ is stably cofibrant, there is a model structure on the
category of $R$-modules where a map is a weak equivalence or fibration
if and only if it is a stable equivalence or stable fibration in
$\topspec $.  
\item If $R$ and $R'$ are stably cofibrant in $\topspec $ and
$f\mathcolon R\xrightarrow{}R'$ is a homomorphism of monoids and a
stable equivalence, then $f$ induces a Quillen equivalence between the
corresponding module categories.
\end{enumerate}
\end{theorem}

Note that, given the good properties of closed inclusions in $\spaces $,
part~2 of Theorem~\ref{thm-top-modules} is an immediate corollary of the
fact that topological symmetric spectra form a monoidal model category.
Also, the proof of part~3 of Theorem~\ref{thm-top-modules} is quite
similar to the proof of Lemma~\ref{lem-level-quillen-equiv}. 

We interpret this theorem as asserting that, in order to study modules
over $R$, one must first replace $R$ by a stably cofibrant monoid $R'$.

Since we do not know whether the monoid axiom is true in topological
symmetric spectra, we cannot expect a model category of topological
symmetric ring spectra.  We do have the following result, also taken
from~\cite{hovey-monoids}.

\begin{theorem}\label{thm-top-monoids}
Let $\topmonoids $ denote the category of monoids in $\topspec $ and
homomorphisms.  Let $\ho \topmonoids $ denote the quotient category of
$\topmonoids $ obtained by inverting all homomorphisms which are stable
equivalences in $\topspec $.  Then the category $\ho \topmonoids $
exists and is equivalent to the homotopy category of monoids of
simplicial symmetric spectra.
\end{theorem}


\providecommand{\bysame}{\leavevmode\hbox to3em{\hrulefill}\thinspace}

\end{document}